\input  amstex

\input epsf

\baselineskip = 14pt \magnification =1100

\centerline{\bf  Rigidity of  Polyhedral Surfaces}

\bigskip

\centerline{\bf Feng Luo}

\medskip
\noindent {\bf Abstract: } {We  study  rigidity of polyhedral
surfaces and the moduli space of polyhedral surfaces using
variational principles. Curvature like quantities for polyhedral
surfaces are introduced. Many of them are shown to determine the
polyhedral metric up to isometry. The action functionals in the
variational approaches are derived from the cosine law and the
Lengendre transformation of them. These include energies used by
Colin de Verdiere, Braegger, Rivin,  Cohen-Kenyon-Propp, Leibon
and Bobenko-Springborn for variational principles on triangulated
surfaces. Our study is based on a set of identities satisfied by
the derivative of the cosine law. These identities which exhibit
similarity in all spaces of constant curvature are probably a
discrete analogous of the Bianchi identity.}

%\midspace{0.1cm}
 %\centerline{\epsfbox{1.1.eps}}
% \midspace{0.1cm}
% \centerline{Figure 1.1}

\medskip

%
%This paper is based on a lecture of the author at the workshop on discrete differential geometry at Overattach on %a
%%March 8, 2006. I was encouraged by the audience to write up the lecture notes. I'd like to thank the organizers %%
%professors Bobenko, Kenyon, Sullivan and
%Zigler for the invitation and for providing a stimulating enviroment at the conference.

\noindent {\bf Contents}

1. Introduction
$................................................................................................................$
1

2. The Derivative Cosine Law
$.........................................................................................$
9

3. Energy Functionals on the Moduli Spaces of Geometric Triangles
$............................$ 15

4. Convexity of  Moduli Spaces of Geometric Triangles
 $..................................................$ 23

5. Polyhedral Surfaces with or without Boundary
$.........................................................$ 29

6. Rigidity and Local Rigidity of Polyhedral Surfaces
 $....................................................$ 31

7. Parameterizations of the Teichmuller Space of a Surface with
Boundary $..................$ 38

8. Moduli Spaces of Polyhedral Surfaces, I, Circle Packing Metrics
 $...............................$ 43

9. Moduli Spaces of Polyhedral Surfaces, II, General Cases
$...........................................$ 46

10. Applications to Teichm\"uller Spaces and Some Open Problems
$................................$ 58

Appendix A.  A Proof of the Uniqueness of the Energy Functions
$..................................$ 63

Appendix B. The Derivative Cosine Law of Second Kind
 $...............................................$ 66

Appendix C. Relationships to the Lobachevsky Function
$...............................................$ 69

Reference
$.........................................................................................................................$
70

\medskip

\noindent \S1. {\bf Introduction}

We use variational principle to study  geometry of polyhedral
surfaces in this paper. Several
 rigidity results are established
for polyhedral surfaces and hyperbolic metrics on surfaces with
boundary.  As one consequence, for each real number $\lambda$, we
produce a natural parameterization $\psi_{\lambda}$ and a cellular
decomposition of the Teichm\"uller space of a surface with
boundary. The coordinate $\psi_0$ was constructed in [Lu2]. The
images of the Teichmuller space under these parameterizations are
convex polytopes and can be explicitly described by our work (for
$\lambda \geq 0$) and the work of Ren Guo [Gu1] (for $\lambda
<0$).

Our study is an attempt to understand the singularity formation in
a variational approach to find constant curvature metrics on
triangulated 3-manifold in [Lu3]. The work carried out in this
paper and the work of [CV1], [Ri], [Le] and others show that
singularity formation on various variational approaches on
triangulated surfaces can be well understood.

\medskip
\noindent 1.1. {\it Polyhedral surfaces their curvatures and the
main results}

A \it Euclidean (or spherical or hyperbolic) polyhedral surface
\rm is a triangulated surface with a metric, called a \it
polyhedral metric, \rm so that each triangle in the triangulation
is isometric to a Euclidean (or spherical or hyperbolic) triangle.
We emphasize that the triangulation is an intrinsic part of a
polyhedral surface. For instance, the boundary of a generic convex
polytope in the 3-dimensional space $\bold E^3$, or $\bold S^3$ or
$\bold H^3$ of constant curvature $0, 1,$ or $ -1$ is a polyhedral
surface. Two polyhedral surfaces are \it triangulation preserving
isometric \rm if there is an isometry between them preserving the
triangulations. A stronger equivalence relation which will be used
in this paper is the following. Two polyhedral metrics on the same
triangulated surface are \it equivalent (or triangulation fixing
isometric) \rm if there is an isometry between them which
preserves each simplex in the triangulation. Thus a Euclidean (or
spherical or hyperbolic) polyhedral surface is determined up to
equivalence by its \it edge length function \rm $l: \{$all
edges$\} \to \bold R_{>0}$ which assigns each edge its length. In
the sequel, we identify the equivalence class of a polyhedral
metric with its edge length function.  The \it discrete curvature
\rm of a polyhedral surface is a function assigning each vertex
$2\pi$ less the sum of the inner angles of triangles at the
vertex.

\medskip

One of the basic problems on polyhedral surfaces is to understand
the relationship between the metric and its curvature.  Since edge
lengths and inner angles of triangles in a polyhedral metric
determine the metric and its discrete curvature, we consider inner
angles of triangles as the basic unit of measurement of curvature.
Using inner angles, we introduce three families of  curvature like
quantities in this paper.  The relationships between the
polyhedral metrics and these curvature like quantities are the
main focus of the study in this paper.

\medskip
\noindent

Suppose $(S, T)$ is a closed triangulated surface so that $T$ is
the triangulation, $E$ and $V$ are the sets of all edges and
vertices. Let $\bold E^2$, $\bold S^2$ and $\bold H^2$ be the
Euclidean, the spherical and the hyperbolic 2-dimensional
geometries.

\medskip
\noindent {\bf Definition.} Given a $K^2$ polyhedral metric $l$ on
$(S, T)$ where $K^2$ $=\bold E^2$, or $\bold S^2$ or $\bold H^2$,
the \it $\phi_{\lambda}$ edge invariant \rm of the polyhedral
metric $l$ is the function $\phi_{\lambda}: E \to \bold R$ sending
an edge $e$ to:

$$ \phi_{\lambda}(e) = \int^{a}_{\pi/2} \sin^{\lambda}(t) dt
+ \int^{a'}_{\pi/2} \sin^{\lambda}(t) dt  \tag 1.1$$ where $a, a'$
are the inner angles facing the edge $e$. See figure 1.1(a).

The \it $\psi_{\lambda}$ edge invariant \rm of the metric $l$ is
the function $\psi_{\lambda}: E \to \bold R$ sending an edge $e$
to
$$ \psi_{\lambda}(e) =\int^{\frac{b+c-a}{2}}_0 \cos^{\lambda}(t)
dt + \int^{\frac{b'+c'-a'}{2}}_0 \cos^{\lambda}(t) dt \tag 1.2$$
where $b,b',c,c'$ are inner angles adjacent to the edge $e$ and
$a,a'$ are the angles facing the edge $e$. See figure 1.1(a).

\vskip.1in

\epsfxsize=5truein \centerline{\epsfbox{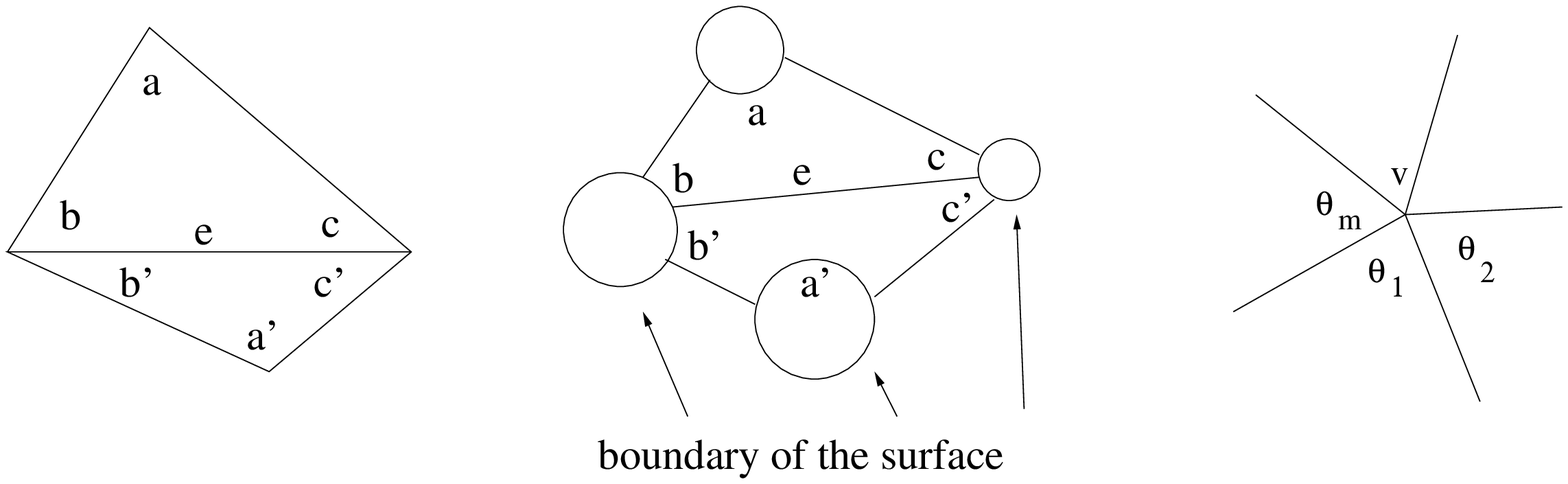}}

\medskip

\medskip

\medskip

%\newpage

%\midspace{0.1cm}
\centerline{(a)  $\quad  \quad \quad \quad \quad \quad \quad \quad
\quad  \quad \quad \quad$ (b)  $\quad \quad \quad \quad \quad
\quad \quad \quad \quad \quad \quad \quad $ (c)}
\centerline{Figure 1.1}

%latex
%\begin{figure}
%\centering
%\includegraphics{1.1.eps}
%\caption{Figure 1.1}
%\end{figure}

The \it $\lambda$-th discrete curvature $k_{\lambda}$ \rm of the
polyhedral metric $l$ on $(S,T)$ is the function $k_{\lambda} : V
\to \bold R$ sending a vertex $v$ to
$$k_{\lambda}(v) =\sum_{i=1}^m \int^{\theta_i}_{\pi/2}
\tan^{\lambda}(t/2) dt  \tag 1.3$$ where $\theta_1$, ...,
$\theta_m$ are all inner angles at vertex $v$. See figure 1.1(c).

\medskip
We remark that $\phi_0$ and $\psi_0$ edge invariants were first
introduced by I. Rivin [Ri] and G. Leibon [Le] respectively. The
0-th discrete curvature $k_0$ differs from the discrete curvature
$k$ by a constant depending on the degree $m$ of the vertex, i.e.,
$k_0(v) = -k(v) + (2-m/2)\pi$. Also the choice of the lower limit
$\pi/2$ in the integrals in (1.1), (1.3) is to make the integral
well defined for all $\lambda$. For $\lambda >-1$, the better
choice should be $0$ instead of $\pi/2$.

There are several well known rigidity theorems of Andreev-Thurston
[An], [Th], Rivin [Ri] and Leibon [Le] which relate the polyhedral
metrics with various curvatures. To state them, let us recall some
concepts. Suppose $(S, T)$ is a triangulated closed surface.
 A \it circle packing metric \rm on $(S,
T)$ is a polyhedral metric $l: E$ $\to \bold R_{>0}$ so that there
is a map, called the radius assignment,  $r: V$ $ \to \bold
R_{>0}$ with  $l(v v') = r(v)+ r(v')$ whenever the edge $vv'$ has
end points $v$ and $v'$.

\medskip
\noindent {\bf Theorem 1.1. (Thurston [Th], Rivin [Ri], Leibon
[Le])} \it Suppose $(S, T)$ is a triangulated closed surface.

(a)([Th]). A Euclidean circle packing metric on $(S, T)$ is
determined up to triangulation fixing isometry and scaling by the
0-th discrete curvature $k_0$.

(b) ([Th]). A hyperbolic circle packing metric on $(S, T)$ is
determined up to triangulation fixing isometry by the 0-th
discrete curvature $k_0$.

(c) ([Ri]). A Euclidean polyhedral metric on $(S, T)$ is
determined up to triangulation fixing isometry and scaling by the
$\phi_0$ edge invariant.

(d) ([Le]). A hyperbolic polyhedral metric on $(S, T)$ is
determined up to triangulation fixing isometry  by the $\psi_0$
edge invariant.
 \rm

\medskip
%The counterparts of Rivin's work for spherical polyhedral
%surfaces, and Leibon's work for  ideal triangulated hyperbolic
%surfaces with boundary have been established recently by [Lu1] and
%[Lu2].

%As we will see,  these  rigidity theorems of Thurston, Rivin, Leibon and those of [Lu1], [Lu2]
%are one case of a continuously family of rigidity theorems for circle packing metrics and polyhedral metrics.

One of the main results in the paper extends theorem 1.1
(a),(b),(d) to,
\medskip
\noindent
{\bf Theorem 1.2.} \it Let $ \lambda \in \bold R $ and $(S, T)$ be  a closed triangulated surface.

(a)  A Euclidean circle packing metric on $(S, T)$ is determined
up to triangulation fixing isometry and scaling by its
$k_{\lambda}$-th discrete curvature.

(b) A hyperbolic circle packing metric on $(S, T)$ is determined
up to triangulation fixing isometry by its $k_{\lambda}$-th
discrete curvature.

(c) If $\lambda \leq -1$,  a Euclidean polyhedral metric on $(S,
T)$ is determined up to triangulation
 fixing isometry and scaling
 by its  $\phi_{\lambda}$ edge invariant.

(d) If $\lambda \leq -1$ or $\lambda \geq 0$, a spherical
polyhedral metric on $(S, T)$ is determined up to triangulation
fixing isometry by its $\phi_{\lambda}$ edge invariant.

(e) If $\lambda \leq -1$ or $\lambda \geq 0$, a hyperbolic
polyhedral surface is determined up to triangulation fixing
isometry by its $\psi_{\lambda}$ edge invariant.

 \rm
\medskip
 For any $\lambda \in \bold R$,  there are
 local rigidity theorems in cases (c), (d), (e) (see theorem 6.2).
%or $\lambda=0$, the above locally rigidity result is a consequence of Rivin's theorem.
%The cases of $\lambda=0, -1$ deserve special attention. In the case of $\lambda=0$, Rivin [Ri] proved that the set of all edge invariants $\phi_0$'s of Euclidean polyhedral metrics on $(S, T)$
%so that  $\phi_0(e) \leq \pi$ for all edges $e$ forms a convex polytope.
Whether these results hold for $\lambda$'s not listed in theorems
1.1 and 1.2 is not clear to us. It deserves a further study. To
the best of our knowledge, theorem 1.2 for the simplest case of
the boundary of a tetrahedron is new.

The counterpart of theorem 1.2(e) for hyperbolic metrics with
totally geodesic boundary on an ideal  triangulated compact
surface is the following. Recall that an \it ideal triangulated
compact surface \rm with boundary $(S, T)$ is obtained by removing
a small open regular neighborhood of the vertices of a
triangulation of a closed surface. The \it edges \rm  of an ideal
triangulation $T$ correspond bijectively to the edges of the
triangulation of the closed surface. Given a hyperbolic metric $l$
with geodesic boundary on an ideal  triangulated surface $(S, T)$,
there is a unique geometric ideal triangulation $T^*$ isotopic to
$T$ so that all edges are geodesics orthogonal to the boundary.
The edges in $T^*$ decompose the surface into hyperbolic
right-angled hexagons. The \it $\psi_{\lambda}$ edge invariant \rm
of the hyperbolic metric $l$ is defined to be the map
$\Psi_{\lambda}: \{$ all edges in $T$\} $\to \bold R$ sending each
edge $e$ to

$$\psi_{\lambda}(e) = \int_0^{\frac{b+c-a}{2}} \cosh^{\lambda}(t) dt
+\int^{\frac{b'+c'-a'}{2}}_0 \cosh^{\lambda}(t) dt \tag 1.4$$
where $a,a'$ are lengths of arcs in the boundary (in the ideal
triangulation $T^*$) facing the edge and $b,b',c,c',$ are the
lengths of arcs in the boundary adjacent to the edge so that
$a,b,c$ lie in a hexagon. See figure 1.1(b).

\medskip
\noindent {\bf Theorem 1.3.} \it  A hyperbolic metric with totally
geodesic boundary on an ideal  triangulated compact surface
 is determined up to triangulation fixing isometry by
its $\psi_{\lambda}$-edge invariant.  Furthermore, if $\lambda
\geq 0$, then the set of all $\psi_{\lambda}$-edge invariants on a
fixed ideal  triangulated surface is an explicit open convex
polytope $P_{\lambda}$ in  a Euclidean space so that
$P_{\lambda}=P_0$. \rm

\medskip

The case when $\lambda <0$ has been recently established by Ren
Guo [Gu1]. He proved that,

\medskip
 \noindent {\bf Theorem 1.4.} (Guo) \it Under the same
assumption as in theorem 1.3, if $\lambda <0$, the set of all
$\psi_{\lambda}$-edge invariants on a fixed ideal  triangulated
surface is an explicit bounded open convex polytope $P_{\lambda}$
in  a Euclidean space. Furthermore, if $\lambda < \mu$, then
$P_{\lambda} \subset P_{\mu}$. \rm

\medskip

 Theorem 1.3 was proved for
$\lambda=0$ in [Lu2] where the open convex polytope $P_0$ was
explicitly described. Evidently for each $\lambda \in \bold R$,
the edge invariant $\psi_{\lambda}$ can be taken as a coordinate
of the Teichm\"uller space of the surface. The interesting part of
the theorem 1.3 is that the images of the Teichm\"uller space in
these coordinates (for $\lambda \geq 0$) are all the same. Whether
these coordinates are related to quantum Teichm\"uller theory is
an interesting topic. See [CF],[Ka],[BL], [Te] for more
information. Combining theorem 1.3 with the work of Ushijima [Us]
and Kojima [Ko], one obtains for each $\lambda \geq 0$ a cell
decomposition of the Teichm\"uller space invariant under the
action of the mapping class group. See corollary 10.6.

Similar results for the moduli spaces of all $\phi_{\lambda}$ or
$\psi_{\lambda}$ edge invariants, or $k_{\lambda}$ discrete
curvatures on a triangulated surfaces are also obtained in this
paper. In Thurston's notes [Th], he showed that the spaces of all
0-th discrete curvatures of Euclidean or hyperbolic circle packing
metrics on a triangulated surface are convex polytopes. One of the
results in this paper states,

\medskip
\noindent {\bf Theorem 8.1.} \it Suppose $\lambda \leq -1$ and
$(S, T)$ is a closed triangulated surface so that $E$ is the set
of all edges. Then,

(a) The space of all $k_{\lambda}$-discrete curvatures of
Euclidean circle packing metrics on $(S, T)$ forms a proper
codimension-1 smooth submanifold in $\bold R^E$.

(b) The space of all $k_{\lambda}$-discrete curvatures of
hyperbolic circle packing metrics on $(S, T)$ is an open
submanifold in $\bold R^E$ bounded by the proper codimension-1
submanifold in part (a).

\rm

\medskip

Theorems 9.1 and theorem 9.8 give descriptions of the moduli
spaces of $\phi_{\lambda}$ and $\psi_{\lambda}$ edge invariants
for $\lambda \leq -1$ for the Euclidean, hyperbolic and spherical
polyhedral surfaces.

%Results similar to theorem 1.2  for
% hyperbolic polyhedral surfaces,
%Teichm\"uller space of surfaces with boundary
%and circle packing type metrics on surfaces are also established in this paper.
%These results generalize the corresponding theorems of
%Thurston, Leibon and [Lu2] to a continuous families of rigidity theorems.

\medskip
\noindent 1.3. {\it The method of proofs and the related works}
\medskip

The proofs of above theorems use variational principles. The use
of variational principle on triangulated surfaces in recent time
appeared in the seminal work of Colin de Verdiere [CV1] in 1991.
We are highly influenced by the works of [CV1], [Th], [Ri], [MaR]
and [Le]. Beside the work of [CV1] and [Ri], variational
principles on triangulated surfaces have also appeared in Braegger
[Br], Leibon [Le], Cohen-Kenyon-Propp [CKP], Bobenko-Springborn
[BS], Springborn [Sp], Schlenkar [Sh], [Lu1], [Lu2] and others.
The energy functions used in [CV1], [Ri], [Br], [Le], [CKP]  are
related to the 3-dimensional volume or its Legendre transform.
Even in the work of [CV1], Colin de Verdiere's energy was
motivated by the 3-dimensional Schlaefli volume formula [CV2].
Very recently in 2004, motivated by the discrete 2-dimensional
integrable system, Bobenko and Springborn [BS] discovered a new
collection of energies for triangulated surfaces.

We observe that all energy functions used by Colin de Verdiere,
Braegger, Rivin, Cohen-Kenyon-Propp, Leibon, Bobenko-Springborn
and Cohen-Kenyon-Propp can be constructed using the cosine law and
the Legendre transformation. Furthermore, we show that these
known energy functions are special cases of continuous families of
energy functions derived from the cosine law.   We also show that
these families are the complete lists of all localized energies
one could construct. All rigidity results of the paper are
consequences of those convex or concave energy functions.

%This probably is not surprising since volume of tetrahedra is
%governed by the Schlaefli formula. To be more precise, volume is
%the integration of the Schlaefli 1-form $\frac{1}{2}\sum_{i, j}
%l_{ij}dx_{ij}$ where $x_{ij}$ and $l_{ij}$ are the dihedral angle
%and the signed length at the ij-th edge. The edge length $l_{ij}$
%is obtained from the dihedral angles $x_{rs}$'s using the cosine
%law twice. In [GL], we observed that the energy function of
%Bobenko-Springborn can also be derived from the cosine law.

\medskip
\noindent
1.4. {\it The derivative cosine law}

Our study is motivated by discretization of 2-dimensional
Riemannian geometry. In the discrete setting, the smooth
 metric is replaced by the polyhedral metric  and the Gaussian curvature is replaced by the discrete curvature. From this point
of view, the relationship between a polyhedral metric and its
curvature is essentially the cosine law for triangles. Thus, the
\it cosine law should be considered as a metric-curvature
relation. \rm  Just like in Riemannian geometry, it is natural to
study the infinitesimal dependence of curvature (inner angles) on
the metric (edge lengths). The result is a collection of
identities which we call the \it derivative cosine law\rm. Among
the most interesting ones are the following. Suppose a triangle in
$\bold S^2$ or $\bold E^2$ or $\bold H^2$, has inner angles
$\theta_1, \theta_2$, $\theta_3$ and opposite edge lengths $l_1,
l_2, l_3$. Then for $i,j,k$ pairwise distinct indices, the
following identities (obtained in [CL]) hold in all geometries
$\bold S^2$, $\bold E^2$ or $\bold H^2$,

$$ \frac{\partial \theta_i /\partial l_j}{\partial \theta_j/\partial l_i} =\frac{\sin \theta_i}{\sin
\theta_j} \quad \quad  \text{and}   \quad \quad  \frac{ \partial
\theta_i/\partial l_j}{\partial \theta_i/\partial l_i} = -\cos
\theta_k.  \tag 1.5$$

\medskip
\noindent
1.5. {\it Two examples}

 We would like to illustrate the use of these identities by two examples.
These examples also show the main techniques and methods used in the paper.
In the first example, given a Euclidean triangle of edge lengths $l_1, l_2, l_3$ and opposite angles $\theta_1, \theta_2, \theta_3$,
the cosine law relating them states,
$$ \cos \theta_i = \frac{ l_j^2 + l_k^2 -l_i^2}{2 l_j l_k}, $$
where $i \neq j \neq k \neq i$.
 Consider $\theta_i
=\theta_i(l_1, l_2, l_3)$ as a smooth function of $l=(l_1, l_2,
l_3)$.  Then identity (1.5) shows that the differential 1-form,
$$\omega =\sum_{i=1}^3 \ln \tan (\theta_i/2) d l_i $$
is closed. This closed smooth 1-form is defined on the space of
all Euclidean triangles (parameterized by the edge lengths) $\bold
E^2(l,3) =\{(l_1, l_2, l_3) | l_i + l_j > l_k\}$. Since the space
$\bold E^2(l,3)$ is convex, the integration $F(l)
=\int^l_{(1,1,1)} \omega$ defines a smooth function on the $\bold
E^2(l,3)$. By definition, this function $F$ satisfies
$$ \frac{\partial F}{ \partial l_i} = \ln \tan(\theta_i/2).  \tag 1.6$$
Its Hessian matrix
$ [\frac{\partial^2 F}{\partial l_r \partial l_s}]_{3 \times 3}
 =[ \frac{1}{\sin \theta_r} \frac{\partial \theta_r}{\partial l_s}]_{3 \times 3}$
can be shown to be congruent to the Gram matrix of the triangle.
Thus the Hessian matrix is  semi-positive definite. It follows
that the function $F$ is convex on $\bold E^2(l, 3)$.  Property
(1.6) says the variation at the i-th edge length (i.e., the
metric) of the function $F$ depends only on the opposite angle
$\theta_i$ (i.e., the curvature). Due to (1.6), we call $F(l_1,
l_2, l_3)$ the $F$-energy of the triangle $(l_1, l_2, l_3)$. A
function with property (1.6) is very useful for variational
framework on polyhedral surfaces. Let  a Euclidean polyhedral
surface $(S, T, l)$  be given so that  $E$ is the set of all edges
and $l: E \to \bold R$ is the edge length function. Define an
``energy" $W(l)$ of the metric $l$ to be the sum of the
$F$-energies of its triangles, i.e., $W(l) = \sum_{ \{a,b,c\} \in
T^{(2)}} F(l_a, l_b, l_c)$. Then the function $W(l)$ is convex in
$l$ since it is the summation of convex functions.  Furthermore,
by property (1.6), we have
$$ \frac{ \partial W(l)}{\partial l_i} = \ln (\tan(a/2))  + \ln (\tan(a'/2))  \tag 1.7$$
where $a$ and $a'$ are the two inner angles facing the i-th edge
$e_i$.  Identities (1.1) and (1.7) show that the gradient of the
convex function $W$ is
$$\bigtriangledown W =(\phi_{-1}(e_1), ...,   \phi_{-1}( e_{|E|}   ) ).$$

On the other hand, it is well known that, if $U: \Omega \to \bold
R$ is a smooth strictly convex function defined on an open convex
set $\Omega$ in $\bold R^n$, then the gradient map
$\bigtriangledown U :\Omega \to \bold R^n$ is a smooth embedding.
The function $W$ is not strictly convex.  With a little extra
work,  we prove that the gradient $\bigtriangledown W$ determines
the metric $l$ up to scaling. This is theorem 1.2(c) for
$\lambda=-1$. Indeed, theorem 1.2(c) is proved in exactly the same
way by using a special collection of closed 1-forms on the space
of all Euclidean triangles.

\medskip
The next example is the seminal work of Colin de Verdiere [CV1] in
1991 who produced an energy function for circle packing metrics on
surfaces. For simplicity, we deal with the case of hyperbolic
metrics. Similar energies for the spherical and Euclidean
triangles were also introduced in [CV1]. Given a hyperbolic
triangle of edge lengths $l_i = r_j + r_k$ and inner angles
$\theta_i$ so that $\theta_i$ is facing the $l_i$-th edge,  Colin
de Verdiere proved that the differential 1-form $\eta_0 =
\sum_{i=1}^3 \frac{\theta_i}{\sinh(r_i)} d r_i =\sum_{i=1}^3
\theta_i du_i$ is closed on the space of hyperbolic triangles $\{
(r_1, r_2, r_3) \in \bold R_{>0}^3\}$ where $u_i=
\int^{r_i}_{\infty} \frac{1}{ \sinh(t)} dt$. Furthermore, its
integration $F(u_1, u_2, u_3) =\int^u \eta_0$ is a strictly
concave function  in $(u_1, u_2, u_3) \in \{ (u_1, u_2, u_3) \in
\bold R^3 | u_i < 0\}$.  By the construction,  one has
$\frac{\partial F}{\partial u_i} = \theta_i$, i.e., the variation
of the energy at the i-th radius depends only on the inner angle
at the i-th vertex. Using the strictly convexity,
 Colin de Verdiere gave a new proof of theorem 1.1(a) and  1.1(b). He
also proved  Thurston's existence of circle packing theorem (in
the tangential case) using the variational principle associated to
the energy $\int^u \eta_0$.

\medskip
\noindent To summarize, in the variational approach on a
triangulated surface $(S, T)$, for each polyhedral metric  on $(S,
T)$, we define the ``energy" of the polyhedral metric to be the
sum of the ``energies" of its (metric) triangles. Thus the key
step is to find all possible ``energy" functions of the triangles.
Depending on geometric problems, a natural condition that one
imposes on the energy $F(l_1, l_2, l_3)$ of a triangle with edge
lengths $l=(l_1, l_2, l_3)$ is that $\frac{\partial F(l)}{\partial
g(l_i)} =f(\theta_i)$ for some non-constant functions $f, g$ for
all indices $i$. Namely, the  variation of the energy with respect
to a scaled metric $(g(l_1), g(l_2), g(l_3))$ is a function of the
curvature (i.e., (1.6)). This problem is equivalent to find all
closed 1-forms of the form
 $w=\sum_{i=1}^3 f(\theta_i) dg(l_i)$ in the space of all geometric triangles
($F$ is then the integration of the 1-form $w$).  This and the
related problem for circle packing radii parameterization are
solved in the paper.

\medskip
\noindent
1.6. {\it  The main technical result }
\medskip

The most general form of the cosine law can be stated as follows.
Suppose a function $y=y(x)$ where $y=(y_1, y_2, y_3) \in \bold
C^3$  and $x=(x_1, x_2, x_3)$ is in some open connected set in
$\bold C^3$ so that $x_i$'s and $y_i$'s are related by
$$ \cos ( y_i) = \frac{\cos x_i + \cos x_j \cos x_k}{\sin (x_j) \sin (x_k)} \tag 1.8$$
where $\{i,j,k\}=\{1,2,3\}$. We say $y=y(x)$ is the \it  cosine
law function\rm.
 Let $r_i =\frac{1}{2}(x_j+x_k-x_i)
$. Then $r=(r_1, r_2, r_3)$ is a new parametrization so that $x_i = r_j + r_k$.

\medskip
\noindent {\bf Theorem 1.5.} \it For the cosine law function
$y=y(x)$, all closed 1-forms of the form $ w =\sum_{i=1}^3
f(y_i)dg(x_i)$ where $f,g$ are two non-constant smooth functions,
are up to scaling and complex conjugation,
$$\omega_{\lambda} =\sum_{i=1}^3 \int^{y_i} \sin^{\lambda}(t) dt  d( \int^{x_i} \sin^{-\lambda -1}(t) dt)
=\sum_{i=1}^3 \frac{\int^{y_i} \sin^{\lambda}(t)
dt}{\sin^{\lambda+1}(x_i) } dx_i$$ for some $\lambda \in \bold C$.
In particular, all closed 1-forms are  holomorphic or
anti-holomorphic.

All closed 1-forms of the form $ \sum_{i=1}^3 f(y_i) dg(r_i)$ where $f,g$ are two non-constant smooth functions, are up to scaling and complex conjugation,
$$\eta_{\lambda} = \sum_{i=1}^3 \int^{y_i} \tan^{\lambda}(t/2) dt  d(\int^{r_i} \cos^{-\lambda-1}(t) dt)
=\sum_{i=1}^3 \frac{\int^{y_i} \tan^{\lambda}(t/2)
dt}{\cos^{\lambda+1}(r_i)} dr_i$$ for some $\lambda \in \bold C$.
In particular, all closed 1-forms are  holomorphic or
anti-holomorphic. \rm

\medskip
 By
specializing theorem 1.5 to various cases of $\bold S^2$, $\bold
E^2$ and $\bold H^2$ and integrating the 1-forms, we obtain
various energy functionals for variational framework on
triangulated surfaces. See theorems 3.2 and 3.4. We have also
identified all those convex or concave energies constructed in
this way. Finally we remark that the most interesting closed
1-forms may have already been discovered by various authors. For
instance, $\eta_0$ was  discovered by Colin de Verdiere [CV1] for
$\bold S^2$, $\bold E^2$ and $\bold H^2$ triangles, $\eta_{-1}$
for hyperbolic triangle was in Leibon's work [Le] and $\eta_{-1}$
for hyperbolic right-angled hexagon was in [Lu2], and
$\omega_{-1}$ for spherical triangle was discovered in [Lu1]. For
Euclidean triangles, the form $\omega_0=\sum_{i=1}^3
\frac{\theta_i}{ l_i} dl_i$ was discovered by Cohen, Kenyon and
Propp [CKP], and its Legendre transformation $\sum_{i=1}^3 \ln l_i
d\theta_i$ was in the work of Rivin [Ri].  The 1-parameter family
of closed 1-forms containing the 1-form in Bobenko-Springborn's
work is in corollary B2 in the appendix B.

\medskip
\noindent 1.7. {\it Surfaces with boundary}

\medskip
The results obtained in this paper can be generalized without
difficulty to compact triangulated surfaces with boundary. Given a
compact triangulated surface $(S, T)$ with boundary, by doubling
across its boundary, we obtain a closed triangulated surface. The
notions of $k_{\lambda}$, $\phi_{\lambda}$ and $\psi_{\lambda}$
invariants can now be defined for polyhedral metrics on $(S, T)$
by using the corresponding concepts on the closed surface. For
simplicity, we will not state the results for triangulated
surfaces with boundary.

\medskip
\noindent 1.8. {\it The organization of the paper} \medskip

The paper is organized as follows.  In section 2, we list some of
the properties of the derivatives of the cosine law and prove
theorem 1.5. In section 3, we deduce various consequences of
theorem 1.5 in the geometric content. In section 4, we study the
convexity of the space of all geometric triangles in various
parameters arising from the energy functions. In \S5, we set up
the framework for variational problems on triangulated surfaces.
\S6 is devoted to prove theorem 1.2. \S 7 is devoted to prove
theorem 1.3 and studies the shapes of the Teichm\"uller space in
the coordinates. In \S8, we study the space of all $\lambda$-th
discrete curvatures of circle packing metrics. The spaces of all
$\phi_{\lambda}$ and $\psi_{\lambda}$ edge invariants are
investigated in \S9. In \S 10, we discuss applications to the
Teichm\"uller space and some open problems. In the appendices, we
give a proof of the uniqueness of the energy functions, derive the
derivative cosine law of the second kind, and recall some of the
known relationships of some energy functions with the Lobachevsky
functions.

\medskip
\noindent 1.9. We would like to thank David Gu and Ren Guo for
discussions related to the topics in the paper. Part of the work
was carried out when we participated the Oberwolfach workshop on
discrete differential geometry. We would like to thank the
organizers, A. Bobenko, R. Kenyon,  and J. Sullivan for the
invitation. The work is supported in part by the NSF.

\medskip
\noindent \S2.  {\bf The Derivative Cosine Law}
\medskip
 A smooth function defined in an open set in $\bold R^n$
is called locally convex (or locally strictly convex) if its
Hessian matrix is semi-positive definite (or positive definite) at
all points.  Let $\{i,j,k\}=\{1,2,3\}$ in this section.

\medskip
\noindent
2.1. {\it The derivative cosine law.}
\medskip
\noindent Given a triangle in $\bold H^2$, $\bold E^2$ or $\bold
S^2$ of inner angles $\theta_1, \theta_2, \theta_3$ and edge
lengths $l_1, l_2, l_3$ so that $\theta_i$ is facing the $l_i$-th
edge, the cosine law expressing length $l_i$ in terms of the
angles $\theta_r$'s is,
$$ \cos(\sqrt{\lambda}l_i) = \frac{ \cos \theta_i + \cos \theta_j \cos \theta_k}{\sin \theta_j \sin \theta_k} \tag 2.1$$
where $\lambda=1,-1,0$ is the curvature of the space $\bold S^2$,
or $\bold H^2$ or $\bold E^2$.
  Another related cosine law
is$$ \cosh(l_i) = \frac{ \cosh \theta_i + \cosh \theta_j \cosh \theta_k}{\sinh \theta_j \sinh \theta_k} \tag 2.2$$
for right-angle hyperbolic hexagon with three non-adjacent edge lengths $l_1, l_2, l_3$ and their opposite edge lengths
$\theta_1, \theta_2, \theta_3$.

Identities (2.1) and (2.2) show that the cosine laws are
specialization of the cosine law function $y=y(x)$ given in (1.8).
The following is a simple derivative calculation. See [Lu3]  for a
proof.

\medskip
\noindent {\bf Theorem 2.1.} \it Suppose the cosine law function
$y=y(x)$ is defined on an open connected set in $\bold C^3$ which
contains a point $(a,a,a)$ so that $y(a,a,a) = (b,b,b)$. Let
$A_{ijk} = \sin y_i \sin x_j \sin x_k$ where
$\{i,j,k\}=\{1,2,3\}$. Then
 $$ \quad \quad \quad  A_{ijk} = A_{jki}.  \tag 2.3 $$
$$A^2_{ijk} =  1- \cos^2 x_i -\cos^2 x_j -\cos^2 x_k - 2 \cos x_i \cos x_j \cos x_k.  \tag 2.4$$
At a point $x$ where $A_{ijk} \neq 0$, then,
$$ \frac{ \partial y_i}{\partial x_i }= \frac{\sin x_i}{A_{ijk}},    \tag 2.5$$
$$ \frac{  \partial y_i}{\partial x_j} = \frac{\partial y_i}{\partial x_i} \cos y_k, \tag 2.6$$
$$  \cos(x_i) = \frac{ \cos y_i - \cos y_j \cos y_k}{\sin y_j \sin y_k}. \tag 2.7$$
\rm

\medskip
\noindent
2.2. {\it Remarks}
\medskip
\noindent {\bf 1.} Formula (2.3) shows that $\frac{\sin y_i}{\sin
x_i}$ is independent of the index $i$. We call it the \it sine
law\rm.
% From now on, we use $A_{123}=\sin(y_1)\sin(x_2)\sin(x_3)$ and $B_{123}=\sin(x_1)\sin(y_2)\sin(y_3)$ for the cosine law function $y=y(x)$. \rm

\noindent
{\bf 2.}  Identity (2.7) can be written in the symmetric form as,
$$\cos(\pi-x_i) = \frac{ \cos(\pi- y_i)+  \cos(\pi- y_j) \cos (\pi- y_k)}{\sin(\pi- y_j) \sin (\pi-y_k)}. \tag 2.8$$
\rm This is the reflection of the duality of the spherical
triangles. Namely, the dual triangle of a spherical triangle has
edge lengths $\pi-\theta_i$ and inner angles $\pi-l_i$. In
particular by (2.6) applied to (2.8), we obtain
 $$\frac{\partial x_i}{\partial y_j}  = -\frac{\partial x_i}{\partial y_i} \cos x_k. \tag 2.9$$

\medskip
\noindent {\bf 3.} Identity (1.5) follows from (2.6) and (2.9).

\medskip
\noindent {\bf 4.} If we consider $(y_i, x_j, x_k)$ as a function
of $(y_j, y_k, x_i)$ in the cosine law, there are similar
derivative identities which we call the derivative cosine laws of
second kind. See appendix B. The energy functions  of
Bobenko-Springborn can be derived from them.

\medskip
\noindent 2.3. {\it The tangent law and the radius
parametrization}
\medskip

In many geometric considerations, we encounter situations (for
instance circle packing or Delaunay triangulations) where the
natural parametrization of the triangle whose i-th edge length (or
angle) $x_i$ is given by $r_j+r_k$, i.e., one uses $(r_1, r_2,
r_3)$ to parameterize $(x_1, x_2, x_3)$ where $r_i
=\frac{1}{2}(x_j+x_k -x_i)$. If $x_i$'s are the edge lengths, then
$r_i$'s are the radii of the pairwise tangent circles whose
centers are the vertices of the triangle. If $x_i$'s are the inner
angles, then $r_i$ is $\pi/2$ less the angle between the i-th edge
and the circumcircle.

\bigskip
\medskip
\vskip.1in

\epsfxsize=3.5truein \centerline{\epsfbox{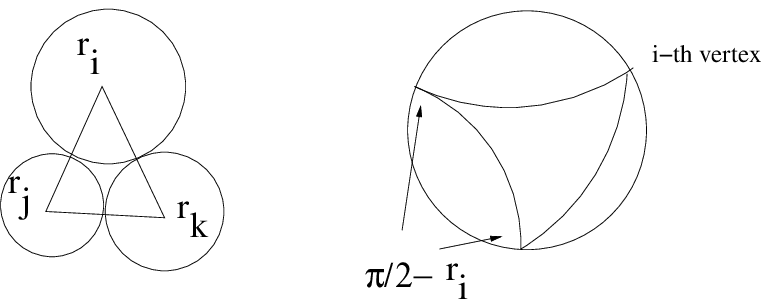}}

\medskip

\medskip

\medskip

%\newpage

%\midspace{0.1cm}
\centerline{(a)  $\quad  \quad\quad \quad \quad \quad  \quad \quad
\quad  \quad \quad \quad$ (b)} \centerline{Figure 2.1}
\medskip
\noindent
{\bf Lemma  2.2}. \it For the cosine law function $y=y(x)$,
 write $x_i = r_j + r_k$, or $r_i =\frac{1}{2}(x_j + x_k -x_i)$, then the following expression is independent of the
indices, $$ \frac{\tan^2(y_i/2)}{\cos^2(r_i)}
=-\frac{\cos(r_1+r_2+r_3)}{\cos(r_1)\cos(r_2) \cos(r_3)}. \tag
2.10$$ Furthermore, there is a quantity
$A=\sin(y_i)\sin(x_j)\sin(x_k)$ independent of indices so that
$$\frac{1}{\cos(r_i)} \frac{\partial y_i}{\partial r_j}= \frac{2 \cos(r_k)}{A \sin(x_k)}      \tag 2.11          $$
and
$$ \frac{\partial y_i}{\partial r_i} =   \frac{2 \sin(x_j+x_k)\cos(r_j)\cos(r_k)}
{A \sin  (x_j) \sin(x_k)}$$

%=\frac{4
%\tan(y_i/2)\cos^2(y_i/2)(\cos^2(y_j/2)+\cos^2(y_k/2)-1)}{B}    $$

In particular,
$$\frac{\partial y_i/\partial r_j}{\partial y_j /\partial r_i} =
\frac{\cos(r_i)}{\cos(r_j)}=\frac{\tan(y_i/2)}{\tan(y_j/2)}. \tag 2.12$$

\rm
\medskip

The identity (2.10) will be used many times in the paper and will
be called the \it tangent law. \rm  In the case of a Euclidean
triangle of edge lengths $r_i+r_j$ and opposite angle $\theta_k$,
identity (2.10) degenerates to state that $r_i \tan(\theta_i/2)$
is independent of the indices $i$ (in fact, the common value is
the radius of the inscribed circle).

\bigskip
\noindent {\bf Proof.} To see (2.10), let us calculate
$\tan^2(y_i/2)$. It is
$$\tan^2(y_i/2) = \frac{ 1-\cos(y_i)}{1+\cos(y_i)}$$
$$=\frac{ \sin(x_j)\sin(x_k)
-\cos(x_i) - \cos(x_j)\cos(x_k)}{ \sin(x_j)\sin(x_k) + \cos(x_i) +
\cos(x_j)\cos(x_k)}$$
$$=-\frac{ \cos(x_i)+ \cos(x_j+x_k)}{\cos(x_i) + \cos(x_j-x_k)}$$
$$= -\frac{ \cos(r_i)\cos(r_1+r_2+r_3)}{\cos(r_j) \cos(r_k)}$$
$$= -\frac{\cos^2(r_i)\cos(r_1+r_2+r_3)}{ \cos(r_1)\cos(r_2)\cos(r_3)}.$$
Thus (2.10) follows.

 Next, let us calculate $\frac{\partial
y_i}{\partial r_j}$ for $i \neq j$. Note that due to
$x_i=r_j+r_k$, we have $\frac{\partial}{\partial r_j}
=\frac{\partial}{\partial x_i }+ \frac{\partial}{\partial x_k}$.
Thus
$$\frac{\partial y_i}{\partial r_j }= \frac{\partial y_i}{\partial x_i} +\frac{\partial y_i}{\partial x_k}$$
$$=\frac{\partial y_i}{\partial x_i}( 1+ \cos(y_j))$$
$$=\frac{\sin(x_i)}{A} \frac{ \cos(x_j)+\cos(x_i-x_k)}{\sin(x_i)\sin(x_k)}$$
$$=\frac{2 \cos(r_i) \cos(r_k)}{A \sin(x_k)}.$$
%On the other hand,
%$$\frac{\partial y_i}{\partial x_i}( 1+ \cos(y_j))=\frac{\sin (y_i)}{B} 2 \cos^2(y_j/2)=\frac{4 \tan(y_i/2)\cos^2(y_i/2)\cos^2(y_j/2)}{B}.$$
This establishes (2.11).

 To see $\frac{\partial y_i}{\partial
r_i}$, we have
$$\frac{\partial y_i}{\partial r_i} = \frac{ \partial y_i}{\partial x_j}+\frac{\partial y_i}{\partial x_k}$$
$$=\frac{\partial y_i}{\partial x_i}( \cos(y_k)+\cos(y_j))$$
$$=\frac{\sin(x_i)}{A}(\frac{\cos(x_k)+\cos(x_i)\cos(x_j)}{\sin(x_i)\sin(x_j)}+  \frac{\cos(x_j)+\cos(x_i)\cos(x_k) }{\sin(x_i)\sin(x_k)})$$
$$=\frac{ \sin(x_j)\cos(x_j)+\sin(x_k)\cos(x_k) + \cos(x_i)(\cos(x_j)\sin(x_k)+\cos(x_k)\sin(x_j))}{A\sin(x_j)\sin(x_k)}$$
$$=\frac{ \sin(x_j+x_k) \cos(x_j-x_k) + \cos(x_i)\sin(x_j+x_k)}{A\sin(x_j)\sin(x_k)}$$
$$=\frac{2 \sin(x_j+x_k) \cos(r_j)\cos(r_k)}{A \sin(x_j) \sin(x_k)}.$$

%This verifies (2.12).

% Finally, we could also express $\frac{\partial y_i}{\partial
%r_i} =  \frac{\partial y_i}{\partial x_j}+ \frac{\partial
%y_i}{\partial x_k} =\frac{\partial y_i}{\partial x_i}(
%\cos(y_k)+\cos(y_j))$ as $\frac{\sin(y_i) (\cos(y_j) +
%\cos(y_k))}{B}$ = $\frac{ 4\tan(y_i/2)
%\cos^2(y_i/2)(\cos^2(y_j/2)+\cos^2(y_k/2)-1)}{B}$. QED
\bigskip

\noindent 2.4. {\it A proof of theorem 1.5}
\medskip
\noindent {\bf Theorem 1.5.} \it For the cosine law function
$y=y(x)$, all closed 1-forms of the form $ w =\sum_{i=1}^3
f(y_i)dg(x_i)$ where $f,g$ are two non-constant smooth functions,
are up to scaling and complex conjugation,
$$\omega_{\lambda} =\sum_{i=1}^3 \int^{y_i} \sin^{\lambda}(t) dt  d( \int^{x_i} \sin^{-\lambda -1}(t) dt)
=\sum_{i=1}^3 \frac{\int^{y_i} \sin^{\lambda}(t) dt}{\sin^{\lambda+1}(x_i)} dx_i$$
for some $\lambda \in \bold C$.

All closed 1-forms of the form $ \sum_{i=1}^3 f(y_i) dg(r_i)$
where $f,g$ are two non-constant smooth functions and
$r_i=\frac{x_j+x_k-x_i}{2}$, are up to scaling and complex
conjugation,
$$\eta_{\lambda} = \sum_{i=1}^3 \int^{y_i} \tan^{\lambda}(t/2) dt  d(\int^{r_i} \cos^{-\lambda-1}(t) dt)
=\sum_{i=1}^3 \frac{ \int^{y_i} \tan^{\lambda}(t/2) dt}{
\cos^{\lambda+1}(r_i)} dr_i$$ for some $\lambda \in \bold C$.

In particular, all closed 1-forms are  holomorphic or
anti-holomorphic. \rm

%%%stopped checking on August 11, 2006

\medskip
\noindent {\bf Proof.} First, we prove that the 1-forms in the
list are closed.
  Indeed, a holomorphic 1-form
 $\omega_{\lambda}=\sum_{i=1}^ 3 f(y_i) g'(x_i) dx_i$ is
closed if and only if $\frac{ \partial (f(y_i) g'(x_i))}{\partial
x_j}=f'(y_i)g'(x_i) \frac{\partial y_i}{\partial x_j}$ is
symmetric in $i,j$. But for $\omega_{\lambda}$, it is equal to
$$f'(y_i)g'(x_i) \frac{\partial y_i}{\partial x_j} =
\frac{1}{A}f'(y_i)g'(x_i) \sin(x_i)\cos(y_k) =
\frac{1}{A}(\frac{\sin(y_i)}{\sin(x_i)})^{\lambda}\cos(y_k)$$
where $A$ is independent of the indices by theorem 2.1.  The above
expression is clearly symmetric in $i,j$ due to the sine law. Thus
the closeness follows.

To see the holomorphic 1-form $\eta_{\lambda}=\sum_{i=1}^3 f(y_i)
d g(r_i)$ is closed, we check if the quantity $f'(y_i) g'(r_i)
\frac{\partial y_i}{\partial r_j}$ is symmetric in $i,j$. It is
equal to $$\frac{\tan^{\lambda}(y_i/2)}{\cos^{\lambda+1}(r_i)}
\frac{\partial y_i}{\partial r_j} =
(\frac{\tan(y_i/2)}{\cos(r_i)})^{\lambda} (\frac{1}{\cos(r_i)}
\frac{\partial y_i}{\partial r_j}) $$
 where
$ \frac{1}{\cos(r_i)} \frac{\partial y_i}{\partial r_j} $ is
symmetric in $i,j$ by (2.11) and $\frac{\tan(y_i/2)}{\cos(r_i)}$
is independent of $i$ by the tangent law.  This shows the
closeness.

The proof that these are  all closed 1-forms is relatively long
and will
 not be used anywhere in the paper.
We  defer it to the appendix A.  The key part of the proof is the
following result concerning the uniqueness of the sine law and the
tangent law (see appendix A for a proof).

\medskip
\noindent {\bf Lemma 2.3.} \it Suppose $y=y(x)$ is the cosine law
function and $f, g$ are two smooth non-constant functions.

(a) If $f(y_i)/g(x_i)$ is independent of the indices for all $x$, then there are
constants $\lambda$, $\mu$, $c_1, c_2$ so that,  $f(t) =c_1\sin^{\lambda}(t)  \sin^{\mu}(\bar t)$
and $g(t) =c_2\sin^{\lambda}(t)\sin^{\mu}(\bar t)$.

(b) If $r_i = \frac{1}{2}(x_j+x_k-x_i)$ where
$\{i,j,k\}=\{1,2,3\}$, and $f(y_i)/g(r_i)$ is independent of the
indices for all $(r_1, r_2, r_3)$, then there are constants
$\lambda, \mu$, $c_1, c_2$ so that $f(t) = c_1\tan^{\lambda}(t)
\tan^{\mu}(\bar t)$ and $g(t) =
c_2\cos^{\lambda}(t)\cos^{\mu}(\bar t)$.

\rm \medskip \noindent 2.5. {\it Remarks on the Hessian matrices
of the energy functions}
\medskip
In theorem 1.5, let $u_i=\int^{x_i} \sin^{-\lambda-1}(t) dt$ in
the case of $w_{\lambda}$ and $u_i=\int^{r_i} \cos^{-\lambda-1}(t)
dt$ for $\eta_{\lambda}$, then  $ w_{\lambda}=\sum_{i=1}^3
\int^{y_i} \sin^{\lambda}(t) dt du_i$ and $\eta_{\lambda}
=\sum_{i=1}^3 \int^{y_i} \tan^{\lambda}(t/2) dt du_i$.
 The Hessian matrices $H=[
\partial ^2F/\partial u_s \partial u_t]$ of the integral $F(u)
=\int^u w_{\lambda}$ (or $F(u)=\int^u \eta_{\lambda}$) has the
following properties: for any two $\lambda$ and $\lambda'$, the
associated Hessian matrices of $  \int^u w_{\lambda}$ and $\int^u
w_{\lambda'}$ are congruent. The same also holds for $\int^u
\eta_{\lambda}$ and $\int^u \eta_{\lambda'}$. Indeed, in the case
of $w_{\lambda}$, we have
$$\frac{\partial^2 F}{\partial u_s \partial u_t}
=\sin^{\lambda}(y_t)\frac{\partial y_t}{\partial x_s}
\frac{\partial x_s}{\partial u_s}$$
$$=\sin^{\lambda}(y_t)\sin^{\lambda+1}(x_s) \frac{\partial
y_t}{\partial x_s}$$
$$=(\frac{\sin(y_t)}{\sin(x_t)})^{\lambda}
(\sin(x_t)\sin(x_s))^{\lambda} (\sin(x_s)\frac{\partial
y_s}{\partial x_t}).  \tag 2.13$$ Let
$q=(\frac{\sin(y_t)}{\sin(x_t)})^{\lambda}$ and $D$ be the $3
\times 3$ diagonal matrix whose $(i,i)$-th entry is
$\sin^{\lambda}( x_i)$. Then (2.13) shows that the Hessian matrix
is $q D [\sin(x_s) \frac{\partial y_s}{\partial x_t}] D$. It
follows that the Hessian matrices for different $\lambda$'s are
congruent.

By the same calculation using the tangent law instead of the sine
law, for the integration of $\eta_{\lambda}$, the (s,t)-th entry
of the Hessian matrix is
$$ =(\frac{\tan(y_t/2)}{\cos(r_t)})^{\lambda} (\cos(r_s)
\cos(r_t))^{\lambda} (\cos(r_t)\frac{\partial y_t}{\partial r_s}).
\tag 2.14$$

It again shows that the Hessian matrices for different $\lambda$'s
are congruent.

\medskip

\medskip
\noindent 2.6. {\it The Legendre transformation.}
\medskip
 The integrals $\int  w_{\lambda}$ in theorem 1.5 are not independent. In fact $\int w_{\lambda}$ and
$\int w_{-\lambda -1}$ are essentially Legendre transformation of
each other.

Let us recall briefly the Legendre transformation. Suppose $U$ and
$ V$ are diffeomorphic connected open sets in $\bold R^n$ so that
the first de Rham cohomology group $H^1_{dR}(U)=0$. Let $x=(x_1,
..., x_n) \in U$ and $y=(y_1, ..., y_n) \in V$ and $y=y(x): U \to
V$ be a diffeomorphism so that its Jacobian matrix is symmetric,
i.e.,
$$ \frac{\partial y_i}{\partial x_j} = \frac{ \partial y_j}{\partial x_i}.$$
Then the differential 1-forms $w_{U} = \sum_{i=1}^n y_i dx_i$ and
$w_V =\sum_{i=1}^n x_i dy_i$ are closed in $U$ and $V$
respectively. We call  $w_{U}$  the \it Legendre transformation
\rm of $w_V$. Their integrations $f(x) = \int^x_a w_{U}$ and $g(y)
=\int_b^y w_V + (a,b)$ where $b=y(a)$ and $(a,b)$ is the dot
product are well defined due to $H^1_{dR}(U)=H^1_{dR}(V) =0$. We
call $g(y)$ the \it Legendre transformation \rm of $f(x)$ (and
vice versa).  It is well known that the Legendre transform of a
strictly convex (or concave) function is  strictly convex (or
concave).

\medskip
\noindent {\bf Proposition 2.4.} \it Let $f_{\lambda} (x) =
\int^x_{(\pi/2, \pi/2, \pi/2)} \sum_{i=1}^3
\frac{\int^{y_i}_{\pi/2} \sin^{\lambda}(t)
dt}{\sin^{\lambda+1}(x_i)} dx_i$. Then the Legendre transformation
of $f_{\lambda}(x)$ is $f_{-\lambda-1}(x) -f_{-\lambda-1}(0)$. In
particular, the Legendre transform of $f_{-1/2}(x)$ is $f_{-1/2}(x
)$ up to adding a constant.\rm

\medskip
\medskip
\noindent
{\bf Proof.} Let $g_{\lambda}(y)$ be the Legendre transformation of $f_{\lambda}(x)$. Then,
$$ g_{\lambda}(y) =\int^y_{(\pi/2,\pi/2,\pi/2)} \sum_{i=1}^3 \int^{x_i}_{\pi/2} \sin^{-\lambda-1}(t) dt  d(\int^{y_i}_{\pi/2} \sin^{\lambda}(t) dt) + c$$
for some constant $c =g_{\lambda}(0,0,0)-g_{\lambda}(\pi/2,\pi/2,\pi/2)$.
 In the above identities, $x$ and $y$ are related by the cosine law (1.8) which we will denote by $y=CL(x)$. Let
$w=(w_1, w_2, w_3)$ and $v=(v_1,v_2, v_3)$ so that $w_i=\pi-y_i$
and $v_i=\pi-x_i$. Then by (2.8), we have $v=CL(w)$. Making a
change of variables $x_i=\pi-v_i$ and $y_i=\pi-w_i$ in the
integral of $g_{\lambda}$, we obtain
%$$g_{\lambda}(y) =
%\int^y_{(\pi/2,\pi/2,\pi/2)} \sum_{i=1}^3 \int^{x_i}_{\pi/2} \sin^{-\lambda-1}(t) dt  d(\int^{y_i}_{\pi/2} \sin^{\lambda}(t) dt)$$
$$g_{\lambda}(y)=\int^{(\pi-w_1, \pi-w_2, \pi-w_3)}_{(\pi/2,\pi/2,\pi/2)} \sum_{i=1}^3 \int^{\pi-x_i}_{\pi/2} \sin^{-\lambda-1}
(\pi-t) d(\pi-t) \sin^{\lambda}(\pi-w_i) d(\pi-w_i)+c$$
$$=\int^{(\pi-w_1, \pi-w_2, \pi-w_3)}_{(\pi/2, \pi/2, \pi/2)} \sum_{i=1}^3 \int^{v_i}_{\pi/2} \sin^{-\lambda-1}(t)dt \sin^{\lambda}(w_i) dw_i+c$$
$$=f_{-\lambda-1}(y)+c$$
since $v=CL(w)$. Thus $g_{\lambda}(y) = f_{-\lambda-1}(y)
-f_{-\lambda-1}(0)$.

\medskip

\medskip

\medskip

\medskip
%\newline

%noindent
\newpage
\noindent
 \S3. {\bf Energy Functionals on the Moduli Spaces of
Geometric Triangles}
\medskip
We will restrict the 1-forms in theorem 1.5 to the moduli spaces
of geometric triangles and determine the convexity of the
integration of the 1-forms on the moduli spaces.

Suppose a triangle in $\bold S^2$, $\bold E^2$ or $\bold H^2$ has
inner angles $\theta_1, \theta_2,$ $\theta_3$ and opposite edge
lengths $l_1, l_2, l_3$.

\medskip
\noindent
3.1. {\it Derivative cosine laws for geometric triangles}
\medskip
Take $(x_1,x_2, x_3)$ and $(y_1, y_2, y_3)$  in theorem 2.1 to be
the inner angles and edge lengths. Then,

\medskip
\noindent {\bf Corollary 3.1.} \it Let $\{i,j,k\}=\{1,2,3\}$.
There is a positive quantity $A$ independent
 of indices so that

(a) ([CL])  $$\frac{\partial \theta_i}{\partial l_j} =
-\frac{\partial \theta_i}{\partial l_i }\cos(\theta_k),$$ and $$
\frac{\partial \theta_i}{\partial l_i} =
\frac{\sin(\theta_i)}{A}>0.$$

(b) For spherical triangles,
$$ \frac{\partial l_i}{\partial \theta_j} =\frac{\partial l_i}{\partial \theta_i} \cos(l_k),$$
and $$\frac{\partial l_i}{\partial \theta_i} = \frac{\sin(l_i)}{A}
>0.$$

(c) For  hyperbolic triangles,
$$\frac{\partial l_i}{\partial \theta_j} =\frac{ \partial l_i}{\partial \theta_i }\cosh(l_k), $$
and
$$\frac{\partial l_i}{\partial \theta_i }=-\frac{\sinh(l_i)}{A} < 0.$$

(d) For a hyperbolic right-angled hexagon of three non-pairwise adjacent edge lengths $l_1, l_2, l_3$ and opposite
edge lengths $\theta_1, \theta_2, \theta_3$,
$$\frac{\partial \theta_i}{\partial l_j} = -\frac{\partial \theta_i}{\partial l_i} \cosh(\theta_k),$$
and $$\frac{ \partial \theta_i}{\partial l_i} =\frac{\sinh(\theta_i)}{A} > 0.$$ \rm

\medskip
The proof is routine using theorem 2.1 by taking care of the
curvature factor $\lambda=\pm 1, 0$ appeared in (2.1).  Note that
$\cos(\sqrt{-1} x) =\cosh(x)$, $\sin(\sqrt{-1} x) = \sqrt{-1}
\sinh(x)$ and $\sinh(\sqrt{-1} x) =\sqrt{-1} \sin(x)$. Using these
relations,  part (a) follows from (2.9) where $x_i=\theta_i$ and
$y_i = \sqrt{\lambda}l_i$. Part (a)
 for  Euclidean triangle was established in [CL] and can be checked easily.  Parts (b) and (c) follow from theorem 2.1.
To see part (d), note that the cosine law for hexagon can be written as,
$$ \cos(\pi-\sqrt{-1} \theta_i) = \frac{ \cos(\sqrt{-1} l_i) +\cos(\sqrt{-1} l_j) \cos(\sqrt{-1} l_k)}{\sin(\sqrt{-1} l_j) \sin(\sqrt{-1} l_k)}.$$
Thus part (d) follows from theorem 2.1.

\medskip
\noindent 3.2. {\it Closed 1-forms on moduli spaces of geometric
triangles, I}
\medskip
Using corollary 3.1 and theorem 1.5, we obtain,
\medskip
\noindent {\bf Theorem 3.2.}  \it  Let $l=(l_1, l_2, l_3)$ and
$\theta=(\theta_1, \theta_2, \theta_3)$ be lengths and angles of a
triangle in $\bold E^2$, $\bold H^2$ or $\bold S^2$. The following
is the complete list, up to scaling, of all closed real-valued
1-forms of the form $\sum_{i=1}^3 f(\theta_i) d g(l_i)$
 for some non-constant smooth functions $f, g$. Let $\lambda \in
 \bold R$ and $u=(u_1, u_2, u_3)$.

\noindent (a) For a Euclidean triangle,

$$ w_{\lambda} = \sum_{i=1}^3 \frac{\int^{\theta_i} \sin^{\lambda}(t) dt}{l_{i}^{\lambda+1}}  dl_i.$$
Furthermore, its integration $\int^u w_{\lambda}$ is locally
convex in variable $u$ where $u_i =\int^{l_i}_1 t^{-\lambda -1}
dt$.

\noindent (b) For a spherical triangle,
$$ w_{\lambda} = \sum_{i=1}^3 \frac{\int^{\theta_i} \sin^{\lambda}(t) dt}{\sin^{\lambda+1}(l_i)} dl_i.$$
The integral $\int^u w_{\lambda}$ is locally strictly convex in
$u$ where
 $u_i =\int^{l_i}_{\pi/2}$  $ \sin^{-\lambda-1}(t) dt$.

\noindent (c) For a hyperbolic triangle,
$$ w_{\lambda} = \sum_{i=1}^3 \frac{\int^{\theta_i} \sin^{\lambda}(t) dt}{\sinh^{\lambda+1}(l_i)} dl_i.$$

\noindent (d) For a hyperbolic right-angled hexagon,
$$ w_{\lambda} = \sum_{i=1}^3 \frac{\int^{\theta_i} \sinh^{\lambda}(t) dt}{\sinh^{\lambda+1}(l_i)} dl_i.$$

In the cases of (b), (c), (d), by taking the Legendre
transformation, we also obtain the complete list of all closed
1-forms of the form $\sum_{i=1}^3 g(l_i) df(\theta_i)$.

\rm

\medskip

 \noindent {\bf Proof.} The closeness of these 1-forms
is evident due to theorem 1.5 and corollary 3.1 except in the case
of Euclidean triangle. In the case of $\bold E^2$, we need to
verify that the expression
$$ \frac{\partial}{\partial l_j}  (\frac{\int^{\theta_i} \sin^{\lambda}(t)dt }{l_i^{\lambda+1}}) \tag 3.1$$ is symmetric in $i,j$.
By corollary 3.1 that $\frac{\partial \theta_i}{\partial l_j}
=-\frac{\partial \theta_i}{\partial l_i} \cos (\theta_k) =
-\frac{\sin(\theta_i)}{A} \cos(\theta_k)$, the expression (3.1)
is equal to $\frac{1}{l_i^{\lambda +1}} \sin^{\lambda}(\theta_i)
 \frac{\partial \theta_i}{\partial l_j}$
$=-(\frac{\sin(\theta_i)}{\l_i})^{\lambda+1}\frac{\cos
(\theta_k)}{A}$ where $A$ is independent of indices. It is
symmetric in $i,j$ due to the sine law.

To verify the convexity, note that if $u_i =g(l_i)$ and $w
=\sum_{i=1}^3 f(\theta_i) du_i$ is closed, then the Hessian of the
function $F(u) =\int^u w$ is $[\frac{\partial^2 F}{\partial u_r
\partial u_s}] =[\frac{f'(\theta_r)}{g'(l_s)}\frac{ \partial
\theta_r} {\partial l_s}]$.

In the cases (a)-(d), by the sine law and the choice of $f,g$, $f'(\theta_i)\sin (\theta_i)g'(l_i) =q$
is a positive function independent of the indices. Thus the $(r,s)$-th entry of the Hessian matrix can be written as,
$$\frac{f'(\theta_r)}{g'(l_s)}\frac{\partial \theta_r}{\partial l_s} =(f'(\theta_r)\sin(\theta_r)g'(l_r))(\frac{1}{g'(l_r)g'(l_s)})(\frac{1}{\sin(\theta_r)}\frac{\partial \theta_r}{\partial l_s})$$
$$=q (\frac{1}{g'(l_r)g'(l_s)})(\frac{1}{\sin(\theta_r)}\frac{\partial \theta_r}{\partial l_s}).$$

This shows that the Hessian matrix can be written as
$q D L  D$ where $D$ is the positive diagonal matrix whose (i,i)-th entry is $\frac{1}{g(l_i)}$ and $L =[\frac{1}{\sin (\theta_r)} \frac{\partial \theta_r}{\partial l_s}]$.
 Recall that given a triangle with inner angles $\theta_1, \theta_2, \theta_3$, its
(angle) Gram matrix  $[a_{rs}]_{3 \times 3}$ satisfies  $a_{ii}=1$
and $a_{ij} = -\cos (\theta_k)$ ($\{i,j,k\}=\{1,2,3\}$). On the
other hand, by corollary 3.1 (a), the matrix $L$ is equal to the
Gram matrix multiplied by the positive function $1/A$. As a
consequence,  the Hessian of the integral of the 1-forms in
(a)-(d) is congruent to the Gram matrix of the triangle. It is
well known that the Gram matrix of a Euclidean triangle is
semi-positive definite of rank 2 and the Gram matrix of a
spherical triangle is positive definite. (See for instance [Lu3]
for proofs). Thus the local convexity of the integrations for
Euclidean and spherical triangles follows. QED

%%%8-12-2006

\medskip
\noindent {\bf Corollary 3.3.}  \it In theorem 3.2(a), the null
space of the Hessian matrix $Hess(F)$ of $F=\int^u w_{\lambda}$ at
a point $u$ is generated by the vector $u$ if $\lambda \neq 0$ and
is generated by $(1,1,1)$ if $\lambda =0$. \rm

\medskip
Indeed, if $w=(w_i, w_j, w_k)$ is in the kernel of the Hessian,
then by definition and the calculation above, we have
$$ \sum_{r,s \in \{i,j,k\}} \frac{ w_r w_s}{g'(l_r) g'(l_s)}
a_{rs} =0$$ where $[a_{rs}]_{3 \times 3}$ is the Gram matrix of
the Euclidean triangle. It is well known that the null space of
the Gram matrix of a Euclidean triangle of edge lengths $(l_i,
l_j, l_k)$ is generated by the length vector $(l_i, l_j, l_k)$. It
follows that there is a constant $c \in \bold R$ so that $w_r =c
l_r g'(l_r) = cl_r^{-\lambda}$ for $r=i,j,k$. Now if $\lambda \neq
0$, then $u_r =-\frac{1}{\lambda} l_r^{-\lambda}$. Therefore, if
$\lambda \neq 0$, $w = -\frac{c}{\lambda} u$ and if $\lambda =0$,
$w=c(1,1,1)$.
\medskip
\noindent
3.3. {\it Closed 1-forms on the moduli spaces of geometric triangles, II}
\medskip
Let $\{i,j,k\}=\{1,2,3\}$ in this subsection. The next result is
the counterpart of theorem 3.2 for triangles parameterized by the
radii. There are two cases to be discussed: (1) the edge lengths
are $l_k=r_i+r_j$ and opposite angles are $\theta_k$ and (2) edge
lengths are $l_i$ and the opposite angles $\theta_i= r_j+r_k$.

We will use the following well known fact from linear algebra.
Given a symmetric matrix $M=[m_{st}]$ so that $m_{ii} \geq \sum_{j
\neq i} |m_{ij}|$ for all indices $i$, then $M$ is semi-positive
definite. If $m_{ii} > \sum_{j \neq i} |m_{ij}|$ for all $i$, then
the matrix $M$ is positive definite. See for instance [CL] lemma
3.10 for a proof. We call $M$ with this property a \it diagonally
dominate \rm matrix.

\medskip

\medskip
\noindent {\bf Theorem 3.4.} \it The following are the complete
list, up to scaling, of all closed real-valued 1-forms of the form
$\sum_{i=1}^3 f(l_i) d g(r_i)$ (where $\theta_i=r_j+r_k$) and
$\sum_{i=1}^3 f(\theta_i) d g(r_i)$ (where $l_i=r_j+r_k$) for some
non-constant smooth functions $f,g$.  Let $\lambda \in \bold R$
and $u=(u_1, u_2, u_3)$.

(a) For a Euclidean triangle of angles $\theta_i$ and opposite
edge lengths $r_j+r_k$,
$$\eta_{\lambda} =\sum_{i=1}^3 \frac{\int^{\theta_i} \cot^{\lambda}(t/2) dt}{ r_i^{\lambda+1}}
dr_i.$$  Its integration $\int^{u} \eta_{\lambda}$ is locally
concave in $u=(u_1, u_2, u_3)$ where $u_i =\int_1^{r_i}
t^{-\lambda-1} dt$.

(b) For a hyperbolic triangle of angles $\theta_i$ and opposite edge lengths $r_j+r_k$,
$$\eta_{\lambda} = \sum_{i=1}^3 \frac{\int^{\theta_i} \cot^{\lambda }(t/2) dt }{\sinh^{\lambda+1}(r_i)} dr_i.$$
 Its integration $\int^{u}
\eta_{\lambda}$ is locally strictly concave in $u$ where $u_i
=\int_1^{r_i} \sinh^{\lambda}(t) dt$.

(c) For a spherical triangle of angles $\theta_i$ and opposite edge lengths $r_j+r_k$,
$$\eta_{\lambda} = \sum_{i=1}^3 \frac{\int^{\theta_i} \cot^{\lambda}(t/2) dt }{\sin^{\lambda+1}(r_i)} dr_i.$$

(d) For a hyperbolic triangle of edge lengths $l_i$ and opposite angles $r_j+r_k$,
$$\eta_{\lambda} = \sum_{i=1}^3 \frac{\int^{l_i} \tanh^{\lambda}(t/2) dt}{ \cos^{\lambda+1}(r_i)} dr_i.$$
 Its integration $\int^{u}
\eta_{\lambda}$ is locally strictly convex in $u$ where $u_i
=\int_1^{r_i} \cos^{\lambda-1}(t) dt$.

(e) For a spherical triangle of edge lengths $l_i$ and opposite angles $r_j+r_k$,
$$\eta_{\lambda} = \sum_{i=1}^3 \frac{\int^{l_i} \tan^{\lambda}(t/2) dt}{ \cos^{\lambda+1}(r_i)} dr_i.$$

(f) For
 a hyperbolic right-angled hexagon of three non-pairwise adjacent edge lengths $l_i$ and opposite edge lengths $r_j+r_k$,
$$\eta_{\lambda} = \sum_{i=1}^3 \frac{\int^{l_i} \coth^{\lambda}(t/2) dt }{\cosh^{\lambda+1}(r_i)} dr_i.$$
 Its integration $\int^{u}
\eta_{\lambda}$ is locally strictly concave in $u$ where $u_i
=\int_1^{r_i} \cosh^{\lambda-1}(t) dt$.\rm

\medskip

\noindent {\bf Proof.} The proof of the uniqueness is essentially
the same as that of theorem 1.5 and will be omitted (see appendix
A).
  The proof of closeness of
 the 1-forms is just a specialization of theorem 1.5 by taking care of the curvature factors.
For (b) and (c), we take $y_i =\pi-\theta_i$ and $x_i=\pi
-\sqrt{\delta} l_i$ in theorem 1.5 where $\delta =\pm 1$ is the
curvature of the space $\bold S^2$
 or $\bold H^2$. For (d) and (e), we take $y_i=\sqrt{\delta} l_i$
and $x_i=\theta_i$. For (f), we take $y_i =\pi - \sqrt{-1} l_i$
and $x_i=\sqrt{-1} \theta_i$ in theorem 1.5. The closeness of the
1-forms together with the convexity in part (a) will be proved
 below.

A short proof of the convexity or concavity of the functions in
cases (a), (b), (d), (f) goes as follows. By remark 2.5 and
(2.13), (2.14), it follows that for any two $\lambda$ and
$\lambda'$, the Hessian matrices of the associated functions
$\int^u \eta_{\lambda}$ and $\int^u \eta_{\lambda'}$ in each case
of (a)-(f) are congruent. Thus, to check the convexity or
concavity in cases (a), (b), (d) and (f), it suffices to verify it
for a specific value of $\lambda$. This has been done by various
authors. For cases (a), (b) and $\lambda =0$, Colin de Verdiere
[CV1] proved the concavity of the function $\int^u \eta_0$. In the
case (d), Leibon [Le] proved the strictly convexity for $\lambda
=0$. In the case (f), we proved it for $\lambda=0$ in [Lu2].

\medskip

Below is a more detailed argument producing proofs of the
convexity or concavity. Moreover, we obtain a complete description
(corollary 3.5) of the null space of the hessian matrix for the
energy in the case (a).
\medskip
\noindent Let $u_i=g(r_i) $. Then the Hessian of the function
$F(u)=\int^u \eta_{\lambda}$ is the  matrix $H=[H_{st}]_{3 \times
3}= [ \frac{ \partial ^2 F(u)}{\partial u_t \partial u_s}]$ where
$ H_{st} =\frac{f'(z_s)}{g'(r_t)} \frac{ \partial z_s}{\partial
r_t}$.

\medskip

In the case (a), we first show that $\eta_{\lambda} =\sum_{i=1}^3
(\int^{\theta_i} \cot^{\lambda}(t/2) dt) r_i^{-\lambda-1} dr_i$ is
closed.  By the definition of the radius $R$ of the inscribed
circle, we have   $\frac{\coth(\theta_i/2)}{r_i}=\frac{1}{R}$.
 Note that by corollary 3.1(a),
$$\frac{\partial \theta_i}{\partial r_j} = \frac{\partial \theta_i}{\partial l_i} +\frac{\partial \theta_i}{\partial l_k}$$
$$=\frac{\partial \theta_i}{\partial l_i} (1 -\cos(\theta_j))$$
$$=\frac{2\sin (\theta_i) \sin^2(\theta_j/2)}{A}$$
$$=\frac{4  \cot(\theta_i/2) \sin^2(\theta_i/2) \sin^2(\theta_j/2)}{A}$$
$$=\frac{ 4 r_i \sin^2(\theta_i/2) \sin^2(\theta_j/2)}{AR}. \tag 3.2$$
By (3.2), for the 1-form $\sum_{i} f(\theta_i) d g(r_i)$ where
$f'(\theta) = \cot^{\lambda}(\theta/2)$ and $g'(r) =
r^{-\lambda-1}$, we have
$$\frac{\partial (f(\theta_i) g'(r_i))}{\partial r_j} = f'(\theta_i) g'(r_i)\frac{ \partial \theta_i}{\partial r_j}$$
$$=\cot^{\lambda}(\theta_i/2) r_i^{-\lambda-1} ( \frac{4 r_i}{AR} \sin^2(\theta_i/2) \sin^2(\theta_j/2))$$
$$= \frac{4}{AR}(r_i \tan(\theta_i/2))^{-\lambda} \sin^2(\theta_i/2) \sin^2(\theta_j/2)$$
$$= \frac{4R^{-\lambda-1}}{A} \sin^2(\theta_i/2) \sin^2(\theta_j/2)$$
is symmetric in $i,j$, i.e., the 1-form $\eta_{\lambda}=\sum_{i} f(\theta_i) dg(r_i)$ in case (a) is closed.

The Hessian $H=[H_{st}]$ of the integration $F(u) =\int^u \eta_{\lambda}$ is

$$H_{st} =\frac{f'(\theta_s)}{g'(r_t)} \frac{ \partial \theta_s}{\partial r_t}$$
$$=(r_s \tan(\theta_s/2))^{-\lambda} (r_s r_t)^{\lambda} (r_t  \frac{\partial \theta_s}{\partial r_t})$$
$$=R^{-\lambda} (r_s r_t)^{\lambda} a_{st}.$$
Since $r_s \tan(\theta_s/2) =R>0$ is independent of the indices,
the above identity
 shows that the Hessian matrix $H$ is congruent to $[a_{st}]$ where, by (3.2),
$$ a_{ij} =r_j  \frac{\partial \theta_i}{\partial r_j} = \frac{4 r_i r_j}{AR} \sin^2(\theta_i/2) \sin^2(\theta_j/2) >0  \tag 3.3$$
and
$$ a_{ii} = r_i  \frac{\partial \theta_i}{\partial r_i} =r_i(
\frac{ \partial \theta_i}{\partial l_j} + \frac{\partial \theta_i}{\partial l_k})$$
$$=- r_i \frac{\partial \theta_i}{\partial l_i}( \cos(\theta_k) + \cos(\theta_j)) <0  \tag 3.4$$
due to $\theta_k + \theta_j  < \pi$ and $\frac{\partial
\theta_i}{\partial l_i}>0$. Furthermore, we have $|a_{ii}| -|
a_{ji}| - |a_{ki}| = -(a_{ii} + a_{ji} + a_{ki}) =- r_i
\frac{\partial (\theta_1+\theta_2 + \theta_3)}{\partial r_i} =0$.
Thus the matrix $[a_{rs}]$ is semi-negative definite due to the
diagonal dominate property above. Thus the Hessian matrix $H$ is
semi-negative definite. This shows that the integration is locally
concave in $u$.

 \medskip
 \noindent
 {\bf Corollary 3.5.} \it In theorem 3.4(a), for a Euclidean triangle of edge lengths
 $l_i = r_j + r_k$, the null space of the Hessian Hess(F) where $F(u) =\int^u \eta_{\lambda}$ at a
 point $u$ is generated by the vector $u$ for $\lambda \neq 0$ and
 by $(1,1,1)$ if $\lambda =0$. \rm

 \medskip

 Indeed, we first note that the null space of the symmetric $3 \times 3$ matrix
 $[a_{st}]$ contains $(1,1,1)$ due to the equality
 $a_{ii}+a_{ij}+a_{ik}=0$. Next, we note that the rank of
 $[a_{st}]$ is 2, due to the fact that $|a_{ii}| > |a_{ij}|$, for
 all $i \neq j$. Thus the null space of $[a_{st}]$ is generated by
 $(1,1,1)$. Now the Hessian matrix $H$ is $R^{-\lambda} [
 r_s^{\lambda} r_{t}^{\lambda} a_{st}]$ by the calculation above.
 It follows that the null space of $H$ is generated by $(1,1,1)$
 if $\lambda =0$ or by $\frac{1}{\lambda} (r_1^{-\lambda}, r_2^{-\lambda},
 r_3^{-\lambda}) =u$ if $\lambda \neq 0$.

\medskip
In case (f), we use lemma 2.2 for $y_i=\pi- \sqrt{-1} l_i$ and
$x_i = \sqrt{-1} \theta_i$. Thus the tangent law in lemma 2.2
shows $ \coth(l_i/2)/ \cosh(r_i)$ is independent of the indices
and $\frac{1}{\cosh(r_i)}\frac{\partial l_i}{\partial r_j}$ is
symmetric in $i,j$. Now the Hessian $H=[H_{st}]$ of the function
$F(u) =\int^u \eta_{\lambda}$ where $f'(l)=\coth^{\lambda}(l/2)$
and $g'(r) = \cosh^{-\lambda-1}(r)$,  is,
$$ H_{st} =\frac{ f'(l_s)}{g'(r_t)} \frac{ \partial l_s}{\partial r_t}$$
$$= \coth^{\lambda}(l_s/2)\cosh^{\lambda+1}(r_t) \frac{ \partial l_s}{\partial r_t}$$
$$= (\frac{\coth(l_s/2)}{ \cosh(r_s)})^{\lambda}(\cosh(r_t)\cosh(r_s))^{\lambda+1} (\frac{1}{\cosh(r_s) } \frac{\partial l_s}{\partial r_t})$$
$$= (\frac{\coth(l_s/2)}{ \cosh(r_s)})^{\lambda}(\cosh(r_t)\cosh(r_s))^{\lambda+1} (\frac{1}{\cosh(r_s) } a_{st}).$$

By the tangent law in lemma 2.2, $\coth(l_s/2)/\cosh(r_s) >0$ is
independent of the indices,  the above identity shows that the
Hessian $H$ is congruent to the matrix $[a_{st}]$ where
$$a_{ij} =\frac{1}{\cosh(r_i)} \frac{ \partial l_i}{\partial r_j}$$
$$ =\frac{1}{\cosh(r_i)} (\frac{ \partial l_i}{\partial \theta_i}+\frac{\partial l_i}{\partial \theta_k})$$
$$= -\frac{1}{\cosh(r_i)} \frac{ \partial l_i}{\partial \theta_i} ( \cosh(l_j) -1).$$

We have
$$ a_{ii} = \frac{1}{\cosh(r_i)} \frac{\partial l_i}{\partial r_i}$$$$ = -\frac{1}{\cosh(r_i)}
\frac{ \partial l_i}{\partial \theta_i} ( \cosh(l_k) +
\cosh(l_j)).$$ By corollary 3.1(d) $\frac{\partial l_i}{\partial
\theta_i}
>0$, this shows that $a_{st} <0$ for all $s, t$ and $|a_{ii}| >
|a_{ij} | + |a_{ik}|$ for all $i$. Thus the matrix $-[a_{st}]$ is
diagonally dominated.  It follows that the Hessian matrix $H$ is
negative definite. Thus the integration $\int^u \eta_{\lambda}$ is
locally strictly concave.

\medskip
In case (b), we take $y_i =\pi-\theta_i$ and $x_i = \pi-\sqrt{-1}
l_i$ in lemma 2.2.  By lemma 2.2 $\cot(\theta_i/2)/ \sinh(r_i)$ is
independent of indices and $\frac{1}{\sinh(r_i)} \frac{ \partial
\theta_i}{\partial r_j}$ is symmetric in $i,j$. The Hessian
$H=[H_{st}]$ of the integral $F(u) = \int^u \eta_{\lambda}$ where
$f'(\theta) =\cot^{\lambda}(\theta/2)$ and $g'(r) =
\sinh^{-\lambda-1}(r)$ is
$$H_{st} =\frac{f'(\theta_s)}{g'(r_t)} \frac{ \partial \theta_s}{\partial r_t}$$
$$
 =\frac{\cot^{\lambda}(\theta_s/2)}{\sinh^{-\lambda-1}(r_t)} \frac{ \partial \theta_s}
{\partial r_t}$$
$$=(\frac{\cot(\theta_s/2)}{ \sinh(r_s)})^{\lambda}(\sinh(r_s) \sinh(r_t))^{\lambda} (\sinh(r_t)  \frac{ \partial \theta_s}{\partial r_t})$$
$$= (\frac{\cot(\theta_s/2)}{ \sinh(r_s)})^{\lambda}(\sinh(r_s) \sinh(r_t))^{\lambda} (\sinh(r_t) a_{st}).$$

By the tangent law $\cot(\theta_s/2)/\sinh(r_s) >0$ is independent
of the indices $s$,
 the Hessian matrix $H$ of $F(u)$ is congruent to $[a_{st}]$ where
$$ a_{ii} = \sinh(r_i) \frac{\partial \theta_i}{\partial r_i}$$ $$
 =\sinh(r_i) (\frac{\partial \theta_i}{\partial l_j} + \frac{ \partial \theta_i}{\partial l_k})$$
$$=-\sinh(r_i) \frac{ \partial \theta_i}{\partial l_i} (\cos(\theta_k) + \cos(\theta_j))$$
$$=-\frac{\sinh(r_i)\sin(\theta_i)}{A}(\cos(\theta_k)+ \cos(\theta_j))$$
and,
$$a_{ij} = \sinh(r_j)\frac{\partial \theta_i}{\partial r_j}$$
$$=\sinh(r_j)(\frac{\partial \theta_i}{\partial l_i} + \frac{\partial \theta_i}{\partial l_k})$$
$$=\sinh(r_j)\frac{\partial \theta_i}{\partial l_i}( 1-\cos(\theta_j))$$
$$ = \frac{2\sinh(r_j)\sin(\theta_i) \sin^2(\theta_j/2)}{A}.$$
From the above identities and $\theta_k + \theta_j < \pi$, we see
that $a_{ii} <0$ and $a_{ij}
>0$. Now $$|a_{ii}| - |a_{ji}| -|a_{ki}| $$
$$= \frac{\sinh(r_i)}{A} ( \sin(\theta_i) (\cos(\theta_k) +\cos(\theta_j)) - 2(\sin(\theta_j) + \sin(\theta_k))\sin^2(\frac{\theta_i}{2}))$$
$$=\frac{4 \sinh(r_i)  \sin(\frac{\theta_i}{2})}{A} (\cos(\frac{\theta_i}{2})( \cos(\frac{\theta_k+\theta_j}{2})\cos(\frac{\theta_k-\theta_j}{2}) - \sin(\frac{\theta_k+\theta_j}{2})\cos(\frac{\theta_k-\theta_j}{2}) \sin(\frac{\theta_i}{2}))$$
$$=\frac{4\sinh(r_i)\sin(\frac{\theta_i}{2})\cos(\frac{\theta_j-\theta_k}{2})}{A}( \cos(\frac{\theta_i}{2}) \cos(\frac{\theta_k+\theta_j}{2}) -\sin(\frac{\theta_i}{2}) \sin(\frac{\theta_k+\theta_j}{2}))$$
$$=\frac{4}{A} \sinh(r_i) \sin(\frac{\theta_i}{2}) \cos(\frac{\theta_j-\theta_k}{2})\cos(\frac{\theta_i+\theta_j+ \theta_k}{2}) >0.$$

This shows that $-[a_{st}]$ is diagonally dominate. Thus $H$ is negative definite and the integration is locally strictly concave.

\medskip
In the case (d), we take $y_i = \sqrt{-1} l_i$ and $x_i =\theta_i
= r_j + r_k$ in lemma 2.2. The tangent law says $\tanh(l_i/2)/
\cos(r_i)$ is independent of the index $i$. Furthermore the
Hessian $H=[H_{st}]$ of the integral
 $f(u) =\int^u \eta_{\lambda}$ where $f'(l) =\tanh^{\lambda}(l/2)$ and $g'(r) =\cos^{-\lambda-1}(r)$ is
$$ H_{st}= \frac{f'(l_s)}{g'(r_t)} \frac{ \partial l_s}{\partial r_t}$$
$$=(\frac{\tanh(l_s/2)}{\cos(r_s)})^{\lambda}(\cos(r_s) \cos(r_t))^{\lambda+1} \frac{1}{\cos(r_s)}\frac{ \partial l_s}{\partial r_t}$$
$$= (\frac{\tanh(l_s/2)}{\cos(r_s)})^{\lambda}(\cos(r_s) \cos(r_t))^{\lambda+1}  a_{st}.$$

Since $\tanh(l_s/2)/\cos(r_s) >0$ is independent of the indices,
thus the Hessian $H$ is congruent to $[a_{st}]$ where, by lemma
2.2,
$$a_{ij} = B' \cosh^2(l_i/2)\cosh^2(l_j/2)$$
and
$$a_{ii} = B' \cosh^2(l_i/2) ( \cosh^2(l_j/2) + \cosh^2(l_k/2) -1)$$
for some positive function $B'$ independent of indices. At the
equilateral triangle with $\cosh^2(l_i/2) = 2$ for all $i$, the
matrix $[a_{st}]$ is  a positive multiple of $[b_{st}]$ where
$b_{ij}=1$ and $b_{ii} = 3/2$. Thus $[a_{st}]$ is positive
definite at this triangle. This shows that the matrix $H$ is
positive definite for one triangle. But on the other hand, the
Hessian $H$ is non-degenerate. Indeed, the Hessian is the Jacobi
matrix of the gradient of $F(u)=\int^u \eta_{\lambda}$ which sends
$u$ to $(\int^{u_1} \tanh^{\lambda} (t/2) dt, \int^{u_2}
\tanh^{\lambda } (t/2) dt, \int^{u_3} \tanh^{\lambda} (t/2) dt)$.
The gradient map is a diffeomorphism since we can solve angle from
the length by the cosine law. It follows that $det(H) \neq 0$ for
all triangles. Since the space of all hyperbolic triangles is
connected and $H$ is symmetric, the signature of $H$ is constant
for all hyperbolic triangles. We noticed above that the signature
for one triangle is (3,0). It follows that $H$ is positive
definite. Thus the integration is locally strictly convex. Another
way to see the positive definiteness of the matrix $[a_{rs}]$ is
observed by Ren Guo. It goes as follows. First we have $a_{ii}
>0$. Next, it is easy to see that all principle $2 \times 2$
submatrices are diagonally dominated. Thus they are positive
definite. Finally, one can show directly that
 the determinant of $[a_{rs}]$ is  positive.

\medskip

\medskip \noindent \S4. {\bf  Convexity of
Moduli Spaces of Geometric Triangles}

\medskip
We determine the convexity of the spaces of all geometric
triangles in various parameterizations in this section.

\medskip
\noindent 4.1. {\it Moduli spaces of geometric triangles}
\medskip
Suppose $K^2$ is one of the three geometries $\bold S^2$, $\bold
E^2$ or $\bold H^2$. We  use $K^2(l, 3)$ (or $K^2(\theta, 3)$) to
denote the space of all $K^2$ triangles parameterized by the edge
lengths (or angles respectively). In particular, we have the
following simple lemma whose proof is obvious.

\noindent {\bf Lemma 4.1.} \it The moduli spaces $K^2(l,3)$ and
$K^2(\theta,3)$ are,

(a)  $\bold E^2(l,3) = \bold H^2(l, 3) =\{ (p_1, p_2, p_3) \in
\bold R^3 | p_i + p_j
> p_k, \{i,j,k\}=\{1,2,3\}\}$.

(b) $\bold S^2(l, 3) =\{ (p_1, p_2, p_3) \in \bold R^3 | p_i + p_j
> p_k, \{i,j,k\}=\{1,2,3\},   \quad  p_1+p_2+p_3 < 2 \pi$\}.

(c) $\bold E^2(\theta,3) =\{ (p_1, p_2, p_3) \in \bold R^3 | p_i >
0, p_1 + p_2+p_3 =\pi\}$.

(d) $\bold H^2(\theta,3) =\{(p_1, p_2, p_3) \in \bold R^3 | p_i >
0, p_1 + p_2+p_3 < \pi\}$.

(e) $\bold S^2(\theta,3)= \{(p_1, p_2, p_3) \in \bold R^3 |  p_1 +
p_2+p_3
> \pi, p_i + p_j < p_k+ \pi, \{i,j,k\}=\{1,2,3\}\}$. \rm

\medskip

Another related moduli space is the space $\bold
H^2(r+r'=\theta,3)$ of all hyperbolic triangles whose angles are
$p_i + p_j $. In this case the moduli space is
$$\bold H^2(r+r'=\theta, 3) =\{(p_1, p_2, p_3) \in \bold R^3 | p_i+p_j >0, i \neq j,  \quad \text{and} \quad p_1+p_2+p_3 <
\pi/2\}.
$$ Note that $\bold H^2(r+r'=\theta, 3) \subset (-\pi/2,
\pi/2)^3$.

It is very interesting to note the relationship between these
spaces and Clebsch-Gordan and quantum Clebsch-Gordan relations
from representation theory of Lie algebra $sl(2)$.

\medskip
\noindent 4.2. {\it New parameterizations of moduli spaces}

\medskip
\noindent By theorems 3.2 and 3.4, we have to consider spaces of
all geometric triangles parameterized by $u=(u_1, u_2, u_3)$ where
$u_i =h(l_i)$ with $h(t) =\int^t \sin^{\lambda}(s) ds$, etc. The
main result in this section addresses the convexity of the moduli
spaces in these coordinates $u$.
\medskip
To begin, let $x = h(t)$ be a smooth diffeomorphism from an
interval $I$ to an interval $J$ so that $h'(t) >0$.  Consider the
diffeomorphism $q=q(p): I^3 \to J^3$ defined by $q_i = h(p_i)$,
$i=1,2,3$. We are interested in the convexity of the images of the
moduli spaces $K^2(l, 3)$ or $K^2(\theta, 3)$, or $\bold
H^2(r+r'=\theta, 3)$ under the diffeomorphism $q=q(p)$.

\medskip
\noindent {\bf Proposition 4.2.} \it Let $\{i,j,k\}=\{1,2,3\}$.

(a) If $\lambda \leq -1$ and $h(t)=-\frac{1}{\lambda}
t^{-\lambda}$ or $h(t) =\int^t_1 \coth^{\lambda+1}(s/2) ds: \bold
R_{>0} \to J$, then the image of $\{(p_1, p_2, p_3) \in \bold R^3
| p_i+p_j > p_k\}$ under $q=q(p)$ is convex.

(b) If $\lambda \geq 0$ and $h(t)=\int^t_0 \cos^{\lambda}(s) ds:
(-\pi/2, \pi/2) \to J$, then the image of the set $\{(p_1,p_2,
p_3) \in \bold R^3|  p_i+p_j >0,  \quad  p_1+p_2+p_3 < \pi/2\}$
under $q=q(p)$ is convex.

(c) If $\lambda \geq 0$ and $h(t) = \int^t_{\pi/2}
\sin^{\lambda}(s) ds: (0, \pi) \to J$, then the image of $\{ (p_1,
p_2, p_3) \in \bold R^3 | p_i+ p_j > p_k, p_1+p_2+p_3 < 2\pi\}$
under $q=q(p)$ is convex.

(d) If $\lambda \geq 0$ and $h(t) = \int^t_{\pi/2}
\sin^{\lambda}(s) ds: (0, \pi) \to J$, then the image of $\{ (p_1,
p_2, p_3) \in \bold R^3 | p_i+ p_j < \pi+ p_k, p_1+p_2+p_3 >
\pi\}$ under $q=q(p)$ is convex.

\medskip

\rm \noindent {\bf Remarks 4.3.} 1.  Note that the moduli spaces
involved in the propositions are, listed in the order, $\bold
E^2(l,3), \bold H^2(l,3),$ $\bold H^2(r+r'=\theta,3)$, $\bold
S^2(l,3)$ and $\bold S^2(\theta, 3)$.

2. It can be shown, from the proof below, that if $\lambda$ does
not satisfy the inequalities in the proposition, then the images
are not convex.

\medskip
\noindent 4.3. {\it A lemma}
\medskip
\noindent We begin with the following lemma which determines the
convexity of a plane under the map $q=q(p)$.

\medskip
\noindent {\bf Lemma 4.4.} \it Suppose $\alpha, \epsilon \in \bold
R$ and $t=g(x): J \to I$ is the inverse of $x=h(t)$. Then the
second derivative and the determinant of the Hessian of the
function $z(x,y)=h(\alpha + \epsilon g(x) + \epsilon g(y))$ are
given by the following identities where $A(t) =\frac{
h''(t)}{h'(t)}$, $u=g(x), v=g(y)$ and $w=\alpha + \epsilon (u+v)$.

$$ \frac{ \partial^2 z}{\partial x \partial x} =\epsilon
\frac{h'(w)}{h'(u)^2} ( \epsilon A(w) -A(u))  \tag 4.1$$ and the
determinant of the Hessian of $z(x,y)$ is
$$ \epsilon^2 \frac{ h'(w)^2}{h'(u)^2 h'(v)^2}( A(u) A(v) -\epsilon A(w) ( A(u) +
A(v))).  \tag 4.2$$ \rm

\medskip
The proof is a routine calculus. Note that $g'(x) =
\frac{1}{h'(t)}$ and $g''(x) = -\frac{h''(t)}{h'(t)^3}$. We have
$$\frac{\partial z}{\partial x} = \epsilon h'(w)g'(x)$$
$$\frac{\partial^2 z}{\partial x \partial y} =\epsilon^2  h''(w) g'(x) g'(y)
= \epsilon^2 \frac{ h''(w)}{h'(u) h'(v)}$$ and
$$\frac{\partial ^2 z}{\partial^2 x} = \epsilon^2  h''(w) g'(x)^2
+ \epsilon h'(w)g''(x)$$
$$=\epsilon^2 \frac{ h''(w)}{h'(u)^2} - \epsilon \frac{ h'(w) h''(u)}{h'(u)^3}$$
$$=\epsilon \frac{h'(w)}{h'(u)^2}( \epsilon A(w)-A(u)).$$

The determinant of the Hessian of $h$ is
$$=\epsilon^2 (\frac{h'(w)}{h'(u) h'(v)})^2( (\epsilon A(w)-A(u))(\epsilon A(w)-A(v)))
-\epsilon^4 \frac{h''(w)^2}{h'(u)^2 h'(v)^2}$$ $$=\epsilon^2
\frac{h'(w)^2}{h'(u)^2 h'(v)^2}( \epsilon^2 A^2(w) - \epsilon A(w)
( A(u) + A(v)) + A(u)A(v) -\epsilon^2 A^2(w))$$
$$= \epsilon^2 \frac{ h'(w)^2}{h'(u)^2 h'(v)^2}( A(u) A(v) -\epsilon A(w) ( A(u) +
A(v))).$$

\medskip
\noindent 4.4. {\it A proof of proposition 4.2}

\medskip
\noindent Let $X$ be the image of the set under $q=q(p)$ in cases
(a)-(d), i.e.,

In the case (a),
$$ X =\{(q_1, q_2, q_3) \in \bold R^3 | q_k < h( g(q_i))+ g(q_j)),
\{i,j,k\}=\{1,2,3\}\}. \tag 4.3$$

In the case (b), due to $h(-t) = -h(t)$ and $h'(t)>0$, $s > -t$ if
and only if $h(s) > -h(t)$,
$$ X =\{(q_1, q_2, q_3) \in \bold R^3 | q_i+q_j
>0, q_k < h( \pi/2 -g(q_i)- g(q_j)), \{i,j,k\}=\{1,2,3\}\}. \tag 4.4$$

In the case (c), $X$ is
$$\{(q_1, q_2, q_3) \in \bold R^3 | q_k < h( g(q_i)+ g(q_j)),
q_3 < h(2 \pi-g(q_1) -g(q_2)),  \{i,j,k\}=\{1,2,3\}\}. \tag 4.5$$

In the case (d), $X$ is
$$\{(q_1, q_2, q_3) \in \bold R^3 | q_k > h(-\pi+ g(q_i))+ g(q_j)),
q_3> h(\pi-g(p_1))-g(p_2)),  \{i,j,k\}=\{1,2,3\}\}. \tag 4.6$$

To prove (a), by (4.3), $X$ is bounded by three surfaces. It
suffices to show that the surfaces bounding $X$ are convex in the
side containing $X$. By symmetry, it remains to show that the
function $z=h(g(x)+ g(y))$ is locally concave in $x,y$. By lemma
4.4 with $\alpha =0$, $\epsilon =1$, $(u, v)=(g(x), g(y)) \in
\bold R^2_{>0}$ and $w=u+v$, we find the determinant of the
Hessian and the sign of the second derivative $\partial^2
z/\partial^2 x$ as follows.

If $h(t) = -\frac{t^{-\lambda}}{\lambda}$, then $A(t)
=-\frac{\lambda+1}{t}$. By (4.2), the determinant of the Hessian
of $z=z(x,y)$ is $$\frac{h'(w)^2}{h'(u)^2 h'(v)^2}( A(u) A(v)
-A(u+v)(A(u)+A(v)))$$ $$=(\lambda+1)^2 \frac{ h'(w)^2}{h'(u)^2
h'(v)^2}( \frac{1}{uv}-\frac{1}{u+v}( \frac{1}{u}
+\frac{1}{v}))=0.$$
 The second derivative $\partial^2 z/\partial^2
x =\frac{h'(w)}{h'(u)^2}(A(w)-A(u))$ =$(\lambda+1)
\frac{h'(w)}{h'(u)^2} \frac{v}{(u+v) u}$. It is non-positive due
to $\lambda \leq -1$, $u,v >0$ and $h'(t)>0$. By symmetry, we see
that $\partial^2 z/\partial^2 y \leq 0$. Thus the Hessian of $z$
is semi-negative definite. This shows that the function
 $z=z(x,y)$ is locally concave. Thus the space $X$ in (4.3) is
 convex.

 \medskip

In the rest of the proof, we will check the signs of $\partial^2
z/\partial^2 x$ (not $\partial^2 z/\partial^2 y$) and the
determinant of the Hessian only.

If $h(t)=\int^t_1 \coth^{\lambda+1}(s/2) ds$, then $h'(t) =
\coth^{\lambda+1}(t/2)$, $h''(t) =-(\lambda+1)
\frac{\coth^{\lambda+1}(t/2)}{\sinh(t)}$. Thus $A(t) =
-\frac{\lambda+1}{\sinh(t)}$. By (4.2) for $\alpha=0, \epsilon
=1$, the determinant of the Hessian of $z$ is
$$\frac{h'(w)^2}{h'(u)^2 h'(v)^2}(\lambda+1)^2[
\frac{1}{\sinh(u)\sinh(v)}
-\frac{1}{\sinh(u+v)}(\frac{1}{\sinh(u)}+\frac{1}{\sinh(v)})]$$
$$=(\lambda+1)^2 \frac{h'(w)^2}{h'(u)^2 h'(v)^2} \frac{ \sinh(u+v)
-\sinh(u) -\sinh(v)}{\sinh(u)\sinh(v) \sinh(u+v)} \geq 0$$ due to
$u, v
>0$ and $\sinh(u+v) > \sinh(u) + \sinh(v)$.

The second derivative $\partial^2 z/\partial^2 x$ can be
calculated by (4.1) to be: $ \frac{ h'(w)}{h'(u)^2}( A(w) -A(u))$
$= -(\lambda+1) \frac{h'(w)}{h'(u)^2}(
\frac{1}{\sinh(u+v)}-\frac{1}{\sinh(u)})$. It is non-positive if
and only if $\lambda \leq -1$ due to $u, v >0$. This shows that
$z$ is locally concave in $x,y$.

This establishes case (a).

\medskip

To see part (b), by (4.4), it suffices to verify that $z=h( \pi/2
-g(x) -g(y))$ is locally concave in $x,y$ where $(u,v)=(g(x),
g(y)) \in (-\pi/2, \pi/2)^2$, and $u+v>0$. The last sets of
conditions on $u,v$ are due to the following. The point $(u,v,
\pi/2-u-v)$ is in the closure of  $\bold H^2(r+r', 3)=\{(p_1, p_2,
p_3)| p_i + p_j
>0, p_1+p_2+p_3 < \pi/2\} \subset [-\pi/2, \pi/2]^3$ so that $(u,v)$ is in the projection of
$\bold H^2(r+r',3)$ to the $p_1, p_2$ coordinates. It follows that
$u, v \in (-\pi/2, \pi/2)$ and $u+v > 0$. Now $h(t) =\int^t_0
\cos^{\lambda}(s) ds$, $h'(t) =\cos^{\lambda}(t)$, $h''(t)
=-\lambda \cos^{\lambda-1}(t) \sin(t)$ and $A(t) =-\lambda
\tan(t)$.

By lemma 4.4 with $\alpha =\pi/2, \epsilon =-1$ and $w
=\pi/2-u-v$, we find the determinant of the Hessian of $z$ to be:
 $$\frac{h'(w)^2}{h'(u)^2 h'(v)^2}( A(u) A(v)
+A(w)(A(u)+A(v)))$$
 $$=\lambda^2  \frac{ h'(w)^2}{h'(u)^2
h'(v)^2}( \tan(u)\tan(v) + \cot(u+v) (\tan(u)+\tan(v)))$$
$$=\lambda^2 \frac{h'(w)^2}{h'(u)^2 h'(v)^2}$$
It is non-negative. By (4.1), we find the second derivative
$\partial^2z/\partial^2 x$ to be
$$\frac{h'(w)}{h'(u)^2}(A(w)+A(u))$$
$$=-\lambda \frac{h'(w)}{h'(u)^2}( \cot(u+v) + \tan( u))$$
$$=-\lambda \frac{h'(w)}{h'(u)^2} \frac{ \cos(u+v) \cos(u) +
\sin(u+v) \sin(u)}{\sin(u+v)\cos(u)}$$
$$=-\lambda
\frac{h'(w)}{h'(u)^2} \frac{\cos(v)}{\sin(u+v)\cos(u)}.$$ Note
that due to $u, v \in (-\pi/2, \pi/2)$ and $u+v \in (0, \pi)$. The
second derivative is non-positive since $\lambda \geq 0$. Thus $z$
is locally concave in $x,y$.  This proves the case (b).

\medskip

In the cases (c) and (d), $h(t) = \int^t_{\pi/2} \sin^{\lambda}(s)
ds$. It follows that $h'(t) =\sin^{\lambda}(t)$, $h''(t)
=\lambda$$ \sin^{\lambda -1}(t)$$ \cos(t)$, and $A(t)=\lambda
\cot(t)$. Due to
$$ \cot(u)\cot(v) -\cot(u+v) (\cot(u)+\cot(v)) =1$$ we have
$$A(u)A(v) - A(u+v)(A(u) + A(v)) =\lambda^2. \tag 4.7$$
Due to $$\cot(u+v) -\cot(u) =-\frac{\sin(v)}{\sin(u+v)\sin(u)}$$
we have $$A(u+v) -A(u) =-\lambda \frac{\sin(v)}{\sin(u+v)\sin(u)}.
\tag 4.8$$

In the case (c), by (4.5), we must verify that (1)
$z=h(g(x)+g(y))$ is locally concave in $x,y$ where $(u, v)=(g(x),
g(y))$ satisfies $u,v \in (0, \pi)$ and $w=u+v \in (0, \pi)$ and
(2) $ z=h(2\pi-g(x) -g(y))$ is locally concave in $x,y$ where
$(u,v) =(g(x), g(y))$ satisfies $u,v \in (0, \pi)$ and $w
=2\pi-u-v \in (0, \pi)$, i.e, $u+v > \pi$.

Now for $z =h(g(x)+g(y))$, by (4.2) with $\alpha=0, \epsilon =1$,
the determinant of the Hessian of $z$ is $$\frac{h'(w)^2}{h'(u)^2
h'(v)^2}( A(u)A(v) -A(u+v)(A(u) + A(v)).$$ By (4.7), it is
$\lambda^2 \frac{ h'(w)^2}{h'(u)^2}{h'(v)^2}$ and is clearly
non-negative. The second derivative $\partial^2 z/\partial ^2 x$
can be calculated by (4.8) to be $\frac{h'(w)}{h'(u)^2}(A(u+v)
-A(u)) = -\lambda \frac{h'(w)}{h'(u)^2}\frac{\sin(v)}{\sin(u+v)
\sin(u)}$. It is non-positive since $\lambda \geq 0$, $u, v, u+v
\in (0, \pi)$. Thus $z$ is locally concave in $x,y$.

To see the local concavity of $z = h(2\pi-g(x) -g(y))$ where
$w=2\pi-u-v \in (0, \pi)$ and $u,v \in (0, \pi)$, $u+v >\pi$, we
use lemma 4.4 for $\alpha=2\pi$ and $\epsilon=-1$. By (4.2) and
(4.7), the determinant of the Hessian is
$$ \frac{h'(w)}{h'(u)^2 h'(v)^2}[ A(u) A(v) + A(w)(A(u)+A(v))]$$
$$=\frac{h'(w)}{h'(u)^2 h'(v)^2} [A(u)A(v) -A(u+v)(A(u)+ A(v))]$$
$$=\frac{h'(w)}{h'(u)^2 h'(v)^2}\lambda^2.$$
It is non-negative. The second derivative $\partial^2 z/\partial^2
x$ is found to be, by (4.1) and (4.8),
$$ \frac{h'(w)}{h'(u)^2}(A(w) + A(u))$$
$$=\frac{h'(w)}{h'(u)^2}(-A(u+v) + A(u))$$
$$ =\lambda \frac{h'(w)}{h'(u)^2} \frac{
\sin(v)}{\sin(u+v)\sin(u)}.$$ It is non-positive due to $\lambda
\geq 0$, $u,v \in(0, \pi)$ and $u+v \in (\pi, 2\pi)$. Thus $z$ is
locally concave in $x,y$.

In the case (d), for $\lambda \geq 0$, by (4.6), we must verify
that (1) $z = h( -\pi+ g(x) + g(y))$ is locally convex in $x,y$
where $(u, v)=(g(x), g(y)) \in (0, \pi)^2$ and $w =-\pi+u+v \in
(0, \pi)$, i.e., $u+v > \pi$ and (2) $z = h(\pi- g(x) -g(y))$ is
locally convex in $x,y$ where $(u, v)= (g(x), g(y)) \in (0,
\pi)^2$ so that $w =\pi-u-v \in (0,\pi)$, i.e., $u+v < \pi$.

For $z=h(-\pi+g(x) + g(y))$, by lemma 4.4 with $\alpha=-\pi,
\epsilon =1$,  $w=-\pi +u+v$ and (4.7), we find the determinant of
the Hessian of $z$ to be
$$ \frac{h'(w)}{h'(u)^2 h'(v)^2}[ A(u) A(v) - A(w)(A(u)+A(v))]$$
$$=\frac{h'(w)}{h'(u)^2 h'(v)^2} [A(u)A(v) -A(u+v)(A(u)+ A(v))]$$
$$=\frac{h'(w)}{h'(u)^2 h'(v)^2}\lambda^2 \geq 0$$

The second derivative $\partial^2 z/\partial^2 x$ is, by (4.1) for
$\alpha=\pi, \epsilon =-1$ and (4.8),
$$ \frac{h'(w)}{h'(u)^2}(A(w) - A(u))$$
$$=\frac{h'(w)}{h'(u)^2}(A(u+v) - A(u))$$
$$ = -\lambda \frac{h'(w)}{h'(u)^2} \frac{
\sin(v)}{\sin(u+v)\sin(u)}.$$ It is non-negative due to $\lambda
\geq 0$, $u, v \in (0, \pi)$ and $u+v \in (\pi, 2\pi)$. Thus
$z=z(x,y)$ is locally convex in $x,y$.

Finally,  for $z = h(\pi-g(x) -g(y))$ to be locally convex in
$x,y$ where $(u, v) = (g(x), g(y)) \in (0, \pi)$ and $w=\pi-u-v
\in (0, \pi)$, i.e., $ u+v < \pi$, we find the determinant of the
Hessian of $z$ to be, by (4.2) and (4.7),
$$ \frac{h'(w)}{h'(u)^2 h'(v)^2}[ A(u) A(v) + A(w)(A(u)+A(v))]$$
$$=\frac{h'(w)}{h'(u)^2 h'(v)^2} [A(u)A(v) -A(u+v)(A(u)+ A(v))]$$
$$=\frac{h'(w)}{h'(u)^2 h'(v)^2}\lambda^2 \geq 0.$$

The second derivative $\partial^2 z/\partial^2 x$ is found to be,
by (4.1) and (4.8),
$$ \frac{h'(w)}{h'(u)^2}(A(w) + A(u))$$
$$=\frac{h'(w)}{h'(u)^2}(-A(u+v) + A(u))$$
$$ = \lambda \frac{h'(w)}{h'(u)^2} \frac{
\sin(v)}{\sin(u+v)\sin(u)}.$$ It is non-negative due to $\lambda
\geq 0$, $u, v,$ $u+v$ are in $(0, \pi)$. This shows that $z$ is
locally convex in $x,y$.

\medskip

\newpage

\noindent \S5. {\bf Polyhedral Surfaces with or without Boundary}
\medskip
We now set up the framework for both triangulated closed surfaces
and ideal  triangulated compact surfaces with boundary.

\medskip
\noindent 5.1. {\it Triangulation and ideal triangulation}

\medskip
By a closed triangulated surface $(S, T)$ we mean the following.
Take a finite disjoint union of Euclidean triangles and identify
all edges of triangles in pairs by homeomorphisms. The quotient is
a closed surface $S$ (possibly disconnected) together with a
triangulation $T$. The cells in $T$ are the quotients of the
vertices, edges and triangles in the disjoint union. A \it
simplicial triangulation \rm of a surface is a triangulation so
that each (closed) cell is homeomorphic to a simplex and the
intersection of any two cells is either empty or a single cell.
 Let $V= T^{(0)}$, $E=T^{(1)}$ and $T^{(2)}$ be the sets of
all vertices, edges and triangles in the triangulation $T$
respectively.  Given a Euclidean triangle $\sigma$ with vertices
$v_1, v_2, v_3$, the \it corner \rm of $\sigma$ at vertex $v_i$ is
the collection of all open sets in the interior of the triangle
having $v_i$ as a limit point. We call the edge $v_jv_k$ (or edges
$v_iv_j$, $v_iv_k$) \it facing \rm (or \it adjacent to\rm) the
corner at $v_i$ where $\{i,j,k\}=\{1,2,3\}$. For a 2-cell, also
called a triangle, $\sigma$ in the triangulation $T$ of a surface,
a \it corner \rm of $\sigma$ is the image of a corner in the
unidentified space. Thus each triangle $\sigma \in T^{(2)}$ has
three corners even if the three vertices of triangle are one
point.
 Every corner in $T$ is
facing an edge and is adjacent to two edges $e, e'$ (it may occur
$e =e'$).
%We use $cor(S, T)$ to denote the set of all corners in a
%closed triangulated surfaces.

For a compact surface with non-empty boundary, the most efficient
way of decomposing it is the \it ideal triangulation. \rm  Here is
the definition. Take a closed triangulated surface $(X, T^*)$. Let
$N(V)$ be a small open regular neighborhood of the union of all
vertices. Then $S = X -N(V)$ is a compact surface with $|V|$ many
boundary components. The set $T =\{ \sigma \cap S | \sigma \in
T^*\}$ is called an \it ideal triangulation \rm of the surface
$S$. The set of all edges (respectively 2-cells, or hexagons) in
$T$ is $T^{(1)}= \{ e \cap S | e$ an edge in $T^*$\} (respectively
$T^{(2)}= \{ \sigma \cap S | \sigma $ a triangle in $T^*$ \}).
Each hexagon in $T$ has three edges counted multiplicity (even
though two of the three edges may be the same). The intersection
of a hexagon with the boundary $\partial S$ consists of three arcs
called \it B-arcs. \rm A B-arc in a hexagon $\sigma^2 \cap S$
corresponds to a corner in the triangle  $\sigma^2$ in $T^*$. The
notion of a B-arc \it facing \rm (or \it adjacent to\rm) an edge
$e$ in $T$ is defined as before. For instance, each B-arc is
facing exactly one edge and is adjacent to two edges counted with
multiplicities.

Another way of introducing ideal triangulation is as follows. A
\it colored hexagon \rm is a hexagon with three pairwise
non-adjacent edges colored by red and the other edges colored by
black. Take a finite disjoint union of Euclidean  colored hexagons
and identify all red edges in pairs by homeomorphisms. The
quotient is a  compact surface (possibly disconnected) with
non-empty boundary together with an ideal triangulation. The
2-cells in the ideal triangulation are quotients of the hexagons.
The quotients of red-edges (respectively black-edges) are the
edges (respectively B-arcs) in the ideal triangulation.

It is well known that each compact surface $\Sigma$ with $\partial
\Sigma \neq \emptyset$ and negative Euler characteristic admits an
ideal triangulation.

In the sequel, we assume that all surfaces are connected. A
triangle (or hexagon) in a triangulation (or ideal triangulation)
is counted to have three vertices $v_i, v_j, v_k$ and edges $e_i,
e_j, e_k$ even if $v_i=v_j$ or $e_i = e_j$.

\medskip
\noindent 5.2. {\it Polyhedral metrics}
\medskip
 Given a closed
triangulated surface $(S, T)$ with $E$ as the set of all edges, a
Euclidean (or hyperbolic, or spherical) polyhedral metric on $(S,
T)$ is characterized by its edge length function $l: E \to \bold
R_{>0}$ so that whenever  edges $e_i, e_j ,e_k$ form a triangle in
$T$, the three numbers $ l(e_i), l(e_j), l(e_k)$ form the edge
lengths of a Euclidean (respectively hyperbolic, or spherical)
triangle.  Let $K^2$ be $\bold E^2$ or $\bold S^2$ or $\bold H^2$
and let $P_{K^2}(S, T)$ be the space of all $K^2$ polyhedral
metrics on $(S, T)$ parameterized by the edge length function. It
follows that $P_{K^2}(S, T)$ is an open convex polytope in $\bold
R^E_{>0}$. Recall that the discrete curvature $k$ of a polyhedral
metric is the map $k: V \to \bold R$ sending each vertex to $2\pi$
minus the sum of all inner angles at the vertex. The \it curvature
map \rm $\Pi : P_{K^2}(S, T) \to \bold R^V$ sending a metric $l$
to its discrete curvature  $k$.

The following question concerning the metric-curvature
relationships may have an affirmative answer.

\medskip
\noindent {\bf Problem 5.1.} \it Suppose $(S, T)$ is a closed
triangulated surface. Let $\Pi: P_{K^2}(S, T) \to \bold R^V$ be
the curvature map and $p \in \bold R^V$.

(a) For $K^2=\bold E^2$ or $\bold H^2$, the space $\Pi^{-1}(p)$ is
either the empty set or a smooth manifold diffeomorphic to $\bold
R^{|E|-|V|}$.

(b) For $K^2=\bold S^2$, the space $\Pi^{-1}(p)$ is either the
empty set or a smooth manifold diffeomorphic to $\bold R^{ |E|-|V|
+ \mu}$ where $\mu$ is the dimension of the group of conformal
automorphisms of a spherical polyhedral metric $l \in
\Pi^{-1}(p)$. \rm

\medskip
Given a spherical polyhedral metric $l$ on $(S, T)$, let $V'$ be
the set of all vertices so that the discrete curvatures at the
vertices are zero. The number $\mu$ above is the dimension of the
group of all conformal automorphisms of the Riemann surface $S-V'$
where the conformal structure is induced by $l$.
  In particular, if the Euler
characteristic of $ S-V'$ is negative, then $\mu=0$.

Using the work of Rivin [Ri] and Leibon [Le], we will prove in
proposition 10.1 that the curvature maps $\Pi: P_{\bold H^2}(S, T)
\to \bold R^V$ and $\Pi: P_{\bold E^2}(S, T) \to \bold R^{|V|-1}$
are submersions. In particular, the preimage $\Pi^{-1}(p)$ is
either empty or a smooth manifold of dimension $|E|-|V|$. A more
detailed discussion of problem 5.1 can be found in subsection
10.1.

\medskip
\noindent 5.3. {\it Hyperbolic metric with geodesic boundary}
\medskip
 Hyperbolic metrics on ideal
triangulated surfaces are related to the Teichm\"uller theory of
surfaces with boundary.
 Here one replaces hyperbolic triangles by
the colored hyperbolic right-angle hexagons. The geometric
realization of a colored hexagon is based on the following well
known lemma.

\medskip
\noindent {\bf Lemma 5.1.} (see [Bu], [IT]) \it For any $l_1, l_2,
l_3 \in \bold R_{>0}$,  there exists a unique colored hyperbolic
right-angled hexagon whose three pairwise non-adjacent red edges
have lengths $l_1, l_2, l_3$. \rm

 \medskip

Given an ideal  triangulated surface $(S, T)$ with $E =T^{(1)}$
and a function $l: E \to \bold R_{>0}$, one uses lemma 5.1 to
produces a hyperbolic metric with totally geodesic boundary on $S$
having $l$ as the edge length function. This metric is constructed
by making each 2-cell $H$ in $T^{(2)}$ with edges $e_i, e_j, e_k$
a colored hyperbolic right-angled hexagon with red edge lengths
$l(e_i), l(e_j), l(e_k)$. Conversely, each hyperbolic metric with
geodesic boundary on $S$ is isometric to one constructed above.
Thus, the Teichm\"uller space of all hyperbolic metrics with
geodesic boundary on $S$, denoted by $Teich(S)$, can be identified
with the space $\bold R_{>0}^E$ by the edge length
parameterization $l: E \to \bold R_{>0}$. See [Us]. The
corresponding prescribing curvature problem (i.e., problem 5.1)
becomes hyperbolic metrics with given boundary lengths. Due to the
work of Fenchel-Nielsen (see [IT]), the space of all hyperbolic
metrics with given boundary lengths is known to be diffeomorphic
to $\bold R^{|E| -r}$ where $r$ is the number of boundary
components of $S$. In [Lu2], a new parametrization of it was
produced and a different proof of it was given. This gives an
evidence to the affirmative solution of problem 5.1.

\medskip

\medskip
\noindent \S6. {\bf Rigidity and Local Rigidity of Polyhedral
Surfaces}

\medskip
 We use the convex or concave energy functions in
theorems 3.2 and 3.4 to prove a collection of rigidity theorems
for polyhedral surfaces.

The main technical tool is the following well known fact from
analysis.

\medskip
\noindent {\bf Lemma 6.1.} \it Suppose $\Omega \subset \bold R^n$
is an open convex set and $W: \Omega \to \bold R^n$ is a smooth
locally strictly convex function. Then the gradient
$\bigtriangledown W: \Omega \to \bold R^n$ is a smooth embedding.
If $\Omega$ is only assumed to be open in $\bold R^n$, then
$\bigtriangledown W: \Omega \to \bold R^n$ is a local
diffeomorphism. \rm

\medskip
\noindent 6.1. {\it The main rigidity result} \medskip \noindent
Suppose $(S, T, l)$ is a $K^2$ polyhedral surface where $K^2=\bold
E^2, \bold S^2$ or $\bold H^2$ and $E$ is the set of all edges in
$T$. Recall that the space of all equivalence classes of $K^2$
polyhedral metrics on $(S,T)$, parameterized by the edge length
functions, is the open convex polytope $P_{K^2}(S, T) \subset
\bold R^E$. Define three maps $\Phi_{\lambda}, \Psi_{\lambda}$ and
$K_{\lambda}$ on $P_{K^2}(S, T)$ as follows. The  map
$\Phi_{\lambda}: P_{K^2}(S,T) \to \bold R^E$ sends a metric to its
$\phi_{\lambda}$ edge invariant defined by (1.1).
 The map
$\Psi_{\lambda}: P_{K^2}(S,T) \to \bold R^E$  sends the metric $l$
to its $\psi_{\lambda}$ edge invariant defined by (1.2). The
 map  $K_{\lambda}: P_{K^2}(S, T)
\to \bold R^V$  sends a polyhedral metric $l$ to its $k_{\lambda}$
discrete curvature defined by (1.3).

Let $\bold R_{>0}$ act on $\bold R^E$ by multiplication. Then
$P_{E^2}(S, T) \subset \bold R^E$ is invariant under the action.
The orbit space, denoted by $P_{ E^2}(S, T)/\bold R_{>0}$, is the
set of all Euclidean polyhedral metrics on $(S, T)$ modulo
scaling. By definition, all
 maps $\Phi_{\lambda}$, $\Psi_{\lambda}$ and $K_{\lambda}$
defined on $P_{\bold E^2}(S,T)$ are homogeneous of degree 0 with
respect to this action, i.e., they satisfy the equation $\phi(k x)
=\phi(x)$ for all $k \in \bold R_{>0}$. We use the same notations
$\Phi_{\lambda}$, $\Psi_{\lambda}$ and $K_{\lambda}$ to denote the
induced maps from $P_{\bold E^2}(S,T)/\bold R_{>0}$ to $ \bold
R^E$ or to $\bold R^V$.  We use $CP_{K^2}(S, T)$ to denote the
space of all $K^2$ circle packing metrics on $(S, T)$ for $K^2=
\bold E^2, \bold H^2$ or $\bold S^2$. The space of all  Euclidean
circle packing metrics modulo scaling is denoted by $CP_{E^2}(S,
T)/\bold R_{>0}$. The curvature map $K_{\lambda}$ is well  defined
on $CP_{E^2}(S, T)/\bold R_{>0}$.

\medskip
\noindent {\bf Theorem 6.2.} \it Suppose $(S, T)$ is a
triangulated closed surface and $\lambda \in \bold R$.

(a) The map $\Phi_{\lambda}: P_{ S^2}(S, T) \to \bold R^E$ is a
local diffeomorphism. It is a smooth embedding if $\lambda \geq 0$
or $\lambda \leq -1$.

(b) The map  $\Phi_{\lambda}: P_{E^2}(S, T)/\bold R_{>0} \to \bold
R^E$ is an immersion.

(c) If $\lambda \leq -1$, then $\Phi_{\lambda}: P_{E^2}(S,
T)/\bold R_{>0} \to \bold R^E$ is a smooth embedding.

(d) For any $\lambda \in \bold R$, $K_{\lambda}: CP_{H^2}(S, T)
\to \bold R^V$ is a smooth embedding.

(e) For any $\lambda \in \bold R$, $K_{\lambda}:  CP_{E^2}(S,
T)/\bold R_{>0} \to \bold R^V$ is a smooth embedding.

(f) The map $\Psi_{\lambda}: P_{ H^2}(S, T) \to \bold R^E$ is a
local diffeomorphism. If $\lambda \leq -1$ or $\lambda \geq 0$, it
is a smooth embedding.

\rm

\medskip

We remark that for $\lambda=0$, theorem 6.2 (d), (e), (f) and (a)
were first proved by Thurston [Th], Leibon [Le] and [Lu2]. Whether
parts (a), (c), (f)  hold for all $\lambda \in (-1, 0)$ is a very
interesting question. Theorem 6.2(c) does not cover the rigidity
result of Rivin.

\medskip
Let $E=\{e_1, ..., e_n\}$ be the set of all edges in the
triangulation $T$.  If $f: E \to X$, then $f_i$ denotes $f(e_i)$.
\medskip

 \noindent 6.2. {\it Proof of theorem 6.2(a)}
\medskip
\noindent Let $h(t) =\int^t_{\pi/2} \sin^{-\lambda-1}(x) dx$ for
$t \in (0, \pi)$. Then $h'(t) >0$ on $(0, \pi)$ and $h(t)$ is
strictly increasing. Given an edge length function $l: E \to (0,
\pi)$ of an  $\bold S^2$ polyhedral metric with $l_i =l(e_i)$,
define $u: E \to \bold R$ to be $u_i = h(l_i)$ and write $u=(u_1,
..., u_n)$. Then the map $u=u(l): P_{\bold S^2}(S, T) \to \bold
R^E$ is a smooth embedding. Let $\Omega=u(P_{S^2}(S, T))$ which is
open in $\bold R^E$. Recall from theorem 3.2 that if $l_i, l_j,
l_k$ are the edge lengths of a spherical triangle with opposite
angles $\theta_i, \theta_j, \theta_k$, then the differential
1-form
$$w_{\lambda} = \int^{\theta_i}_{\pi/2}
\sin^{\lambda}(t)dt  du_i + \int^{\theta_j}_{\pi/2}
\sin^{\lambda}(t)dt  du_j + \int^{\theta_k}_{\pi/2}
\sin^{\lambda}(t)dt  du_k$$ is closed and its integration $F(u_i,
u_j, u_k) =\int^{(u_i, u_j, u_k)}_{(\pi/2, \pi/2, \pi/2)}
w_{\lambda}$ is a locally strictly convex function in $(u_i, u_j,
u_k)$. Define an energy function $W: \Omega \to \bold R$ by
$$ W(u) = \sum_{ \{e_a, e_b, e_c\} \in T^{(2)}}  F(u_a, u_b, u_c)$$ where
the sum is over all triangles $\{e_a, e_b, e_c\}$ in the
triangulation $T$ with edges $e_a, e_b, e_c$. By the construction,
$W$ is locally strictly convex and
$$\frac{ \partial W}{\partial u_i} = \int^{\alpha}_{\pi/2}
\sin^{\lambda}(t) dt + \int^{\beta}_{\pi/2} \sin^{\lambda}(t) dt$$
where $\alpha, \beta$ are the angles facing the edge $e_i$,  i.e.,
$$\bigtriangledown W =(\phi_{\lambda}(e_1), ...,
\phi_{\lambda}(e_n)).$$

By lemma 6.1 applied to $W$ on $\Omega$, the first part of theorem
6.2(a) follows.

To prove the global rigidity for $\lambda \leq -1$, by lemma 6.1,
it suffices to prove that $\Omega$ is convex.  The convexity of
$\Omega$ follows from proposition 4.2(c). Indeed, since $-\lambda
-1 \geq 0$, by proposition 4.2(c), for the map
$h(t)=\int^t_{\pi/2} \sin^{-\lambda -1}(s) ds$, the image of
$\bold S^2(l, 3) =\{ (l_1, l_2, l_3) | l_i + l_j > l_k,
l_1+l_2+l_3 < 2\pi, \{i,j,k\}=\{1,2,3\}\}$ under $u=u(l)$ is an
open convex set $X$ in $\bold R^3$. Be definition, $\Omega$ is the
intersection of the open convex set $X^{T^{(2)}} \subset \bold R^{
3 |T^{(3)}|}$ with affine subspaces. (Each of the affine subspace
is of the form $x_i=x_j$ in the Euclidean space.) It follows that
$\Omega$ is convex in the case of $\lambda \leq -1$.

In the case of $\lambda \geq 0$, to prove that $\Phi_{\lambda}$ is
an embedding, we need to use the Legendre transform of the energy
function and the notion of angle structures. The proof is longer
and more complicated. We defer the proof to subsection 6.10 so
that it will coincide with the same method of the proof for part
(f) with $\lambda \geq 0$.

\medskip
\noindent 6.3. {\it A proof of theorem 6.2(b)}
\medskip
For a $\bold E^2$ polyhedral metric $l: E \to \bold R_{>0}$, let
$u_i =h(l_i)$ where $u_i = -\frac{1}{\lambda} l_i^{-\lambda}$ if
$\lambda \neq 0$ and $u_i = \ln( l_i)$ if $\lambda =0$.  The space
of all $\bold E^2$ polyhedral metrics on $(S, T)$, parameterized
by the edge length function, is an open convex polytope $P_{
E^2}(S, T)$ in $\bold R^E_{>0}$. The map $u=u(l)$ sends $P_{
E^2}(S, T)$ onto an open set $\Omega \subset \bold R^E$. Since $
h(\mu l_i) = \mu^{-\lambda} h(l_i)$ for $\lambda \neq 0$ and
$h(\mu l_i) =h(l_i) + \ln(\mu)$ for $\lambda=0$, the space
$\Omega$ has the following property. If $\lambda \neq 0$, then for
any positive number $c \in \bold R_{>0}$, $c \Omega =\{ cx | x \in
\Omega\} =\Omega$. If $\lambda=0$, then for any $c \in \bold R$,
$\Omega +c(1,1,..., 1) =\{ x+c(1,1,..., 1) | x \in \Omega \}
=\Omega.$ It follows that the space $P_{E^2}(S, T)/\bold R_{>0}$
of all Euclidean polyhedral metrics modulo scaling is
homeomorphic, under the map $u=u(l)$, to  the set $\Omega \cap
P_{\lambda}$ where $P_{\lambda} =\{ u=(u_1, ..., u_n) \in \bold
R^E | \sum_{i=1}^n u_i = -\lambda\}$.

\medskip
By theorem 3.2(a), if $l_i, l_j, l_k$ are the edge lengths of a
Euclidean triangle with corresponding angles $\theta_i, \theta_j,
\theta_k$, then the differential 1-form
$$\omega_{\lambda} =\int^{\theta_i}_{\pi/2} \sin^{\lambda}(t) dt du_i+
\int^{\theta_j}_{\pi/2} \sin^{\lambda}(t) dt du_j +
\int^{\theta_k}_{\pi/2} \sin^{\lambda}(t) dt du_k$$ is closed and
its integration $F(u_i, u_j, u_k) =\int^{(u_i, u_j, u_k)}_{(\pi/2,
\pi/2, \pi/2)} \omega_{\lambda}$ is a locally convex function
whose Hessian has null space generated by $(u_i, u_j, u_k)$ if
$\lambda \neq 0$ and by $(1,1,1)$ if $\lambda =0$.

Now for $u=(u_1, ..., u_n) \in \Omega$, define the energy function
$$W(u) = \sum_{\{e_a, e_b, e_c\} \in T^{(2)}} F(u_a, u_b, u_c)$$ where the
sum is over all triangles $\{e_a, e_b, e_c\}$ in $T$ with edges
$e_a, e_b, e_c$. By the construction, the function $W: \Omega \to
\bold R$ is locally convex whose gradient $\bigtriangledown W$ is
the $\phi_{\lambda}$ edge invariant. To establish theorem 6.2(b),
it suffices to prove that the map $\bigtriangledown W$ restricted
to $\Omega \cap P_{\lambda}$ is locally injective.

\medskip
\noindent {\bf Lemma 6.3.} \it The restriction map $W |: \Omega
\cap P_{\lambda} \to \bold R$ is locally strictly  convex. In
particular, the associate gradient map $\bigtriangledown (W|):
\Omega \cap P_{\lambda} \to P_0  =\bold R^{n-1}$ is a local
diffeomorphism. \rm

\medskip
\noindent {\bf Proof.} Let $A: \bold R^n \to P_0$ be the
orthogonal projection. Then by definition, $\bigtriangledown (W|)
=A( \bigtriangledown(W)|)$ on the subspace $\Omega \cap
P_{\lambda}$. Furthermore, the restriction of the Hessian $H(W)$
to $P_0 \times P_0$ is the Hessian matrix $H(W|)$ of $W|$. Thus it
suffices to show that the null space of the Hessian Hess(W) at a
point $u$ is transverse to the plane $P_0$. To see this, take a
vector $v=(v_1, ..., v_n) \in \bold R^n$ so that
$$ v Hess(W)|_u v^t=0$$
where $v^t$ is the transpose of the row vector $v$. By definition
of the function $W$, the above is equivalent to
$$ \sum_{\{e_a, e_b, e_c\} \in T^{(2)}} (v_a, v_b, v_c) Hess(F)|_{(u_a,
u_b,u_c)} (v_a, v_b, v_c)^t =0.$$ Each term in the summation is
non-negative due to convexity of $F$. It follows that
$$(v_a,v_b,v_c) Hess(F)|_{(u_a,u_b, u_c)} (v_a, v_b, v_c)^t =0.$$ By corollary 3.3, there is a
constant $C_{\{a,b,c\}}$ depending only on the triangle
$\{e_a,e_b, e_c\}$ so that if $\lambda \neq 0$,
$$(v_a,v_b, v_c) = C_{\{a,b,c\}} (u_a,u_b, u_c) \tag 6.1$$
and if $\lambda =0$,
$$(v_a,v_b, v_c) = C_{\{a,b,c\}}(1,1,1). \tag 6.2$$
Since the surface is assumed to be connected, by comparing above
equations (6.1) (or (6.2)) at two triangles sharing the same edge
and that $u_i \neq 0$ when $\lambda \neq 0$, we conclude that
$C_{\{a,b,c\}}$ is a constant independent of the choice of the
triangles. Thus there is a constant $C \in \bold R$ so that $v=Cu$
if $\lambda \neq 0$ or $v=C(1,1,...,1)$ if $\lambda =0$. This
vector $v$ is not in the subspace $P_0$ unless $v=0$. This shows
that $W|: \Omega \cap P_{\lambda} \to \bold R$ is locally strictly
convex. QED

\medskip
By lemmas 6.1 and 6.2, it follows that the composition $A
(\bigtriangledown (W|) : \Omega \cap P_{\lambda} \to \bold
R^{n-1}$ is a local diffeomorphism.  Thus $\bigtriangledown (W)|:
\Omega \cap P_{\lambda} \to \bold R^{n}$ is locally injective.
This proves theorem 6.2(b).

\medskip
\noindent 6.4.  {\it A proof of theorem 6.2(c)}

\medskip
To prove part (c), it suffices to show that if $\lambda \leq -1$,
then the open set $\Omega$ is convex in $\bold R^E$. Then theorem
6.2(c) follows from lemma 6.1.

Let $X$ be the image of $E^2(l,3)=\{(l_1, l_2, l_3) | l_i + l_j
> l_k$ for $\{i,j,k\} =\{1,2,3\}$ \} under the map $u=(u_1, u_2,
u_3)$ where $u_i =-\frac{1}{\lambda} l_i^{-\lambda}$.  By
proposition 4.2(a), $X$ is open convex in $\bold R^3$. By
definition, the open set $\Omega$ is the
 intersection of $X^{ T^{(2)}}$ with some
affine spaces. It follows that for $\lambda \leq -1$, $\Omega$ is
convex. Thus theorem 6.2(c) holds.

\medskip

\medskip
\noindent 6.5. {\it A proof of theorem 6.2(d)}

\medskip
The proof is straight forward due to the strictly convexity of the
energy functional in theorem 3.4(b) (by replacing $\lambda$ by
$-\lambda$). Namely, given a hyperbolic triangle of edge lengths
$l_k=r_i+r_j$, $\{i,j,k\}=\{1,2,3\}$, and opposite angles
$\theta_k$, the differential 1-form
$$\eta_{\lambda} = \sum_{i=1}^3 \int^{\theta_i}_{\pi/2}
\tan^{\lambda}(t/2) dt du_i$$ is closed where $u_i =\int^{r_i}_1
\sinh^{\lambda-1}(t) dt$. Furthermore, the integration $F(u_1,
u_2, u_3) =\int^u \eta_{\lambda}$ is locally strictly concave in
$u=(u_1, u_2, u_3)$. By definition, we have
$$\frac{\partial F}{\partial u_i} =\int^{\theta_i}_{\pi/2}
\tan^{\lambda}(t) dt.$$

Let $V=\{ v_1, ..., v_m\}$ be the set of all vertices in the
triangulation. For a hyperbolic circle packing metric $r: V \to
\bold R_{>0}$, define $u=u(r): V \to \bold R$ by $u_i =h(r_i)$
where $h(x) =\int^{x}_1 \sinh^{\lambda-1}(t) dt$. The image of
$\bold R^V_{>0}$ under the map $u=u(r)$ is the open cube $I^V$
where $I =h(\bold R_{>0})$.  Define a smooth function $W$ on $I^V$
by
$$ W(u_1, ..., u_m) =\sum_{ \{v_a, v_b, v_c\} \in F^{(2)}} F(u_a,u_b, u_c) \tag 6.3$$
where the sum is over all triangles with vertices $v_a, v_b, v_c$.

By the construction, $W$ is locally strictly convex in the open
convex set $I^V$ and its gradient is the $\lambda$-th discrete
curvature $k_{\lambda}$.

Thus theorem 6.2(d) follows from lemma 6.1 applied to the energy
function $W$.

\medskip
\noindent 6.6. {\bf Remark.}  For $\lambda=0$, the above proof was
first given by Colin de Verdiere [CV1]. Our proof uses exactly the
same method pioneered by him.

\medskip
\noindent 6.7. {\it A proof of theorem 6.2(e)}
\medskip
\noindent The proof of part(e) is essentially the same as that of
parts (b),(c). We sketch the main steps. First, by theorem 3.4(a)
(by replacing $\lambda$ by $-\lambda$), the integration
$F(u)=\int^u \eta$ of the closed 1-form
$$\eta=\sum_{i=1}^3 \int^{\theta_i} \tan^{\lambda}(t/2) dt du_i$$
is locally convex where $u_i =\frac{1}{\lambda} r_i^{\lambda}$ for
$\lambda \neq 0$ and $u_i =\ln r_i$ for $\lambda=0$. Furthermore,
by corollary 3.5, the null space of the Hessian of $F(u)$ at a
point $u$ is generated by $u$ if $\lambda \neq 0$ and by $(1,1,1)$
if $\lambda =0$. Now for  a Euclidean circle packing metric $r: V
\to \bold R_{>0}$, define a new function $u: V \to \bold R$ by
$u_i= h(r_i)$ where $h(t) = \frac{1}{\lambda} t^{\lambda}$ for
$\lambda\ \neq 0$ and $h(t) =\ln (t)$ for $\lambda =0$. We write
$u=(u_1(r), ..., u_m(r)) \in \bold R^E$. The image of $\bold
R_{>0}^V$ under $u=u(r)$ is an open convex cube $I^V$ where $I
=h(\bold R_{>0})$.
 Define a function $W$ on $I^V$ by the same
formula (6.3). Then this function $W$ is concave with gradient
$\bigtriangledown W$ equal to the $\lambda$-th discrete curvature.
The space $CP_{E^2}(S,T)/\bold R_{>0}$ of all  circle packing
metrics modulo scaling is homeomorphic to   $I^V \cap P_{\lambda}$
where $P_{\lambda} = \{ (u_1, ..., u_m) \in \bold R| \sum_{i=1}^k
u_i = \lambda$\} under $u=u(l)$. Now due to corollary 3.5 and by
the same argument as the one used in the proof of part 6.2(b), the
function $W|: I^V \cap P_{\lambda} \to \bold R$ is strictly
concave. Thus by the convexity of $I^V \cap P_{\lambda}$, the map
$\bigtriangledown W |: I^V \cap P_{\lambda} \to \bold R^V$ is an
embedding.

\medskip
\noindent 6.8. {\it A proof theorem 6.2(f) when $\lambda \leq -1$}

\medskip
\noindent Take the Legendre transformation of the integration of
the 1-form in theorem 3.4(d) (replacing $\lambda$ by $-\lambda
-1$). We obtain a locally strictly convex function $F(u)$ defined
as follows. For a hyperbolic triangle of edge lengths $l_i, l_j,
l_k$ and opposite angles $\theta_i =r_j+r_k$, the differential
1-form
$$ \eta^*_{\lambda} = \sum_{i=1}^3 \int^{r_i}_0 \cos^{\lambda}(t)
dt du_i$$ is closed where $u_i =\int^{l_i}_1
\coth^{\lambda+1}(t/2) dt$. The integration $F(u_1, u_2,
u_3)=\int^{(u_1, u_2, u_3)}_{(0,0,0)} \eta^*_{\lambda}$ is locally
strictly convex in $u$ by theorem 3.4(d). Furthermore, by the
construction,
$$ \frac{\partial F(u)}{\partial u_i} =\int^{r_i}_0
\cos^{\lambda}(t) dt.$$  For each hyperbolic metric $l: E \to
\bold R_{>0}$, define $u=(u_1, ..., u_n)$ by $u_i = h(l_i)$ where
$h(t) =\int^t_1 \coth^{\lambda+1}(s/2) ds$. Let $\Omega \subset
\bold R_{>0}^E$ be the image of $P_{\bold H^2}(S, T)$ under the
map $u=u(l)$. Define an energy function $W(u): \Omega \to \bold R$
by the formula
$$ W(u_1, ..., u_n) =\sum_{ \{e_a, e_b, e_c\} \in T^{(2)}} F(u_a, u_b, u_c)$$ where the sum is over all triangles with edges $e_a$, $e_b$
and $e_c$.

By the construction, the function $W: \Omega \to \bold R$ is
locally strictly convex so that
$$ \frac{\partial W}{\partial u_i} =\psi_{\lambda}(e_i)$$ i.e.,
$\bigtriangledown W$ is the $\psi_{\lambda}$ edge invariant. By
lemma 6.1 and the fact that $\Omega$ is open, $\bigtriangledown W:
\Omega \to \bold R^E$ is locally injective. It follows that a
hyperbolic metric is locally determined by its $\psi_{\lambda}$
edge invariant. For $\lambda \leq -1$, the open set $\Omega$ is
convex in $\bold R^E$. Indeed, by the same argument as in
subsection 6.2,
 the convexity of
$\Omega$ follows from the convexity of $X =\{(u_1,u_2, u_3) |
u_i=h(l_i), l_i+ l_j > l_k\}$ for $\lambda \leq -1$ in proposition
4.2(a). Thus by lemma 6.1, if $\lambda \leq -1$, $\bigtriangledown
W: \Omega \to \bold R^E$ is a smooth embedding.

\medskip

\medskip
\noindent 6.9. {\it A proof of theorem 6.2(a) in the case $\lambda
\geq 0$}

\medskip
To prove that the map $\Phi_{\lambda} : P_{\bold S^2}(S, T) \to
\bold R^E$ is an embedding for $\lambda \geq 0$, since
$\Phi_{\lambda}$ is known to be a local diffeomorphism, it
suffices to prove that $\Phi_{\lambda}$ is injective.

Recall that the space of all spherical triangles parameterized by
its inner angles is $\bold S^2(\theta, 3) =\{ (\theta_1, \theta_2,
\theta_3)| \theta_1+\theta_2+\theta_3 > \pi, \theta_i + \theta_j <
\theta_k+\pi, \{i,j,k\}=\{1,2,3\} \}$. Let $h(t) =\int^t_{\pi/2}
\sin^{\lambda}(s) ds$ and $u_i =h(\theta_i)$ for $i=1,2,3$. Then
by proposition 4.2(d) and $\lambda \geq 0$, the image of $\bold
S^2(\theta, 3)$ under the map $(u_1, u_2, u_3)=( u_1(\theta_1),
u_2(\theta_2), u_3(\theta_3))$ is an open convex set $X$ in $\bold
R^3$. By taking the Legendre transform of the closed 1-form
$w_{\lambda}$ in theorem 3.2(b), we obtain a closed 1-form
$$ w=\sum_{i=1}^3 (\int^{l_i}_{\pi/2} \sin^{-\lambda -1}(s)
ds)du_i$$ so that the integration $F(u) = \int^u_0 w$ is a locally
strictly convex function in $u$ defined on $X$. Furthermore, by
definition, $$\frac{\partial F}{\partial u_i} =\int^{l_i}_{\pi/2}
\sin^{-\lambda -1}(s) ds \tag 6.4$$

Let $C(S, T)$ be the set of all corners in the triangulated
surface $(S, T)$ (see \S5.1). A \it $\lambda$-angle structure \rm
on $(S, T)$ is a map $u: C(S, T) \to \bold R$ so that if $c_1,
c_2, c_3$ are three corners of a triangle in $T^{(2)}$ then
$(u(c_1), u(c_2), u(c_3)) \in X$. Geometrically, a $\lambda$-angle
structure can be identified with the realization of each 2-cell in
$T^{(2)}$ by a spherical triangle. However, these spherical
triangles may not have the same edge lengths at each edge of
$T^{(1)}$. A $\lambda$-angle structure $u$ is called \it geometric
\rm if the lengths of each edge in $T^{(1)}$ are the same in the
two spherical triangles adjacent to it. This is the same as the
following condition. There is a spherical polyhedral metric $l$ on
$(S, T)$ so that $u(c)= \int^{\theta}_{\pi/2} \sin ^{-\lambda
-1}(s) ds$ where $\theta$ is the inner angle at the corner $c$ in
the metric $l$.

The \it edge invariant \rm of a $\lambda$-angle structure $u$ is
the map from the set of all edges $E$  to $\bold R$ sending an
edge $e$ to $u(c)+ u(c')$ where $c, c'$ are the corners facing the
edge $e$. By definition, the  $\phi_{\lambda}$ invariant is the
edge invariant of the $\lambda$-angle structure associated to the
spherical polyhedral metric.

Let $A$ be the space of all $\lambda$-angle structures on $(S,
T)$. The space $A$ is an open convex subset of $\bold R^{C(S, T)}$
affine homeomorphic to $X^{T^{(2)}}$.  Fix  a map $\phi: E \to
\bold R$. Let $A(\phi)$ be the subspace of $A$ consisting of
$\lambda$-angle structures whose edge invariants are $\phi$. By
definition, $A(\phi)$ is convex.

\medskip
\noindent {\bf Claim 6.4.} \it For any $\phi: E \to \bold R$,
there is at most one geometric $\lambda$-angle structure in
$A(\phi)$. \rm

\medskip
As a consequence, the map $\Phi_{\lambda}: P_{\bold S^2}(S, T) \to
\bold R^E$ is injective.

To see the claim, let us label the set of all corners $C(S, T)$ by
$\{ c_1, ...., c_p\}$. The value of a $\lambda$-angle structure
$u$ at $c_i$ is denoted by $u_i$. Define  a function $W: A \to
\bold R$ by
$$ W(u) =\sum_{\{i,j,k\} \in T^{(2)}} F(u_i, u_j, u_k)  $$ where the
sum is over all triangles $\{i,j,k\}$ with three corners $c_i,
c_j, c_k$. By definition and the basic property of $F$,  the
function $W$ is strictly locally convex defined on the convex set
$A$. Thus its restriction $W|: A(\phi) \to \bold R$ is again a
strictly locally convex function defined on the convex set
$A(\phi)$. In particular, the function $W|$ has at most one
critical point in $A(\phi)$. On the other  hand, the critical
points of $W|$ are exactly the geometric $\lambda$-angle
structures by the lemma below. Assuming this lemma, we see that
the claim 6.4 follows.

\medskip
\noindent {\bf Lemma 6.5.} \it The critical  points of $W|:
A(\phi) \to \bold R$ are exactly equal to the geometric
$\lambda$-angle structures in $A(\phi)$. \rm

\medskip

\noindent {\bf Proof.} Suppose $u^*=(u_1^*, ..., u^*_p)$ is a
critical point of $W|$. By the Lagrangian multiplier's method,
there exists a function $\alpha: E \to \bold R$ (the multiplier)
so that
$$ \frac{ \partial W}{\partial u_i} (u^*) = \alpha(e)  \tag 6.5$$
when the i-th corner $c_i$ is facing the edge $e$. By (6.4),
equation (6.5) shows that
$$\int^{l}_{\pi/2} \sin^{-\lambda -1}(s) ds = \int^{l'}_{\pi/2}
\sin^{-\lambda -1}(s) ds  \tag 6.6$$ where $l$ and $l'$ are the
lengths of the edge $e$ in the two spherical triangles  adjacent
to $e$. But (6.6) implies $l=l'$. Thus at the critical point
$u^*$, we can glue all spherical triangles in the $\lambda$-angle
structure isometrically to produce a spherical polyhedral metric,
i.e., $u^*$ is geometric.

Conversely, since the constraint equations $u(c)+u(c') =\phi(e)$
are linear equations, the critical points of $W|$ are identical to
the solutions of equation (6.5). This shows that geometric
$\lambda$-angle structures are critical points.

It follows that $\Phi_{\lambda}$ is injective when $\lambda \geq
0$.

The above proof follows the same strategy used by Rivin in [Ri]
for characterization geometric structures among angles structures.

\medskip
\noindent 6.10. {\it A proof of theorem 6.2(f) for $\lambda \geq
0$}

\medskip

The proof is essentially the same as that of theorem 6.2(a) when
$\lambda \geq 0$ in subsection 6.9.  We replace the space $\bold
S^2(\theta, 3)$ by $\bold H^2(r+r'=\theta, 3)$, replace the energy
function by $F(u) = - \int^{u}_{0} \sum_{i=1}^3 (\int^{l_i}
\tanh^{-\lambda-1}(s/2) ds) du_i $ where $u_i=\int^{r_i}_1
\cos^{\lambda}(s) ds$ in theorem 3.4(e) (by replacing $\lambda$ by
$\lambda-1$). By proposition 4.2(b), the image $X$ of $\bold
H^2(r+r'=\theta, 3)$ under the map $u=(u_1,u_2, u_3)$ where $u_i
=\int^{r_i}_1 \cos^{\lambda}(s) ds$ is open convex in $\bold R^3$
for $\lambda \geq 0$. We use the same definition of angle
structures as in subsection 6.9 and the edge invariant in the same
way. Then due to the convexity of $X$ and strictly locally
convexity of $F$, the same claim 6.4 still holds in this case.
Thus the map $\Psi_{\lambda}$ is injective when $\lambda \geq 0$.

\medskip
%\noindent 6.11. {\bf Remark.} In the work of Thurston for $\bold
%E^2$ and $\bold H^2$ circle packing metrics on $(S, T)$, the set
%of all images of the curvature is explicitly described as open
%convex polytopes. Similar results for $\lambda \leq -1$ exist and
%will be discussed in  section 7.

\medskip
\noindent \S7. {\bf Parameterizations of the Teichm\"uller Space
of a Surface with Boundary}

\medskip
Suppose $S$ is a compact surface with negative Euler
characteristic and non-empty boundary. For each $\lambda \in \bold
R$,  we produce a parameterization $\psi_{\lambda}$ of the
Teichm\"uller space of the surface $S$ in this section. In the
case $\lambda=0$, this parameterization was first found in [Lu2].

Suppose $T$ is an ideal triangulation of $S$.
 Let $E=\{ e_1, ..., e_m\}$ be the
set of all edges in $T$. If $f: E \to X$ is a function, we use
$f_i$ to denote $f(e_i)$. By [Us] or the discussion in section 5,
equivalent classes of hyperbolic metrics with geodesic boundary on
$S$ are in one-one correspondence with the edge length functions
$l: E \to \bold R_{
>0}$. We  identify the
Teichm\"uller space $Teich(S)$ with $\bold R_{>0}^E$ by the edge
length function.

Let $g_{\lambda}(x) = \int^x_0 \cosh^{\lambda}(t) dt$. For a
hyperbolic metric $l$ on $S$, recall that the \it
$\psi_{\lambda}$-edge invariant \rm of the metric,
$\psi_{\lambda}: E \to \bold R$,
 is
defined to be
$$ \psi_{\lambda}(e) = g_{\lambda}( \frac{b+c-a}{2}) +
g_{\lambda}(\frac{b'+c'-a'}{2})  \tag 7.1$$ where $b,c,b'c'$ are
the lengths of the B-arcs adjacent to the edge $e$ and $a,a'$ are
the lengths of the B-arcs facing the edge $e$. See figure 1.1 (b).
Denote $\Psi_{\lambda}: \bold R_{>0}^E \to \bold R^E$ the map
sending a metric $l$ to its $\psi_{\lambda}$ edge invariant.

\medskip
\noindent {\bf Theorem 7.1.} \it Suppose $(S, T)$ is an ideal
triangulated surface. For any $\lambda \in \bold R$, the map
$\Psi_{\lambda}: Teich(S) \to \bold R^E$ is a smooth embedding.
\rm

\medskip
Given an ideal triangulated (or triangulated) surface $(S, T)$, an
\it edge cycle \rm is an edge loop in the 1-skeleton of the dual
cellular decomposition of $(S, T)$. To be more precise, an edge
cycle consists of edges $\{ e_{i_1}, ..., e_{i_r}\}$ in $E$ and
2-cells $\{H_{1}, ..., H_{r}\}$ in $T^{(2)}$ so that for all
indices $s$, $e_{i_{s}}$ and $e_{i_{s+1}}$ lie in the 2-cell
$H_{i}$ where $e_{i_{r+1}}=e_{i_1}$. Since in most of the cases,
there is at most one 2-cell adjacent to two edges $e, e'$ (except
a degree two vertex), we will simply use   $\{ e_{i_1}, ...,
e_{i_r}\}$ to denote the edge cycle. It is understood that the
2-cells $H_{i}$ are part of the definition of the edge cycle.

\medskip
\noindent {\bf Theorem 7.2.} \it Let $\lambda \geq 0$. For an
ideal triangulated surface $(S, T)$, $ \Psi_{\lambda}(Teich(S))
=\{ z \in \bold R^E | $ for each edge cycle $\{ e_{i_1}, ...,
e_{i_m}\}$, $\sum_{s=1}^m z(e_{i_s})
>0$\}.  Furthermore,  the image $
\Psi_{\lambda}(Teich(S))$ is an open convex polytope independent
of the parameter $\lambda \geq 0$.\rm

\medskip

For $\lambda <0$, similar result has been established by Ren Guo
[Gu1] (theorem 1.4). His proof also works for all $\lambda$ and
shows that $\Psi_{\lambda}(Teich(S))$ is an explicit bounded open
convex polytope so that $\Psi_{\lambda}(Teich(S)) \subset
\Psi_{\mu}(Teich(S))$ if $\lambda < \mu$.
%\medskip
%\noindent {\bf Theorem 7.3 (Guo).} \it Let $\lambda <0$. For an
%ideal triangulated surface $(S, T)$, $ \Psi_{\lambda}(Teich(S)) =\{ z
%\in \bold R^E | $ for each edge cycle $\{ e_{i_1}, ...,
%e_{i_m}\}$, $\sum_{s=1}^m z(e_{i_s})
%>0$,  and $z(e_i) < 2 \int_0^{\infty} \cosh^{\lambda}(t) dt$\}. In particular, the image is an open convex polytope
%in $\bold R^E$. \rm

\medskip
\noindent 7.1. {\it Remarks} 1. Theorem 7.2 was proved for
$\lambda =0$ in [Lu2].

2. Whether these new coordinates is related to the quantum
Teichm\"uller space ([FC], [Ka], [BL], [Te]) is an interesting
question.

3. An edge cycle $\{ e_{i_1}, ..., e_{i_m}\}$ is called \it
fundamental \rm if each edge appears at most twice. It is proved
in [Lu2] that the convex set $\{ z \in \bold R^E | $ for each edge
cycle $\{ e_{i_1}, ..., e_{i_m}\}$, $\sum_{s=1}^m z(e_{i_s})
>0$\}
is defined by a finite set of linear inequalities $\sum_{s=1}^m
z(e_{i_s}) >0$ where $\{e_{i_1}, ..., e_{i_m}\}$ is a fundamental
edge cycle. Thus, $\Psi_{\lambda}(Teich(S))$ is an open convex
polytope in $\bold R^E$.

\medskip
\noindent 7.2. {\it A proof of theorem 7.1}

\medskip
The proof of theorem 7.1 is a simple application of the strictly
convexity of the energy functions introduced in \S3. By theorem
3.4(f) (by replacing $\lambda$ by $-\lambda-1$) and Legendre
transformation, for a colored hyperbolic right-angled hexagon of
red edge lengths $l_i, l_j, l_k$ and opposite black edge lengths
$\theta_i , \theta_j, \theta_k$ where $\theta_i = r_j + r_k$, $i
\neq j \neq k \neq i$, the following 1-form
$$ \omega = \sum_{i=1}^3 \int_0^{r_i} \cosh^{\lambda}(t) dt d(
\int_1^{l_i} \tanh^{\lambda+1} (t/2) dt)$$ is closed.   Let $u_i
=\int_1^{l_i} \tanh^{\lambda+1}(t/2) dt$ and $u=(u_1, u_2, u_3)$.
Then
$$w=\sum_{i=1}^3 \int^{r_i}_0 \cosh^{\lambda}(t) dt du_i, $$ and the
integration $F(u) =\int^u_{(1,1,1)} w$ is strictly concave in $u
\in I^3$ where $I = h(\bold R_{>0})$ and $h(x) = \int_1^x
\tanh^{\lambda+1}(t/2) dt$. Furthermore, $\frac{\partial
F}{\partial u_i} = \int^{r_i}_0 \cosh^{\lambda}(t)dt$.

For a hyperbolic metric $l: E \to \bold R_{>0}$  on $(S, T)$, let
$u: E \to \bold R$ be $u(e) =\int_1^{l(e)} \tanh^{\lambda+1}(t/2)
dt$. Then the set of all possible values of $u$ forms the open
convex cube $I^E$. Define an energy function $W: I^E \to \bold R$
by
$$ W(u) =\sum_{ \{e_a, e_b, e_c\} \in T^{(2)}} F(u_a, u_b, u_c)$$ where
$u_i =u(e_i)$ and the sum is over all hexagons with edges $e_a,
e_b, e_c$. By definition, $W$ is smooth and strictly concave in
$I^E$. Furthermore, by the construction of $F$,
$$\frac{ \partial W}{\partial u_i} =\psi_{\lambda}(e_i)  $$
i.e., $\bigtriangledown W =\Psi_{\lambda}$. By lemma 6.1 on
gradient of strictly convex function, we conclude that the map
$\bigtriangledown W: I^E \to \bold R^E$ is a smooth embedding.
This proves theorem 7.1.

\medskip
\noindent 7.3. {\it Degenerations of hyperbolic hexagons}

\medskip
\noindent
 {\bf Lemma 7.4.} \it Suppose a hyperbolic right-angled
hexagon has three non-pairwise adjacent edge lengths $l_1, l_2,
l_3$ and opposite edge lengths $\theta_1, \theta_2, \theta_3$ so
that $\theta_i = r_j + r_k$, $\{i,j,k\}=\{1,2,3\}$. Then the
following hold.

(a)  $\lim_{\theta_i \to 0} l_j(\theta_1, \theta_2, \theta_3)
=\infty$ for $j \neq i$ so that the convergence is uniform in
$\theta=(\theta_1, \theta_2, \theta_3)$. To be more precise, for
any $M>0$, there is $\epsilon
>0$, so that if $ \theta_i < \epsilon$, then $l_j > M$ for all choices of $\theta_j, \theta_k$.

(b) $\lim_{l_i \to 0} |r_i(l_1, l_2, l_3)| =\infty$ so that the
convergence is uniform in $l$, i.e.,  for any $M>0$, there is
$\epsilon >0$ so that if $l_i < \epsilon$ then $|r_i| >M$ for all
choices of $l_j, l_k$.

(c) Suppose a sequence of hexagons satisfies that $|r_1|, |r_2|,
|r_3|$ are uniformly bounded. Then $\lim_{l_i \to \infty}
\theta_j(l) \theta_k(l) =0$ so that the convergence is uniform in
$l$. \rm

\medskip
%\vskip.1in

%\epsfxsize=3truein \centerline{\epsfbox{7.1.eps}}

\noindent {\bf Proof.} For (a), we use the cosine law that
$$\cosh(l_j) =\frac{ \cosh(\theta_j) +
\cosh(\theta_i)\cosh(\theta_k)}{\sinh(\theta_i) \sinh(\theta_k)}$$
$$ \geq \frac{ \cosh(\theta_i)
\cosh(\theta_k)}{\sinh(\theta_i)\sinh(\theta_k)}$$
$$\geq \coth(\theta_i).$$
Since $\lim_{\theta_i \to 0} \coth(\theta_i) =\infty$, it follows
from the above inequality that $\lim_{\theta_i \to 0} l_j =\infty$
and the convergence is uniform in $\theta$.

\medskip

For (b), we use the tangent law (2.10) for hexagons that
$$\tanh^2(l_i/2) =
\frac{\cosh(r_j)\cosh(r_k)}{\cosh(r_i)\cosh(r_i+r_j+r_k)}$$
$$= \frac{1}{ \cosh( r_i)( \cosh(r_i) ( 1+\tanh(r_j)
\tanh(r_k)) +\sinh(r_i) ( \tanh(r_j) +\tanh(r_k)))}$$
$$\geq \frac{1}{ \cosh(r_i) ( \cosh(r_i)( 1+1) + |\sinh(r_i)|(
1+1))}$$
$$\geq \frac{1}{4 \cosh^2(r_i)}.$$
It follows that $\cosh^2(r_i) \geq \frac{1}{4 \tanh^2(l_i/2)}$.
Thus part (b) follows and the convergence is uniform in $l$.

\medskip
For part (c), by the assumption that $|r_i|$'s are uniformly
bounded, it follows that $\theta_i = r_j + r_k$ are uniformly
bounded from above. Now the cosine law says that
$$ \cosh(l_i) =\frac{ \cosh(\theta_i) + \cosh(\theta_j)
\cosh(\theta_k)}{\sinh(\theta_j) \sinh(\theta_k)}$$
$$\leq \frac{C}{\sinh(\theta_j)\sinh(\theta_k)}$$
for some constant $C$. Thus $\sinh(\theta_j) \sinh(\theta_k) \leq
\frac{C}{\cosh(l_i)}$. Since $\theta_j$ and $\theta_k$ are
uniformly bounded from above, it follows that $\lim_{ l_i \to
\infty} \theta_j \theta_k=0$ and the convergence is uniform in
$l$. QED

\medskip
\noindent 7.4. {\it A proof of theorem 7.2.}

\medskip
\noindent Let $\Psi_{\lambda}: \bold R_{>0}^E \to \bold R^E$ be
the map sending a hyperbolic metric $l \in \bold R^E_{>0}$ to its
$\psi_{\lambda}$ edge invariant. Let $\Omega$ be the convex set
$\{ z \in \bold R^E |$ whenever $e_{n_1}, ..., e_{n_r}$ form an
edge cycle, $\sum_{j=1}^r z(e_{n_j}) >0$\}.
 First
$\Psi_{\lambda}(\bold R_{>0}^E) \subset \Omega$. Indeed, fix a
hyperbolic metric $l \in \bold R^E_{>0}$. For an edge cycle
$\{e_{n_1}, ..., e_{n_r}\}$, let $a_j$ be the length of the B-arc
adjacent to $e_{n_j}$ and $e_{n_{j+1}}$ in the hexagon $H_j$
containing both $e_{n_j}$ and $e_{n_{j+1}}$.  Denote the lengths
of B-arcs in $H_j$ facing $e_{n_j}$ and $e_{n_{j+1}}$ by $b_j$ and
$c_j$. Then by definition, the contribution to $\sum_{j=1}^r
\psi_{\lambda}(e_{n_j})$ from the B-arcs inside $H_j$ is
$$ \int_0^{\frac{a_j+b_j -c_j}{2}} \cosh^{\lambda}(t) dt +
\int_0^{\frac{a_j+c_j -b_j}{2}} \cosh^{\lambda}(t)dt. \tag 7.2$$
It is positive due to the following lemma.

\medskip
\noindent {\bf Lemma 7.5.} \it For $a, b \in \bold R$, $
\int^{a}_0 \cosh^{\lambda}(t) dt + \int^{b}_0 \cosh^{\lambda}(t)
dt >0$ if and only if $a + b > 0$.\rm

\medskip

Indeed, the function $f(x) =\int_0^x \cosh^{\lambda}(t) dt$ is
strictly increasing in $x$ and $f$ is odd, i.e., $f(-x)=-f(x)$.
Thus if $a+b>0$, i.e., $a> -b$, then
 $ \int^{a}_0 \cosh^{\lambda}(t) dt + \int^{b}_0 \cosh^{\lambda}(t)
dt >   \int^{-b}_0 \cosh^{\lambda}(t) dt + \int^{b}_0
\cosh^{\lambda}(t) dt =0$. Conversely, by the same argument, if
$a+b < 0$, then $ \int^{a}_0 \cosh^{\lambda}(t) dt + \int^{b}_0
\cosh^{\lambda}(t) dt < \int^{-b}_0 \cosh^{\lambda}(t) dt +
\int^{b}_0 \cosh^{\lambda}(t) dt =0$.

By (7.2) and lemma 7.5 with $a =\frac{1}{2}(a_j+b_j-c_j)$ and
$b=\frac{1}{2}(a_j+c_j-b_j)$ so that $a+b>0$, we see that (7.2) is
positive. Thus the total sum $\sum_{j=1}^r \psi_{\lambda}(e_{n_j})
 >0$. This shows that $\Psi_{\lambda}(\bold R^E_{>0}) \subset
 \Omega$.

By theorem 7.1, $\Psi_{\lambda}(\bold R^E_{>0})$ is open in
$\Omega$. We claim that $\Psi_{\lambda}(\bold R^E_{>0})$ is also
closed in $\Omega$. Since $\Omega$ is connected, it follows that
$\Psi_{\lambda}(\bold R^E_{>0}) =\Omega$.

To see the closeness, take a sequence $l^{(m)} \in \bold R^E_{>0}$
so that $\Psi_{\lambda}(l^{(m)})$ converges to a point $w \in
\Omega$. We claim that $l^{(m)}$ contains a subsequence converging
to a point $p \in \bold R^E_{>0}$. By taking subsequence, we may
assume that $\lim_{m \to \infty} l^{(m)} = p \in [0, \infty]^E$
and the lengths of each B-arc in the metric $l^{(m)}$ converge in
$[0, \pi]$. It remains to show that for each edge $e \in E$, $p(e)
\in (0, \infty)$.

Suppose otherwise that there is an edge $e \in E$ so that $ p(e)
\in \{0, \infty\}$. We will derive a contradiction below.

Recall that the \it r-coordinate \rm of a B-arc $x$ is
$\frac{1}{2}(\alpha + \beta -\gamma)$ where $\gamma$ is the length
of the B-arc $x$ and $\alpha$ and $\beta$ are the lengths of the
other two B-arcs in the hexagon containing $x$.  By definition,
the edge invariant $\psi_{\lambda}(e)$ is $\int^{u}_0
\cosh^{\lambda}(t)dt +\int^{u'}_0 \cosh^{\lambda}(t) dt$ where $u,
u'$ are the r-coordinates of the B-arcs facing the edge $e$.

\medskip
 We claim,

 \medskip
 \noindent
{\bf Claim 7.6.}  {\it The r-coordinates $u(x)$ of each $B$-arc
$x$ in the metrics $l^{(m)}$ are bounded.}
\medskip
 If otherwise, say
$|u(x)|$ tends to infinity in the metrics $l^{(m)}$ as $m$ tends
to infinity, we will derive a contradiction as follows. Let $x'$
be the B-arc so that $x, x'$ are both facing an edge $e$. By the
assumption that $w \in \Omega$, $\psi_{\lambda}(e) = \int^{u(x)}_0
\cosh^{\lambda}(t) dt + \int^{u(x')}_0 \cosh^{\lambda}(t) dt$ is
finite. Since $\lambda \geq 0$ implies $\int^{\infty}_0
\cosh^{\lambda}(t)dt = \infty$, we must have $u(x')$ tends to
infinity in the metrics $l^{(m)}$ so that $u(x)$ and $u(x')$ have
the different signs. Say $u(x)$ tends to $-\infty$.
 Let $r_1$ and $r_2$ be the r-coordinates of the other
two B-arcs in the hexagon $H$ which contains $x$. Since $r_i+u(x)
\geq 0$ for $i=1,2$ by the definition of r-coordinate, we obtain
$\lim_{m \to \infty} r_i = \infty$ for $i=1,2$ in metrics
$l^{(m)}$.

In summary, we obtain two rules governing the r-coordinates in the
metrics $l^{(m)}$.

\noindent {\it Rule I}. If a B-arc has r-coordinate converging to
$-\infty$, then the other two B-arcs in the same hexagon have
r-coordinates converging to $+\infty$.

\medskip
\noindent {\it Rule II}. If $a$ and $b$ are two different B-arcs
facing an edge so that the r-coordinate of $a$ converges to $\pm
\infty$, then the r-coordinate of $b$ converges to $\mp \infty$.

\medskip
\vskip.1in

\epsfxsize=4truein \centerline{\epsfbox{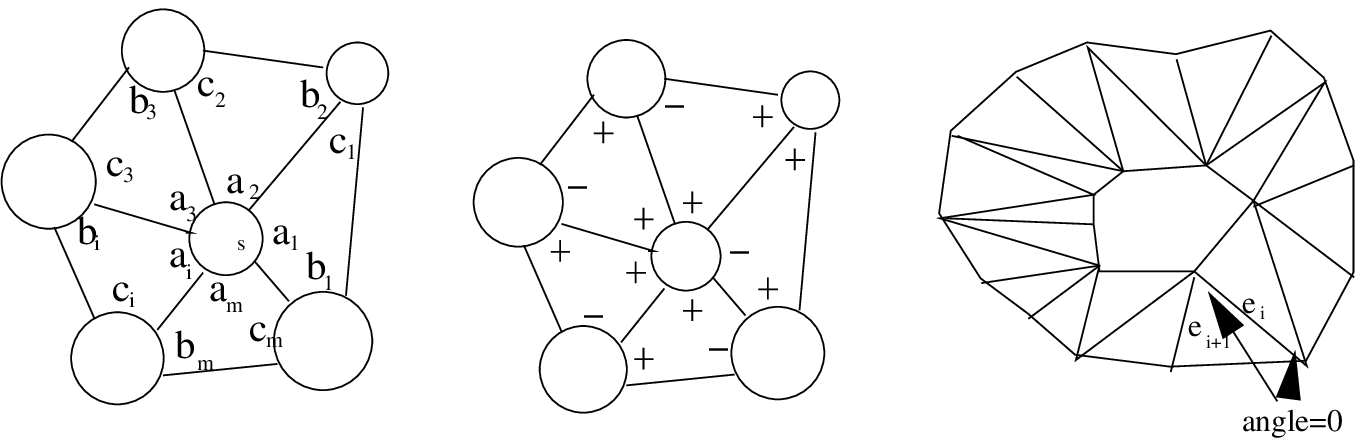}}

\centerline{(a)  $\quad \quad \quad \quad \quad \quad \quad \quad
\quad \quad$ (b)  $\quad  \quad \quad \quad \quad \quad \quad
\quad \quad$ (c)}

\centerline{Figure 7.1}

\medskip
 We claim that these two rules are contradicting to each
other on an ideal  triangulated surface $(S, T)$. Indeed, by the
assumption that $u(x)$ is unbounded and rule II, we find a B-arcc,
say $x$ itself, whose r-coordinate converges to $-\infty$. Let $s$
be the boundary component of the surface $S$ which contains the
B-arc $x$.
 Let
us label the edges ending at $s$ cyclically by, say,  $e_1, e_2,
..., e_k$ and the hexagon containing $e_i$ and $e_{i+1}$ be $H_i$
(with $e_{k+1}=e_1$). The r-coordinates of the B-arcs in $H_i$ are
denoted by $a_i$, $b_i$ and $c_i$ with $a_1=u(x)$ so that

(1) the B-arc of r-coordinate $a_i$ is in the boundary component
$s$ and is adjacent to both $e_i$ and $e_{i+1}$,

(2) the B-arc of r-coordinate $b_i$ is facing the edge $e_{i+1}$,

(3) the B-arc of r-coordinate $c_i$ is facing the edge $e_{i}$.

Then by the assumption, we have $\lim_{m \to \infty} a_1
=-\infty$. By rule I applied to $H_1$,  $\lim_{m \to \infty} b_1
=\lim_{m \to \infty} c_1 =\infty$. By rule II applied the two
B-arcs $b_1$ and $c_2$ facing the edge $e_2$, we have $\lim_{m \to
\infty} c_2 =-\infty$. By rule I applied to hexagon $H_2$, we have
$\lim_{m \to \infty} a_2 =\lim_{m \to \infty} b_2=\infty$.
Inductively, we obtain, for $i=2,..., k$, $\lim_{m \to \infty} a_i
=\lim_{m \to \infty} b_i =\infty$ and $\lim_{m \to \infty} c_i
=-\infty$. Finally, apply rule II to  two B-arcs with
r-coordinates $c_1$ and $b_{m}$ facing $e_1$ we obtain a
contradiction to rule II due to  $\lim_{m \to \infty} c_1 =\infty$
and $\lim_{m \to \infty} b_m =\infty$. See figure 7.1(a), (b).
This establishes claim 7.6.

By claim 7.6 and lemma 7.4(b), we see that $p(e) \in (0, \infty]$
for all edges $e$. Indeed, if an edge $e$ has $p(e) =0$, then by
lemma 7.4(b), the r-coordinates of the B-arcs adjacent to $e$
tends to infinity which is ruled out by claim 7.6 that
$r$-coordinates are bounded.

\medskip
It remains to show that $p(e) < \infty$ for all $e \in E$. Suppose
otherwise that $p(e) =\infty$ for some $e \in E$. By the above
claim, all r-coordinates of B-arcs in the metrics $l^{(m)}$ are
uniformly bounded. Let $H$ be a hexagon containing the edge $e$.
By the assumption $\lim_{m \to \infty} l^{(m)}(e) =\infty$ and
lemma 7.4 (c) applied to $H$, after taking a subsequence, the
length of one of the B-arcs adjacent to $e$ tends to 0. Say the
B-arc is $x$. Then by lemma 7.4(a) applied to $H$ with length of
$x$ tends to zero, the length of the other edge $e'$ adjacent to
$x$ tends to infinity. Let $H'$ be the hexagon adjacent to $H$
along the edge $e'$. We can then apply the same argument to $H'$.
In this way, we produce an edge cycle $\{e_1=e, ..., e_k\}$ so
that

(1) $\lim_{m \to \infty} l^{(m)}(e_i) =\infty$ for all $i$,

(2) $e_i$ and $e_{i+1}$ lie in a hexagon $H_i$ so that
$e_{k+1}=e_1$ for  $i=1,2,..., k$,

(3) the length $a_i^{(m)}$  of the B-arc $a_i$ in $H_i$ adjacent
to $e_i$ and $e_{i+1}$ converges to 0 (in the metrics $l^{(m)}$),

(4) the lengths of all B-arcs are bounded.

By definition, the sum of the $\psi_{\lambda}$ edge invariant at
$e_1, ..., e_k$ is
$$ \sum_{i=1}^k \psi_{\lambda}(e_i) = \sum_{i=1}^k (\int_0^ {\frac{a^{(m)}_i +b^{(m)}_i -b^{(m)}_{i+1}}{2}}
\cosh^{\lambda}(t) dt   + \int_0^ {\frac{a^{(m)}_i +b^{(m)}_{i+1}
-b^{(m)}_{i}}{2}} \cosh^{\lambda}(t) dt ). \tag 7.3$$

Since $\lim_{m \to \infty} a^{(m)}_i=0$, $b_i^{(m)}$ bounded and
$\cosh(-t)=\cosh(t)$, it follows that $$\lim_{m \to \infty}
\sum_{i=1}^k \psi_{\lambda}(e_i) =0.$$
 This
contradicts the assumption that $w \in \Omega$. Thus we conclude
that $\Psi_{\lambda}(\bold R^E_{>0})$ is closed in $\Omega$.

\medskip

\medskip

\noindent \S8. {\bf The Moduli Spaces of Polyhedral Surfaces, I,
Circle Packing Metrics}

\medskip
In Thurston's notes [Th], he gave a description of the spaces of
all discrete curvatures $k$ (and hence $k_0$) of Euclidean or
hyperbolic circle packing metrics on a triangulated surface and
showed that they are convex polytopes. (Also see [Ga],  [MiR],
[MaR], [CL]). The goal of this section is to give a description of
the spaces of all $k_{\lambda}$ curvatures when $\lambda \leq -1$.

\medskip
Let $(S, T)$ be a triangulated closed surface so that $E$ and $V$
are sets of all edges and vertices. Let $CP_{K^2}(S, T)$ be the
space of all circle packing metrics on $(S, T)$ in $K^2$ geometry
where $K^2=\bold S^2, \bold E^2$ or $\bold H^2$. Recall that
$K_{\lambda}: CP_{K^2}(S, T) \to \bold R^V$ sends a circle packing
metric to its $k_{\lambda}$-th discrete curvature.

\medskip
\noindent {\bf Theorem 8.1.} \it Suppose $\lambda \leq -1$.

(a) The space $K_{\lambda}(CP_{E^2}(S, T))$ is a proper
codimension-1 hypersurface in $\bold R^V$.

(b) The space $K_{\lambda}(CP_{H^2}(S, T))$ is an open set bounded
by $K_{\lambda}(CP_{E^2}(S, T))$ in $\bold R^V$. \rm

\medskip
\noindent 8.1. {\it Degeneration of Euclidean and hyperbolic
triangles}
\medskip
\noindent The following result on degeneration of triangles will
be used to analysis the singularities appeared in the variational
framework. Part of the lemma was proved already in [MaR] and [Th].

\medskip
\noindent {\bf Lemma 8.2.} (See also [Th], [MaR]). \it Suppose  a
Euclidean or hyperbolic triangle has edge lengths $l_1, l_2, l_3$
and angles $\theta_1, \theta_2, \theta_3$ where $\theta_i$ is
facing the edge of length $l_i$. Let $\{i,j,k\}=\{1,2,3\}$ and
$l_i = r_j+ r_k$.

(a) If the triangle is hyperbolic, then $\lim_{r_i \to \infty}
\theta_i(r_1, r_2, r_3)=0$ so that the convergence is uniform,
i.e., for any $\epsilon >0$, there is $M$ so that when $r_i>M$,
then  $\theta_i(r_1, r_2, r_3) < \epsilon$ for all choices of
$r_j$ and $r_k$.

(b) If the triangle is hyperbolic and  $l_i \to \infty$, then
after taking a subsequence, one of $l_j$ or $l_k$, say $l_j$,
tends to $\infty$, so that the angle $\theta_k$ between $l_i$-th
and $l_j$-th edges tends to zero.

(c) Suppose $r_i \geq c$ for a fixed constant $c>0$. Then
$\lim_{r_j \to 0} \theta_i(r_1, r_2, r_3) =0$ and the convergence
is uniform, i.e., for any $\epsilon >0$, there is $\delta =\delta(
\epsilon, c)>0$ so that for all $(r_1, r_2, r_3)$ with $r_j <
\delta$ and $r_i > c$, $\theta_i(r_1, r_2, r_3) < \epsilon$.

\rm \vskip.1in

\epsfxsize=3.5truein \centerline{\epsfbox{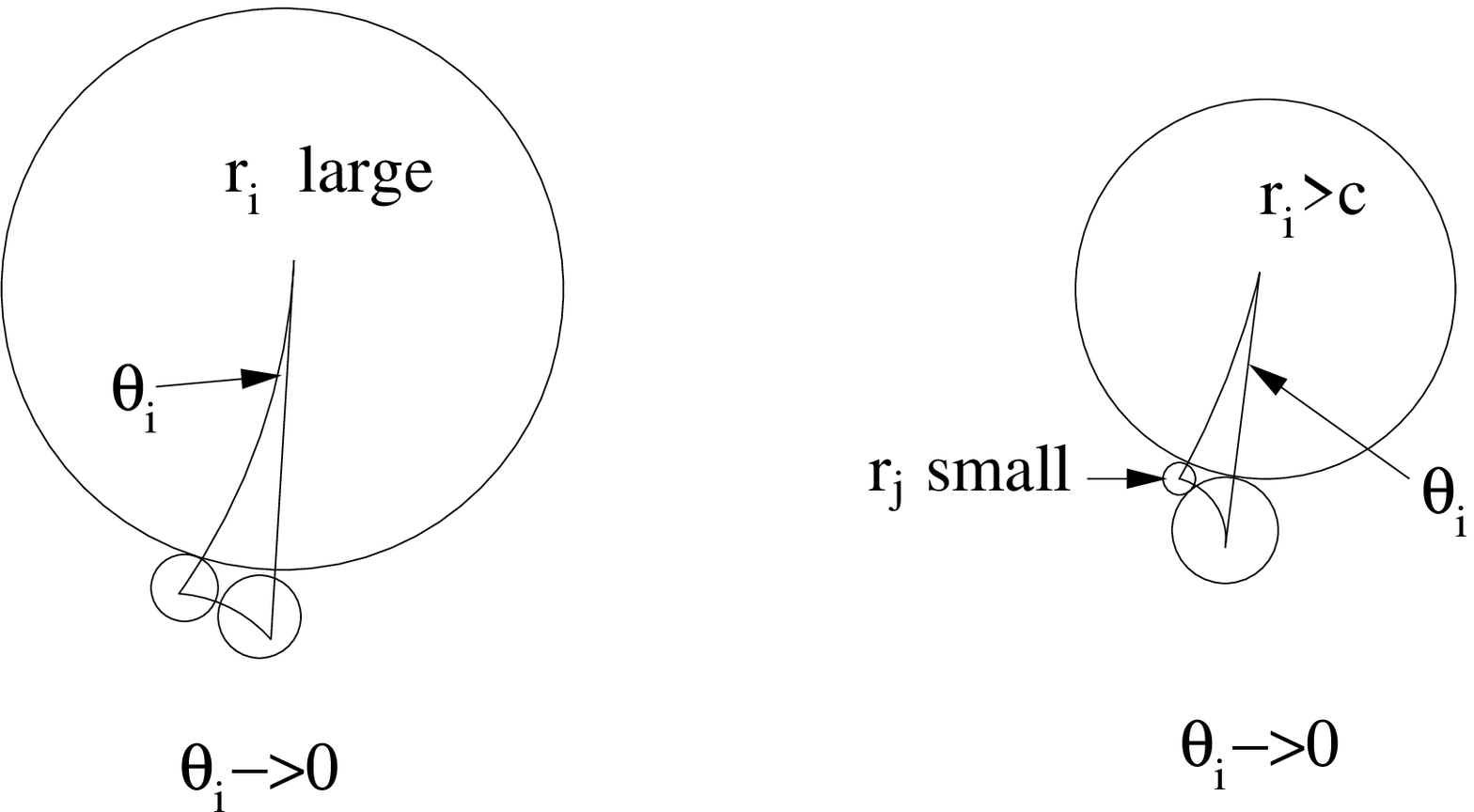}}

\medskip

 \centerline{Figure 8.1}

\medskip
\noindent {\bf Proof.} To see (a), recall the tangent law for
hyperbolic triangle (2.10) says,
$$ \tan^2(\theta_i/2) = \frac{\sinh(r_j) \sinh(r_k)}{\sinh(r_i)
\sinh(r_1+r_2+r_3)} \tag 8.1$$ Due to $\sinh(x+y) \geq 2
\sinh(x)\sinh(y)$ for $x, y >0$, it follows from (8.1) that
$$ \tan^2(\theta_i/2) \leq
\frac{\sinh(r_j) \sinh(r_k)}{\sinh(r_i) \sinh(r_j+r_k)}$$
$$\leq \frac{1}{2 \sinh(r_i)}.$$
Thus part (a) holds.

Another simple way to see part (a) is to put the i-th vertex to be
the Euclidean center of the Poincare disk model. The large radius
$r_i$ means the Euclidean diameter of the hyperbolic disk  $C$ of
radius $r_i$ centered at the origin is almost 1. This forces the
Euclidean diameter of any hyperbolic disk tangent to $C$  very
small. Thus the angle $\theta_i$ is very small no matter how one
chooses the radii $r_j$ and $r_k$.

\medskip
Part (b) follows from part (a). Indeed, since $l_i = r_j+r_k$ and
$l_i$ tends to infinity, one of $r_j$ or $r_k$ must tend to
infinity after taking a subsequence. Say $r_k$ tends to infinity.
Then due to $l_j \geq r_k$, $l_j$ converges to infinity. By part
(a), $\theta_k$ tends to 0.

\medskip
To see part (c) for hyperbolic triangles, using (8.1), we obtain
$$ \tan^2(\theta_i/2) \leq
\frac{\sinh(r_j)}{\sinh(r_i)} \frac{\sinh(r_k)}{ \sinh(r_j+r_k)}$$
$$\leq  \frac{\sinh(r_j)}{\sinh(c)}
\frac{\sinh(r_k)}{\sinh(r_j+r_k)}$$
$$\leq \frac{\sinh(r_j)}{\sinh(c)}.$$
Thus part (c) follows.

\medskip
To see part (c) for Euclidean triangles, recall that the radius of
the inscribed circle of a Euclidean triangle is $R =\sqrt{
\frac{r_1r_2r_3}{r_1+r_2+r_3}}$. Thus by $\tan(\theta_i/2)
=\frac{R}{r_i}$, we obtain,
$$\tan^2(\theta_i/2) = \frac{ r_j r_k}{r_i(r_1+r_2+r_3)}$$
$$\leq \frac{r_j}{c}$$
Thus we obtain the uniform convergence of $\theta_i$ to 0.

\medskip
\noindent 8.2. {\it A proof of theorem 8.1 (a)}
\medskip
We identify the space $CP_{E^2}(S, T)$ of all Euclidean circle
packing metrics with $\bold R^V_{>0}$ by the radius parameter. Let
$X =\{ r \in \bold R^V_{>0} | \sum_{v \in V} r(v) =1$\} be the
space of all normalized circle packing metrics. By definition and
theorem 6.2, $K_{\lambda}(CP_{E^2}(S, T)) = K_{\lambda}(X)$,
$K_{\lambda}|: X \to \bold R^V$ is an embedding and its image is a
codimension-1 smooth submanifold. It remains to show that when
$\lambda \leq -1$, $K_{\lambda}(X)$ is a closed subset of $\bold
R^V$. To this end, take a sequence of points $\{ r^{(m)}\}$ in $X$
so that $K_{\lambda}(r^{(m)})$ converges to a point $w \in \bold
R^V$. We will prove that $\{r^{(m)}\}$ contains a convergent
subsequence in $X$.

Since the space $X$ is bounded, by taking a subsequence if
necessary, we may assume that $r^{(m)}$ converges to a point $p$
in the closure $\bar{X}$ of $X$ and the inner angles of each
corner in metrics $r^{(m)}$ converge. If $p \in X$, we are done.
If otherwise, the set $I =\{ v \in V | p(v) =0\}$ is non-empty and
$I \neq V$ due to the normalization condition. Since the surface
$S$ is connected, there exists a triangle $\sigma \in T^{(2)}$
with vertices, say $v_1, v_2, v_3$, so that $p(v_2)=0$ and $p(v_1)
>0$.

We claim that $\lim_{ m \to \infty} k_{\lambda}(v_1) =-\infty$ in
the metrics $r^{(m)}$. This will contradict the assumption that
the limit is $w(v_1) \in \bold R$.

To see the claim, consider those triangles $\tau$ having $v_1$ as
a vertex. Let $\theta$ be the inner angle in triangle $\tau$ at a
vertex $v_1$. By definition,

$$ k_{\lambda}(v_1) =\sum_{\tau} \int^{\theta}_{\pi/2}
\tan^{\lambda}(t/2) dt.  \tag 8.2$$ where the sum is over all such
triangles.

We now analysis the angle $\theta$. If $v_j$ and $v_k$ are the
other two vertices of $\tau$, then there are two cases: (1) both
$p(v_j)$ and $p(r_k)$ are positive or (2) one of $p(v_j), p(v_k)$
is zero. In the case (1), the triangle $\tau$ is non-degenerated
since $p(v_1) >0$ and thus $\theta \in (0, \pi)$. The contribution
of $\int^{\theta}_{\pi/2} \tan^{\lambda}(t/2) dt$ to the sum (8.2)
is finite. In the case (2), say $p(r_j)=0$, by lemma 8.1(c), the
angle $\theta$ in the metrics $r^{(m)}$ converges to $0$  as $m$
tends to infinity. Thus the contribution of the term from $\tau$
to (8.2) is negative infinity (i.e., $\int^{0}_{\pi/2}
\tan^{\lambda}(t) dt =-\infty$, due to $\lambda \leq -1$). By the
choice of $v_1, v_2, v_3$, $p(v_2)=0$ and $p(v_1) >0$, it follows
that case (2) exists. This establishes the claim and hence the
proof of theorem 8.1(a).

\medskip
\noindent 8.3. {\it A proof of theorem 8.1(b)}

\medskip
We again identify $CP_{H^2}(S, T)$ with $\bold R^V_{>0}$ by the
radius parameter.  By theorem 6.2, the map $K_{\lambda}: \bold
R^V_{>0} \to \bold R^V$ is an embedding. The goal is to prove the
image $K_{\lambda}(\bold R^V_{>0})$ is an open region bounded by
$K_{\lambda}(CP_{E^2}(S, T))$, i.e., boundary points of
$K_{\lambda}(\bold R^V_{>0})$  are in
$K_{\lambda}(CP_{E^2}(S,T))$. To this end, take a sequence $\{
r^{(m)} \}$ converging to a boundary point $p \in [0, \infty]^V$
of $\bold R^V_{>0}$ so that $K_{\lambda}(r^{(m)})$ converges to a
point $w \in \bold R^V$.  We may assume, after taking a
subsequence,  that the inner angles of each corner in metrics
$r^{(m)}$ converge.  We will show that $w \in
K_{\lambda}(CP_{E^2}(S, T))$.

Since the point $p$ is in the boundary of $\bold R^V_{>0}$, there
are three possibilities: (1) there is a vertex $v$ so that $p(v) =
\infty$, (2)  $p(v) < \infty$ for all $v \in V$ and there are
$v_1, v_2 \in V$ so that $p(v_2) =0$ and $p(v_1)>0$, (3) $p(v)=0$
for all $v \in V$.

In the first case, say $p(v_1) =\infty$. Then by lemma 8.1(a), all
angles $\theta$  at vertex  $v_1$ converge to 0 uniformly. It
follows that the $\lambda$-th discrete curvature at $v_1$,
$k_{\lambda}(v_1) = \sum_{\theta} \int_{\pi/2}^{\theta}
\tan^{\lambda}(t/2) dt$ diverges to $-\infty$ due to $\lambda \leq
-1$, i.e., $w(v_1)$ is infinite, contradicting $w \in \bold R^V$.

In the case (2),  then exactly the same argument used in
subsection 8.2 works in this case due to that fact that lemma
8.2(c) holds for Euclidean and hyperbolic triangles. Thus we
conclude that $w$ is infinite contradicting $w \in \bold R^V$.

The only case left is that $p(v)=0$ for all $v \in V$. In this
case, the metric $r^{(m)}$ are degenerating to Euclidean circle
packing metrics after a scaling. By theorem 8.1(a) that
$K_{\lambda}(CP_{E^2}(S, T))$ is closed in $\bold R^V$, it follows
that either $w$ is infinite or $w$ is in $K_{\lambda}(CP_{E^2}(S,
T))$.

\medskip

\medskip
\noindent \S9. {\bf Moduli Spaces of Polyhedral Metrics, II,
General Cases}

\medskip
We give descriptions of the spaces of all $\phi_{\lambda}$ and
$\psi_{\lambda}$ edge invariants on a triangulated surface in the
case $\lambda \leq -1$. The results of Rivin and Leibon on the
spaces of all $\phi_0$ and $\psi_0$ invariants of Delaunay
polyhedral metrics will be revisited and reproved using different
methods.

One crucial step in the proofs below is to analyze degenerations
of geometric triangles. Recall that for $K^2=\bold H^2, \bold S^2$
or $\bold E^3$, the set $K^2(l, 3) \subset [0, \infty]^3$ denotes
the space of all triangles in $K^2$ parameterized by the edge
lengths $l=(l_1, l_2, l_3)$. A point $l \in \partial K^2(l, 3)
=\overline{K^2(l, 3)} - K^2(l, 3) \subset [0, \infty]^3$ is called
a \it degenerated triangle. \rm  The inner angles $\theta_1,
\theta_2, \theta_3$ of a degenerated triangle are not well defined
in general. Here is our convention of the inner angles. Take a
sequence $l^{(n)}$ in $K^2(l,3)$ converging to $l$ so that their
inner angles $(\theta^{(n)}_1, \theta^{(n)}_2, \theta^{(n)}_3)$
converge to $\theta=(\theta_1, \theta_2, \theta_3) \in [0,
\pi]^3$. Then we call $\theta_1, \theta_2, \theta_3$ the inner
angles of the degenerated triangle $l$. Note that $\theta_i$'s
depend on the choice of the converging sequences. However, in many
cases, even though each individual vector  $(\theta_i,\theta_j,
\theta_k)$ is not well defined, there are relations among their
entries which are valid for all choices of convergent sequences
$l^{(n)}$. For instance, $\theta_1+\theta_2+\theta_3 =\pi$ for all
degenerated hyperbolic triangles $l \in [0, \infty)^3$.

Recall that $\Phi_{\lambda}$ and $\Psi_{\lambda} : P_{K^2}(S, T)
\to \bold R^E$ are the maps sending a metric to its
$\phi_{\lambda}$ and $\psi_{\lambda}$ edge invariants. As a
convention, if $X$ is a subset of $R^n$, then $-X =\{ x \in \bold
R^N | -x \in X\}$.

We will prove, among other things, the following theorem.

\medskip
\noindent {\bf Theorem 9.1.} \it Suppose $(S, T)$ is a closed
triangulated surface so that $E=\{e_1, ..., e_n\}$ is the set of
all edges.  Let $\lambda \leq -1$.

(a) The space $\Phi_{\lambda}(P_{E^2}(S, T))$ is a proper smooth
codimension-1 submanifold in $\bold R^E$.

(b) The space $\Psi_{\lambda}(P_{H^2}(S, T))$ is an open set
bounded by $-\Phi_{\lambda}(P_{E^2}(S, T))$ and the following
linear inequalities:

$$ \sum_{ i=1}^k z(e_{n_i}) >0 \tag 9.1$$
where $z \in \bold R^E$ and $\{e_{n_1}, ..., e_{n_k}\}$ is an edge
cycle. \rm

\medskip
Note that a finite set of linear inequalities in (9.1) suffices.
Call an edge cycle $\{e_{n_1}, ..., e_{n_k}\}$ fundamental if
every edge appears at most twice. It is proved in [Lu3] that $\{
z\in \bold R^E | $(9.1) holds for all edge cycle\} is equal to $\{
z \in \bold R^E |$ (9.1) holds for all fundamental edge cycles\}.

\medskip
\noindent 9.1. {\it Degenerations of Euclidean triangles and
polyhedral metrics}

\medskip
Suppose $\sigma$ is a Euclidean triangle of edge lengths $l_1,
l_2, l_3$ so that their opposite angles are $\theta_1, \theta_2,
\theta_3$.  There are two cases that  Euclidean triangles $(l_1,
l_2, l_3)$ with $l_1+l_2+l_3=1$ degenerate:

(a1) one of $l_i=0$, and $l_j=l_k>0$ or,

(a2) $l_1, l_2, l_3
>0$ and $l_i=l_j+l_k$ for some $\{i,j,k\}=\{1,2,3\}$.

\medskip

\vskip.1in

\epsfxsize=3truein \centerline{\epsfbox{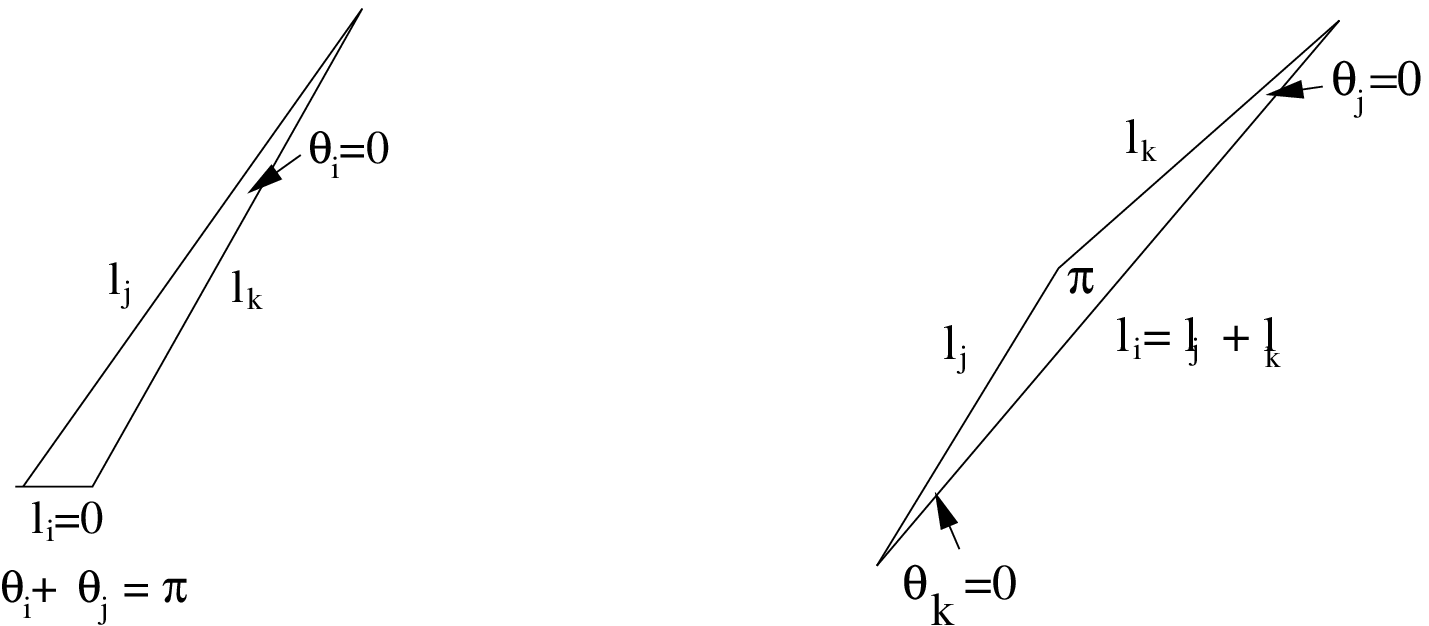}}

\medskip
\centerline{Figure 9.1}

\medskip
 In the case (a1), the angle $\theta_i=0$ is well defined.
 In the case (a2), $\theta_i=\pi$,
$\theta_j=\theta_k=0$. See figure 9.1. We call a degenerated
triangle in the case (a2) a \it $(\pi, 0, 0)$-angled triangle. \rm

A polyhedral metric $l \in P_{E^2}(S, T)$ is called \it normalized
\rm if $\sum_{ e \in E} l(e) =1$. Let $X$ be the set of all
normalized polyhedral metrics on $(S, T)$, i.e., $$X =P_{E^2}(S,
T) \cap \{ z \in \bold R^E | \sum_{ e\in E} l(e) =1\}. \tag 9.2$$
Due to $\phi_{\lambda}( l) = \phi_{\lambda}( cl)$ for $ c \in
\bold R_{>0}$, $\Phi_{\lambda}(X) =\Phi_{\lambda}(P_{E^2}(S, T))$.

The closure $\overline{X}$ of $X$ in $\bold R^E$ is compact. A
point $l \in \partial X =\overline{X}-X$ is called a \it
degenerated polyhedral metric \rm on $(S, T)$. There are two types
of degenerations: (1) there exists an edge $e$ so that $l(e)=0$,
or (2) $l(e) >0$ for all edges $e$ and there is a triangle
$\sigma$ with edges $e_i, e_j, e_k$ so that $$l(e_i) = l(e_j) +
l(e_k). \tag 9.3$$

In the case (1), let
 $\sigma$ be a triangle adjacent to an edge $e'$ of length zero so that one of the edge
 length of $\sigma$ is positive. Then the angle $\theta$ of
 $\sigma$ facing $e'$ is 0.
In the case (2), if $a,b,c$ are the angles in $\sigma$ facing
$e_i, e_j, e_k$ respectively, then $a=\pi$, $b=c=0$.
 Furthermore, the $\phi_{\lambda}$ edge
invariants are well defined on the case (2) degenerated metrics.
This shows that the map $\Phi_{\lambda}$ can be extended
continuously to $X \cup Y$ where $Y$ consists of all case (2)
degenerated polyhedral metrics.

 \medskip
 \noindent
 9.2. {\it A proof of theorem 9.1(a)}

 \medskip
 \noindent
By theorem 6.2, $\Phi_{\lambda}|: X \to \bold R^E$ is an embedding
and its image is a smooth codimension-1 submanifold. Thus it
suffices to show that $\Phi_{\lambda}(X)$ is a closed subset of
$\bold R^E$. To this end, take a sequence $\{l^{(m)}\}$ of points
in $X$ converging to a point $p \in \overline{X} -X$ so that the
angles of each corner in metrics $l^{(m)}$ in $T$ converge. In
particular, we may assume that  $\Phi_{\lambda}( l^{(m)})$
converges to a point $w$ in $[-\infty, \infty]^E$. We will show
that one of the coordinate of $w$ is infinite. Suppose otherwise
that $w \in \bold R^E$. We will derive a contradiction as follows.

Recall that the set of all edges $E =\{e_1, ..., e_n\}$.  By the
classification of degenerated polyhedral metrics in subsection 9.1
and normalization $\sum_{ e \in E} p(e) =1$, there are two cases:
(1) $p(e)>0$ for all edges $e$ and there is a triangle $\sigma$
with edges $e_i, e_j, e_k$ so that $p(e_i) = p(e_j) + p(e_k)$, or
(2) there is an edge $e$ so that $p(e) =0$.

In the argument below, all angles and lengths are measured in the
metric $p$.

 In the case (1), let the inner angles of the
triangle $\sigma$ be $\alpha, \beta, \gamma$ where $(\alpha,
\beta, \gamma)=(0, \pi, 0)$ so that $\alpha$ faces $e_j$. Now let
$\sigma'$ be the triangle adjacent to $\sigma$ along $e_j$ and
$\alpha'$ be the angle in $\sigma'$ facing $e_j$. Then by
definition,
$$\phi_{\lambda}(e_j) = \int^{\alpha}_{\pi/2} \sin^{\lambda}(t)
dt + \int_{\pi/2}^{\alpha'} \sin^{\lambda}(t) dt$$ is finite. Due
to the divergence of $\int^0_{\pi/2} \sin^{\lambda}(t) dt =
-\infty$ for $\lambda \leq -1$ and $\alpha=0$, it follows that
$\alpha'=\pi$. Thus, the inner angles of $\sigma'$ must be
$0,0,\pi$. To summary, we obtain the following rule: for $\lambda
\leq -1$, if $\alpha, \alpha'$ are two angles facing an edge so
that $\alpha=0$, then $\alpha'=\pi$. Now using this rule to
triangle $\sigma'$, we obtain a third $(0, 0, \pi)$-angled
triangle $\sigma''$ adjacent to $\sigma'$. Since there are only
finitely many triangles in $T$, by keep using this rule, we obtain
an edge cycle $\{ e_{n_1}, ..., e_{n_k}\}$ so that $e_{n_i},
e_{n_{i+1}}$ lie in a triangle $\sigma_i$ ($e_{n_{k+1}}=e_{n_1})$
and the angle of $\sigma_i$ facing $e_{n_i}$ is $\pi$. The inner
angles of $\sigma_i$ are $\pi, 0, 0$. We call such an edge cycle a
\it $(\pi, 0, 0)$-angled edge cycle\rm.

\medskip
\noindent {\bf Lemma 9.2.} \it There are no $(\pi, 0, 0)$-angled
edge cycles in a degenerated $K^2$ polyhedral metric for
$K^2=\bold E^2$, or $\bold H^2$, or $\bold S^2$. \rm

\medskip
\vskip.1in

\epsfxsize=4truein \centerline{\epsfbox{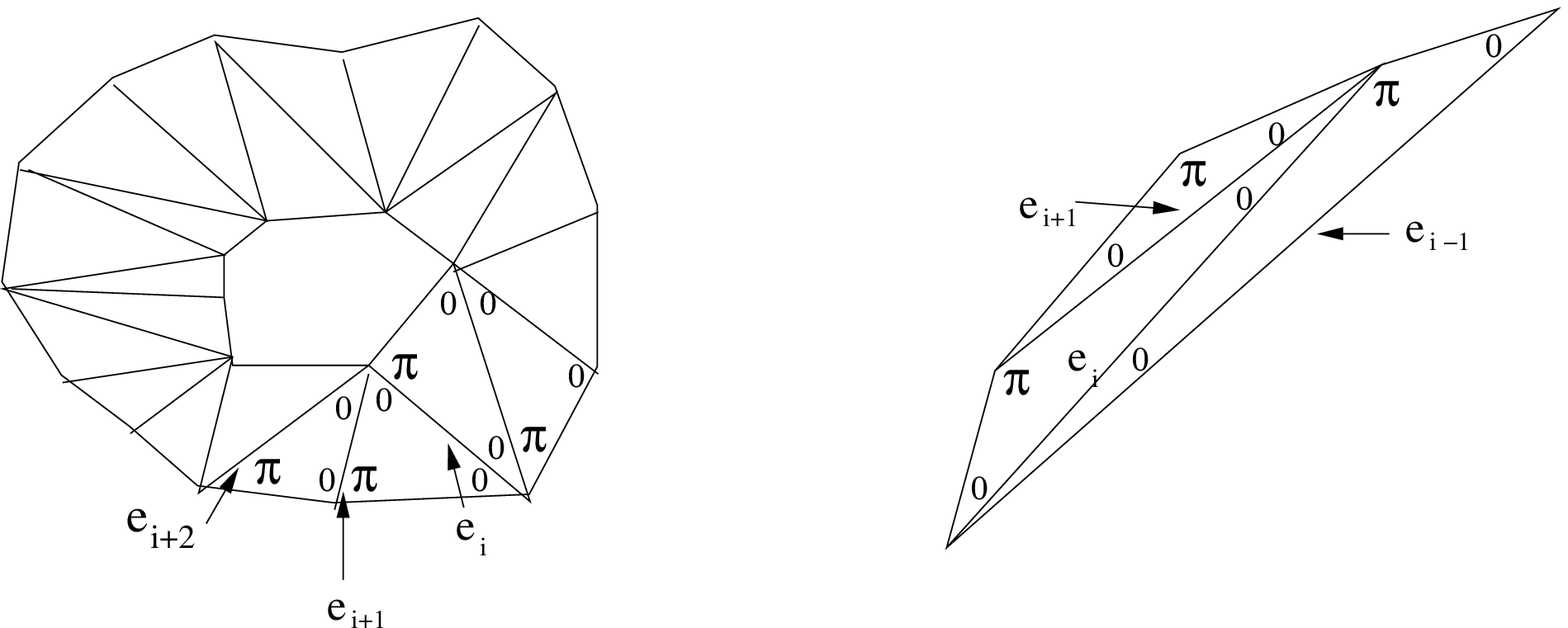}}

\centerline{Figure 9.2}\medskip

\medskip
\noindent {\bf Proof.} Suppose otherwise that such an edge cycle
exists. Take a sequence of non-degenerated polyhedral metrics
converging to the degenerated metric. We obtain a
(non-degenerated) polyhedral metric $l$ on $(S, T)$ so that the
inner angle of $\sigma_i$ facing $e_{n_i}$ is larger than the
other two angles in $\sigma_i$. Using the fact that in a Euclidean
(or hyperbolic or spherical) triangle, larger angle faces the edge
of longer length, we see that the length $l(e_{n_i})$ of $e_{n_i}$
is strictly larger than the length of $l(e_{n_{i+1}})$ of
$e_{n_{i+1}}$. Thus, we obtain
 $$l(e_{n_1}) > l(e_{n_2}) > ....> l(e_{n_k}) >
l(e_{n_{k+1}}) = l(e_{n_1}).$$ This is absurd. qed

\medskip
By this lemma, we conclude that case (1) does not occur.

\medskip
In the case (2) that some edge $e \in E$ has length $p(e) =0$,
 there must be some $e' \in E$ so that $p(e')>0$ due to the
normalization assumption $\sum_{x \in E} p(x) =1$.  It follows
that there is a triangle $\sigma$ with two edges $e, e'$ so that
$p(e)=0$ and $p(e')>0$. Thus the inner angle $\alpha$ of $\sigma$
facing $e$ must be 0. By the same argument as above, if $\alpha'$
is the other angle facing $e$, then $\alpha'=\pi$. This implies
that the triangle $\sigma'$ containing $\alpha'$ must have inner
angles $(\pi, 0, 0)$. By the same argument as above, we produce a
$(\pi, 0, 0)$-angled edge cycle. By lemma 9.2, this is impossible.
This ends the proof of theorem 9.1.

\medskip

\noindent 9.3. {\it Rivin's work revisited} \medskip

In a very influential paper [Ri], Rivin proved that the space
$\Phi_0(P_{E^2}(S, T)) \cap [-\pi, 0]^E$ of  $\phi_0$ invariants
of Euclidean Delaunay polyhedral metrics on a triangulated surface
 forms an explicit convex polytope. The goal of
this subsection is to extend his theorem slightly and give a
different proof of Rivin's theorem.

\medskip
\noindent {\bf Theorem 9.3.} \it Let $(S, T)$ be a triangulated
closed surface. The space $\Phi_0(P_{E^2}(S, T)) \subset \bold
R^E$ is in the affine plane $$\Cal A= \{ z \in \bold R^E | \sum_{
e \in E} z(e) = \pi( |T^{(2)}| - |T^{(1)}|) \}  \tag 9.4$$ bounded
by the following inequalities and codimension-1 submanifolds:

(a) For any proper subset $I \subset E$ so that no triangle has
exactly two edges in $I$,
$$ \sum_{ e \in I} z(e) > \pi |F_I| - \pi |I| \tag 9.5$$ where
$F_I =\{ \sigma \in T^{(2)} |$ all edges of $\sigma$ are in $I$\}.

(b) The hypersurfaces $\Cal W(i=j+k)$. Each of them is  the
$\Phi_0$ image of the codimension-1 submanifold $\{ z \in \bold
R^E_{>0} | z(e_i) = z(e_j)+z(e_k)$ where $e_i, e_j, e_k$ form the
edges of an triangle\}$ \cap \overline{ P_{E^2}(S, T)}$. \rm

\medskip
Numerical calculation shows that $\Phi_0(P_{E^2}(S,T))$ is not
convex in general. Metrics in  the hypersurface $\Cal W(i=j+k)$
are non-degenerated with respect to a different triangulation
obtained by
 the diagonal switch surgery operation on the $T$.

\medskip
\noindent {\bf Proof.} Identity (9.4) is the Gauss-Bonnet theorem
for Euclidean polyhedral surfaces. It follows that
$\Phi_0(P_{E^2}(S, T)) \subset \Cal A$. Let $X$ be the space of
all normalized polyhedral metrics defined by (9.2). We have
$\Phi_0(P_{E^2}(S, T))=\Phi_0(X)$ by definition.
 By Rivin's rigidity
theorem, the map $\Phi_0: X \to \Cal A$ is an embedding. It
follows that $\Phi_0(X)$ is open in $\Cal A$ by dimension
counting. To prove the theorem, we must analysis the boundary of
$\Phi_0(X)$. We will show that if $l^{(m)}$ is a sequence of
polyhedral metrics in $X$ converging to a boundary point $p \in
\overline{X} -X$, then $\Phi_0(l^{(m)})$ contains a subsequence
converging to a point either in $\Cal W(i=j+k)$ or in an affine
surface where one of the inequalities in condition (9.5) becomes
equality. Furthermore, we will prove that (9.5) holds for all
non-degenerated polyhedral metrics.

To this end, let us assume, after taking a subsequence that angles
of each corner in metrics $l^{(m)}$ in $T$ converge and
$\Phi_0(l^{(m)})$ converges to a point $w \in \bold R^E$. There
are two cases which could occur for the degenerated metric $p$ :
(1) $p(e)
>0$ for all $e \in E$ and there is a triangle with edges $e_i,
e_j, e_k$ so that $p(e_i) = p(e_j) + p(e_k)$, (2) the set $I =\{ e
\in E | p(e) =0\} \neq \emptyset$ and $I \neq E$. In the case (1),
we have $w \in \Cal W(i=j+k)$. In the case (2), by the triangular
inequality, there is no triangle $\sigma$ having exactly two edges
in $I$. Let $F_I$ be the set of all triangles with all edges in
$I$ and $G_I$ be the set of all triangles with exactly one edge in
$I$.
 By definition,
$$ \sum_{ e \in I} \phi_0(e) =  \sum_{
\sigma \in F_I} (a+b+c) + \sum_{ \sigma \in G_I} a- \pi |I| \tag
9.6$$ where the first sum is over triangles $\sigma$  in $F_I$
with inner angles $a,b,c$ and the second sum is over all triangles
$\sigma$ in $G_I$ with an inner angle $a$ facing an edge in $I$.
But the angle $a=0$ for triangles in $G_I$ by definition. It
follows from (9.6) that $\sum_{e \in I} \phi_0(e) = \pi
(|F_I|-|I|)$. This shows that the point $w$ is in an affine
surface defined by an equality from condition (9.5).

\medskip
The above argument also shows that condition (9.5) holds for
non-degenerated polyhedral metrics due to $G_I \neq \emptyset$,
$a>0$ and $\sum_{ \sigma \in G_I} a >0$ for non-degenerated
metrics. qed

\medskip
\noindent {\bf Corollary 9.4 (Rivin [Ri]).} \it The  space
$\Phi_0(P_{E^2}(S, T)) \cap [-\pi, 0]^E \subset \Cal A$ is the
convex polytope defined by condition (9.5) in theorem 9.3 and
inequalities $-\pi \leq z(e) \leq 0$ for all $e \in E$. \rm

\medskip
\noindent {\bf Proof.} We will use the same notations as above. It
suffices to show that the condition $\Cal W(i=j+k)$ does not arise
in the limits of $\phi_0$ edge invariants of Delaunay polyhedral
metrics $\Phi_0(P_{E^2}(S, T) )\cap [-\pi, 0]^E$. Suppose
otherwise that there is a sequence of metrics $\{ l^{(m)}\}$ in
$P_{E^2}(S, T) \cap \Phi_0^{-1}([-\pi, 0]^E)$ so that the sequence
$ l^{(m)}$ converges to a degenerated polyhedral metric $p \in
\Cal W(i=j+k)$. By definition, $p(e) \in (0, \infty)$ for all $e
\in E$ and there is a triangle $\sigma$ with edges
 $e_i, e_j, e_k$  so that $p(e_i)=p(e_j)+p(e_k)$.  Let the two inner angles facing the edge $e_i$ be $a$ and
$a'$ so that $a$ is in the triangle $\sigma$. Then $a=\pi$ and
inner angles of $\sigma$ are $\pi, 0, 0$. Buy definition
$\phi_0(e_i) =a+a' -\pi \leq 0$. It follows that $a'=0$. Since the
only degenerated triangles in $p$ are $(0, 0, \pi)$-angled
triangles, this implies the triangle $\sigma'$ adjacent to
$\sigma$ along $e_i$ must have inner angles $0, 0, \pi$. To
summary, we see that the Delaunay condition that $\phi_0(e) \in
[-\pi, 0]$ forces the propagation of $(0,0,\pi)$-angled triangles.
By keep using this propagation rule, we obtain a
$(0,0,\pi)$-angled edge cycle in the degenerated metric $p$. But
by lemma 9.2, this is impossible.

\medskip
\noindent 9.4. {\it Degeneration of hyperbolic polyhedral metrics}

\medskip

Let $l_1, l_2, l_3$ be the edge lengths of a hyperbolic triangle
so that the opposite angles are $\theta_1, \theta_2, \theta_3$. A
point in the boundary of $H^2(l, 3) =\{(l_1, l_2, l_3) | l_i + l_j
> l_k, \{i,j,k\}=\{1,2,3\} \}$ in $[0, \infty]^3$ is called a \it degenerated hyperbolic
triangle. \rm For degenerated triangles, the inner angles
$(\theta_1, \theta_2, \theta_3)$ depend on the choice of
convergence sequences. Let $\{i,j,k\}=\{1,2,3\}$. There are four
types of degenerated triangles:

(a1) some $l_i =\infty$;

(a2) all edge lengths $l_r$'s are finite so that some $l_i=0$ and
some $l_j>0$;

(a3) all $l_i=0$; and

(a4) all edge lengths are in $\bold R_{>0}$ and $l_i=l_j+l_k$ for
some $\{i,j,k\}=\{1,2,3\}$.

\medskip
\vskip.1in

\epsfxsize=4.5truein \centerline{\epsfbox{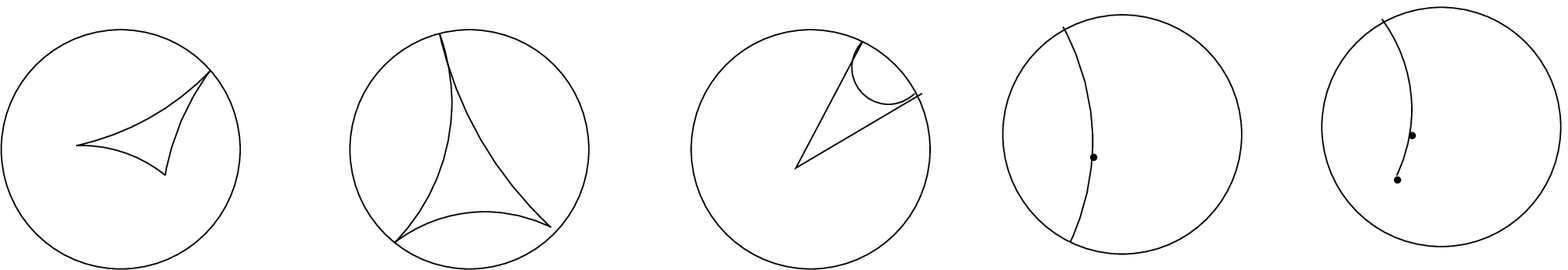}}

\medskip
\centerline{ Figure 9.3. degenerated hyperbolic triangles with one
infinite edge length}

\medskip
In the case (a1), due to the triangle inequality and lemma 8.2(b),
one of $l_j, l_k$, say $l_j =\infty$ so that $\theta_k=0$. In the
case (a2) that $l_i=0$ and $l_j=l_k >0$, we conclude that
$\theta_i=0$ and $\theta_j+\theta_k=\pi$. In the case (a3), we
have $\theta_i+\theta_j+\theta_k =\pi$.  In the case (a4),
$\theta_i=\pi$, $\theta_j =\theta_k =0$. See figure 9.3.

\medskip
Given a triangulated surface $(S, T)$, a  point $p$ in the
boundary of $P_{H^2}(S, T) \subset [0,\infty]^E$ is called a \it
degenerated hyperbolic metric. \rm  There are four types of
degenerated metrics: type I, $p(e) =\infty$ for some $e \in E$;
type II, $p(e) < \infty$ for all $e \in E$ so that  $p(e')=0$ and
$p(e'')
>0$ for some $e', e'' \in E$; type III, all $p(e) =0$; and type
 IV, $p(e) \in (0, \infty)$ for all $e
\in E$ and there is a triangle $\sigma$ with edges $e_i, e_j, e_k$
so that $p(e_i) =p(e_j)+p(e_k)$.

\medskip
\noindent {\bf Lemma 9.5.} \it (a).  In the type I degeneration
that $p(e)=\infty$ for some edge $e$, there exists an edge cycle
$\{ e_{n_1}, ...., e_{n_k}\}$ so that the lengths of $e_{n_i}$ are
infinite and the angle between $e_{n_i}, e_{n_{i+1}}$ in the
triangle $\sigma_i$ is 0.

(b). In the type II case, take a triangle $\sigma$ with two edges
$e, e'$ so that $p(e)=0$ and $p(e')>0$. Then  the angle facing $e$
in $\sigma$ is zero and the sum of all inner angles of $\sigma$ is
$\pi$.

(c). In the type III case, the sum of all inner angles of each
triangle is $\pi$.

(d). In the type IV case, the inner angles of $\sigma$ are $\pi,
0, 0$ so that the angle $\pi$ is facing edge $e_i$. \rm

\medskip

Indeed, the parts (b), (c), (d) of the lemma follow from the above
analysis of degenerations of hyperbolic triangles. In part (a),
take a triangle $\sigma$ adjacent to $e$.  Then lemma 8.2(c) shows
that there is another edge $e'$ in $\sigma$ so that $p(e')=\infty$
and the angle between $e,e'$ in $\sigma$ is 0. Now consider the
triangle $\sigma'$ adjacent to $\sigma$ along $e'$ and using lemma
8.2 (c) again. By keeping using lemma 8.2(c), we obtain an edge
cycle $\{ e_{n_1}, ...., e_{n_k}\}$ so that the lengths of
$e_{n_i}$ are infinite and the angle between $e_{n_i},
e_{n_{i+1}}$ in the triangle $\sigma_i$ is 0. We call it a \it
$(\infty, \infty, 0)$ edge cycle. \rm

\medskip
\noindent 9.5. {\it Proof of theorem 9.1(b)}

\medskip

By theorem 6.2, $\Psi_{\lambda}: P_{H^2}(S, T) \to \bold R^E$ is a
smooth embedding and its image is an open subset of $\bold R^E$.
To prove theorem 9.1(b), we must find the boundary points of
$\Psi_{\lambda}(P_{H^2}(S, T))$ in $\bold R^E$. Take a sequence of
points $\{ l^{(m)} \}$ converging to a boundary point $p$ of
$P_{H^2}(S, T)$ in $[0, \infty]^E$ so that angles of each corner
in metrics $l^{(m)}$ converge and $\Psi_{\lambda}( l^{(m)})$
converges to $w$ in $\bold R^E$. We will show that either $w$ is
in $-\Phi_{\lambda}(P_{E^2}(S, T))$ or in an affine surface
defined by the equality case of (9.1) for some edge cycle.
Furthermore, we will prove that (9.1) holds.  There are four types
of degenerations of $p$ according to subsection 9.4. If $p$ is  of
type  I that there is an edge $e$ so that $p(e) =\infty$, then by
lemma 9.5(a) there exists an edge cycle of type $(\infty, \infty,
0)$. Then by the same argument in proof of theorem 7.2 and
identity (7.3) where $\cosh(t)$ is replaced by $\cos(t)$ and
$\lim_{m \to \infty} a_i^{(m)}=0$, we conclude that $\sum_{i=1}^m
\psi_{\lambda}(e_{n_i}) =0$ along the edge cycle. Thus the point
$w$ is in the surface defined by an equality in (9.1). The proof
also shows that (9.1) holds for all edge cycles due to lemma 7.5
for $\cos(t)$ instead of $\cosh(t)$ in the region $t \in [-\pi/2,
\pi/2]$.

If $p$ is of type II that $p(e) <\infty$ for all $e \in E$ so that
there are two edges $e', e''$ with $p(e')=0$ and $p(e'')>0$, we
find a triangle $\sigma$ with three edges $e_i, e_j, e_k$ so that
$p(e_i)=0$ and $p(b_j) = p(e_k) >0$. Let the inner angles of
$\sigma$ be $a,b,c$ so that $a$ is facing $e_i$. Let $\sigma'$ be
the triangle adjacent to $\sigma$ along $e_i$ so that the inner
angles are $a', b', c'$ with $a'$
 facing $e_i$. Then by the choice of $\sigma$,
$a=0$ and $b+c=\pi$. On the other hand,
$$\psi_{\lambda}(e_i) =\int_0^{\frac{ b+c-a}{2}} \cos^{\lambda}(t) dt + \int_0^{\frac{
b'+c'-a' }{2}} \cos^{\lambda}(t) dt$$ By the assumption
$\frac{b+c-a}{2}$ is
 $\pi/2$. Due to the
divergence of $\int^{\pi/2}_0 \cos^{\lambda}(t) dt$ for $\lambda
\leq -1$ and the assumption that $\psi_{\lambda}(e_i)$ is finite,
we must have $b'+c'-a' =-\pi$. By the assumption that $a', b', c'
\geq 0$ and $a'+b'+c' \leq \pi$, we must have $(a',b',c') = (\pi,
0, 0)$. Now by the same argument applied to $b'=0$, we produce a
new $(0, 0, \pi)$ angled triangle adjacent to $\sigma'$. In this
way, we obtain a $(0,0,\pi)$-angled edge cycle in the
triangulation. By lemma 9.2,  this is impossible, i.e., type II
degenerated metric $p$ does not occur.

If $p$ is of type III that all edge lengths are zero, then each
triangle degenerates to a Euclidean triangle. Evidently
 if $a+b+c=\pi$, then $(a+b-c)/2 =\pi/2 -c$. Thus
$\int^{(a+b-c)/2}_0 \cos^{\lambda}(t) dt = \int^{\pi/2 -c}_0
\cos^{\lambda}(t) dt =- \int^{c}_{\pi/2}\sin^{\lambda}(t) dt$.
Thus $\psi_{\lambda}(e) = -\phi_{\lambda}(e)$. It follows that $w$
is in  the image $-\Phi_{\lambda}(P_{E^2}(S, T))$.

In the type IV degeneration that all $p(e) >0$ and there is a
triangle $\sigma$ with edges $e_i, e_j, e_k$ so that $p(e_i) =
p(e_j) + p(e_k)$. Then the inner angles of $\sigma$ are $\pi, 0,
0$. By the same argument as in the type II degeneration, due to
$\lambda \leq -1$, we see that the triangle adjacent to $\sigma$
along $e_i$ (also, $e_j, e_k$) has inner angles $0, 0, \pi$. It
follows that there must be a $(0, 0, \pi)$ edge cycles in $p$.
This again contradicts lemma 9.2.

\medskip
\noindent 9.6. {\it Leibon's work revisited}

\medskip
Leibon proved in [Le] that the space $\Psi_0(P_{H^2}(S, T) \cap
(0, \pi)^E$  of all $\psi_0$ invariants of Delaunay hyperbolic
polyhedral metrics is a convex polytope. We establish a
generalization of Leibon's theorem.

\medskip

\noindent {\bf Theorem 9.6}. \it The space $\Psi_0(P_{H^2}(S, T))$
is an open set in $\bold R^E$ bounded by the following set of
inequalities and hypersurfaces. Let $z \in \bold R^E$.

(a) For each edge cycle $\{e_{n_1}, ..., e_{n_k}\}$, $\sum_{i=1}^k
z(e_{n_i}) >0$.

(b) For any subset $I$ of $E$ with the property that no triangle
has exactly two edges in $I$, let $F'_I$ be the set of all
triangles having at least one edge in $I$, then

$$\sum_{e \in I} z(e) < \frac{\pi |F'_I|}{2}$$

(c) The hypersurface $\Cal W(i=j+k)$ which is the image under
$\Psi_0$ of the codimension-1 submanifold $ \{z \in \bold R^E_{>0}
| z(e_i)=z(e_j) + z(e_k)$ where $e_i, e_j, e_k$ are the edges of a
triangle\} $\cap \overline{ P_{H^2}(S, T)}$. \rm

\medskip
Points in $\Cal W(i=j+k)$ are non-degenerated polyhedral metrics
in a new triangulation obtained by the diagonal switch surgery
operation on $T$.

\medskip
\noindent {\bf Proof.} The proof follows the same argument used in
the proof of theorem 9.1(b). We will use the same notations as in
subsection 9.5. First by Leibon's rigidity theorem, $\Psi_0$ is a
smooth embedding. It follows that $\Psi_0(P_{H^2}(S, T))$ is an
open connected set in $\bold R^E$. We need to determine its
boundary.  Take a sequence of points $\{ l^{(m)} \}$ converging to
a boundary point $p$ of $P_{H^2}(S, T)$ in $[0, \infty]^E$ so that
the angles of each corner in the metric $l^{(m)}$ converge and
$\Psi_{\lambda}( l^{(m)})$ converges to $w$ in $\bold R^E$. We
will show that either $w$ lies in a surface defined by the
equality cases of conditions (a) or (b) or $w$ is in the
hypersurface $\Cal W(i=j+k)$. Furthermore, we prove that $(a)$ and
$(b)$ hold.

There are four types of degenerated metrics $p$ as shown in
subsection 9.4. All edge lengths and angles are measured in the
degenerated metric $p$ below.

In the type I case that $p(e) =\infty$ for some edge $e$, then we
obtain a $\{\infty, \infty, 0 \}$ edge cycle $\{ e_{n_1}, ...,
e_{n_m} \}$ according to lemma 9.5(a). By the definition of
$\psi_0$, the summation $\sum_{i=1}^m \psi_0(e_{n_i})$ is equal to
the summation $\sum_{i=1}^m a_i$ where $a_i$ is the angle between
$e_{n_i}$ and $e_{n_{i+1}}$ in the triangle containing both edges.
By the choice of the $\{0, 0, \infty  \}$ edge cycle, $a_i=0$.
Thus $\sum_{i=1}^m \psi_0( e_{n_i}) =0$. This shows that $p$ is in
the surface defined by the equality case of condition (a) for some
edge cycle. It also shows that condition (a) holds for all
hyperbolic polyhedral metrics in $P_{H^2}(S, T)$ since $a_i >0$
for non-degenerated triangles.

In the type II and III cases that $p(e) < \infty$ and some
$p(e')=0$, let  $I=\{ e \in E | p(e)=0\}$. By the triangular
inequalities, there is no triangle with exactly two edges in $I$.
Take a triangle $\sigma$ with inner angles $\theta_i, \theta_j,
\theta_k$ so that one of the edge of $\sigma$ is in $I$. If all
edges of the triangle are in $I$, then the sum
$\theta_i+\theta_j+\theta_k=\pi$. In this case the sum
$$(\theta_j+\theta_k-\theta_i)/2 + (\theta_i+\theta_k-\theta_j)/2 +
(\theta_i+\theta_j-\theta_k)/2 = \pi/2 \tag 9.7$$

If  only one edge of $\sigma$ is in I, say the edge $e_i$ facing
$\theta_i$ is in $I$, then by the assumption $p(e_i)=0$ and
$p(e_j) = p(e_k)
>0$. It follows that $\theta_i=0$ and $\theta_j + \theta_k =\pi$.
Thus
$$(\theta_j+\theta_k-\theta_i)/2 =\pi/2 \tag 9.8$$

Now the summation $\sum_{ e \in I} \psi_0(e)$ can be expressed as
$$ \sum_{ \sigma; e_i, e_j, e_k \in I}((\theta_j+\theta_k-\theta_i)/2 + (\theta_i+\theta_k-\theta_j)/2 +
(\theta_i+\theta_j-\theta_k)/2) +\sum_{ \sigma; e_i \in I, e_j
\notin I} (\theta_j+\theta_k-\theta_i)/2,$$ where first part
consists of sum over all triangles $\sigma$ whose three edges are
in $I$ and the second part consists of sum over all triangles with
exactly one edge $e_i$ in $I$. In the first part, the contribution
of $z(e)$'s from each triangle is $\pi/2$ due to (9.7). In the
second part and the contributions of $z(e)$ from each triangle is
again $\pi/2$ due to (9.8). It follows that $\sum_ {e\in I}
\psi_0(e) = \pi |F'_I|/2$, i.e., $p$ lies in a surface defined by
the equality case of condition (b). This also shows that the
inequality in condition (b) holds for all metrics in $P_{H^2}(S,
T)$ since (9.7) and (9.8) become strictly less than $\pi/2$ for
non-degenerated hyperbolic triangles.

In the type IV case,  by definition, $w \in \Cal W(i=j+k)$. qed
\medskip

\medskip
\noindent {\bf Corollary 9.7. (Leibon)} \it The space
$\Psi_0(P_{H^2}(S, T)) \cap (0, \pi]^E$ is a convex polytope
defined by condition (b) in theorem 9.6. \rm

\medskip
\noindent {\bf Proof.} It suffices to show that for the Delaunay
condition that $\phi_0(e) \in (0, \pi]$, both constraints (a) and
(c) are not necessary.

First of all, we show that condition (c) $\Cal W(i=j+k)$ does not
arise in the limits of Delaunay polyhedral metrics. Suppose
otherwise that there is a sequence of metrics $\{ l^{(m)}\}$
converging to $p$ in $P_{H^2}(S, T) \cap \Psi_0^{-1}((0, \pi]^E)$
so that the angles of each corner in metrics $l^{(m)}$ converge
and the sequence $ \Psi_0(l^{(m)})$ converges to a point $w \in
\Cal W(i=j+k)$. In the degenerated metric $p$, let $a,b,c$ be the
inner angles in the triangle $\sigma$ facing the edges $e_i, e_j,
e_k$ and $a',b', c'$ be angles of the triangle $\sigma'$ adjacent
to $\sigma$ along $e_i$ so that $a,a'$ are facing $e_i$. Then
$(a,b,c)=(\pi,0,0)$ and $(b+c-a)/2=-\pi/2$. Since $\psi_0(e_i)
=\frac{1}{2}(b+c-a+b'+c'-a') \geq 0$ and $(b'+c'-a')\leq \pi$, it
follows that $(b'+c'-a')/2=\pi/2$. This in turn implies that
$\{a',b',c'\} =\{ 0, 0, \pi \}$ with $a'=0$. In summary,  the
Delaunay condition that $\psi_0(e) \in [0, \pi]$ forces the
propagation of $(0,0,\pi)$ angled triangles. By keep using this
propagation rule, we construct a $(0,0,\pi)$ angled edge cycle in
the degenerated metric $p$. But by lemma 9.2, this is impossible.

Finally, it is clear that condition (a) follows from Delaunay
condition that $\psi_0(e) >0$.

\medskip

\noindent 9.7. {\it The moduli space of spherical polyhedral
surfaces}

\medskip

In this section we investigate the space of all spherical
polyhedral metrics on  $(S, T)$ in terms of the $\phi_{\lambda}$
edge invariant where $\lambda \leq -1$ or $\lambda=0$.

\medskip

We begin with a discussion of the degenerated spherical triangles.
Recall that the space of all spherical triangles in the edge
length parameterization is $\bold S^2(l, 3) =\{ (l_1, l_2, l_3)
\in \bold R^3 | l_i + l_j > l_k, $ and $ l_1 + l_2 + l_3 < 2\pi$
where $\{i,j,k\}=\{1,2,3\}$\}. It is an open set in $[0, \pi]^3$.
A point $l=(l_1, l_2, l_3)$ in the boundary $\partial \bold S^2(l,
3) \in \bold R^3$ is called a \it degenerated spherical triangle
\rm of edge lengths $l_1, l_2, l_3$. Let $\theta_1, \theta_2,
\theta_3$ be inner angles of $l$ (the vector $(\theta_1, \theta_2,
\theta_3)$ depends on the choice of convergent sequences). Since
the closure $\overline{ \bold S^2(l, 3)} $ is defined by the
inequalities $l_i + l_j \geq l_k$ and $l_1+l_2+l_3 \leq 2 \pi$, it
follows that if $l_i=0$ then $l_j=l_k$ and if $l_i=\pi$ then
$l_j+l_k=\pi$. We classify degenerated spherical triangles into
six types (assume $\{i,j,k\}=\{1,2,3\}$ below):

(a1) $l=(0,0,0)$. In this case $\theta_i+\theta_j+\theta_k=\pi$;

(a2) $l_i =0$ and $l_j=l_k=\pi$. In this case,
$\theta_i=\theta_j+\theta_k-\pi$;

(a3) $l_i=0$ and $l_j=l_k \in (0, \pi)$. In this case,
$\theta_i=0$ and $\theta_j+\theta_k=\pi$;

(a4) $l_i=\pi$ and $l_j+l_k =\pi$ so that $l_j, l_k \in (0, \pi)$.
In this case $\theta_i=\pi$ and $\theta_j=\theta_k$;

(a5)  $(l_1, l_2, l_3) \in (0, \pi)^3$ and $l_i = l_j + l_k$ for
some $i,j,k$. In this case $\theta_i=\pi$ and
$\theta_j=\theta_k=0$;

(a6) $(l_1, l_2, l_3) \in (0, \pi)^3$ and $l_1+l_2+ l_3 =2\pi$. In
this case all $\theta_i=\pi$.

\medskip
\vskip.1in

\epsfxsize=4.5truein \centerline{\epsfbox{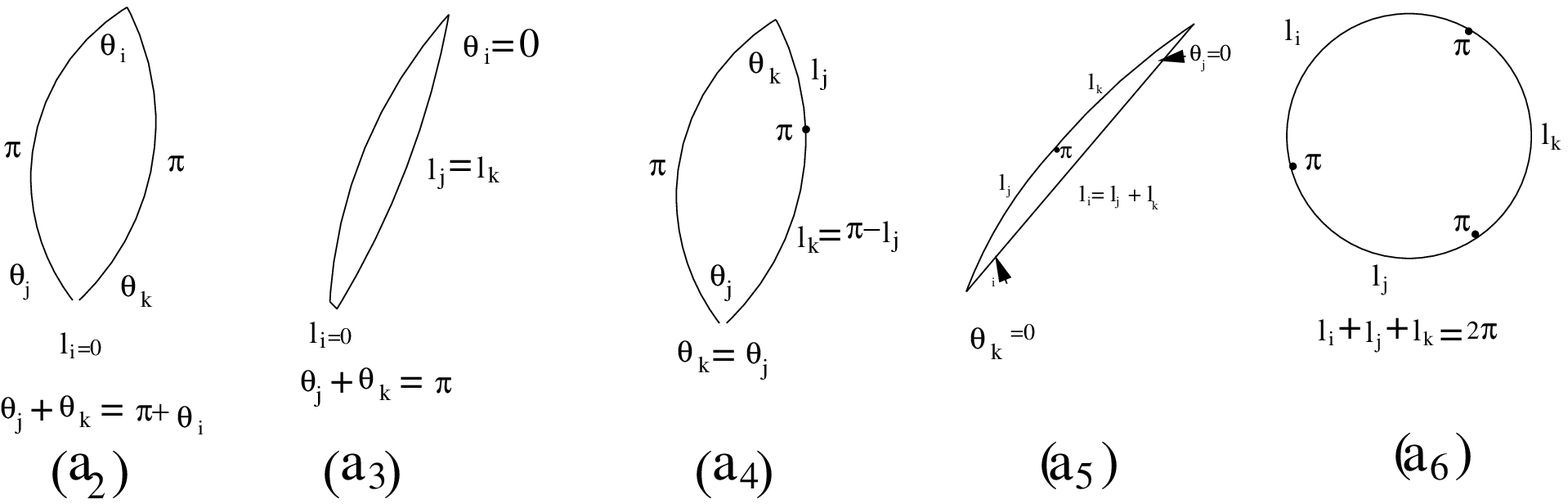}}

\medskip
\centerline{Figure 9.4}
 Note that in the last two cases (a5) and
(a6), each inner angle $\theta_i$ is well defined.

A \it degenerated \rm spherical polyhedral metric $l$ on a
triangulated surface $(S, T)$ is a point in the boundary of
$P_{S^2}(S, T) \subset \bold R^E$.  A degenerated spherical
polyhedral metric is called a \it bubble \rm if all triangles in
the metric are of types (a1) and (a2). Since a type (a2) triangle
is represented by a region in the 2-sphere bounded by two
geodesics of length $\pi$, i.e.,  a \it secant\rm,  geometrically
a bubble polyhedral surface is obtained by taking a finite set
(may be empty)  of secants and points and identify edges in pairs
and identify vertices. See figure 9.5.
\medskip
\epsfxsize=2.5truein \centerline{\epsfbox{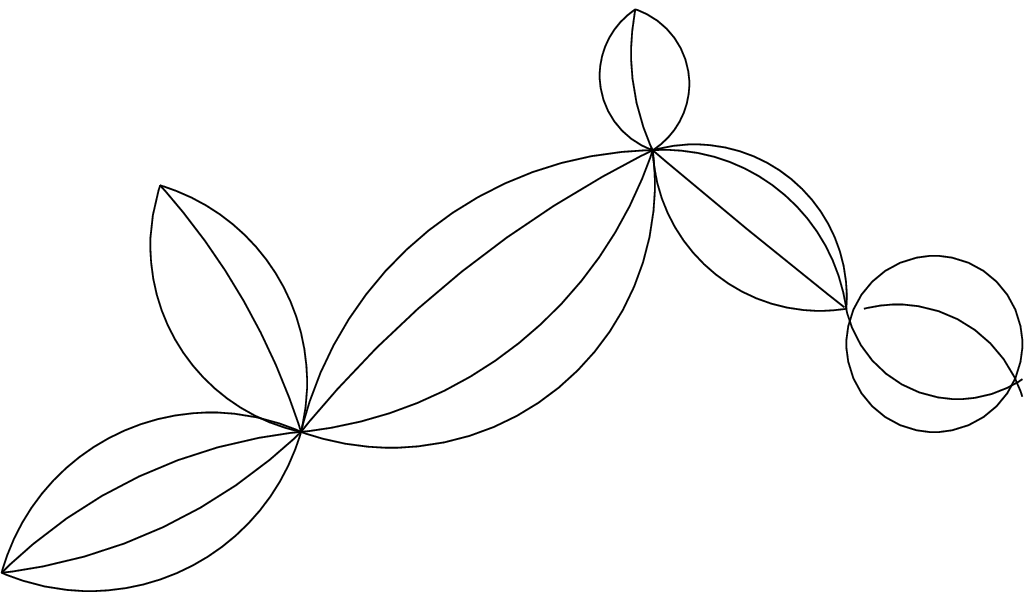}}

\medskip
\centerline{Figure 9.5}

A degenerated spherical polyhedral surface is called \it removable
\rm if all simplices in the metric are either in $\bold S^2(l, 3)$
or of types (a5) or (a6).  For a removable degenerated metric, the
curvature at each point is well defined and the metric becomes a
non-degenerated polyhedral metric in a different triangulation of
the surface.  By the discussion above, the $\phi_{0}$ edge
invariant is well defined on types (a5) and (a6) triangles. Let
$\Cal W(i=j+k)$ be the image $\Phi_{0}$(Y) where $Y=\{ z \in
(0,\pi)^E | z(e_i)=z(e_j)+ z(e_k)$ where $e_i, e_j, e_k$ form the
edges of a triangle $\} \cap \overline{ \bold P_{S^2}(S, T)}$.
Similarly, let $\Cal W(i+j+k)$ be the image $\Phi_{0}$(Z) where
$Z= \{ z \in (0,\pi)^E | z(e_i)+z(e_j)+ z(e_k)=2\pi$ where $e_i,
e_j, e_k$ form the edges of a triangle\} $\cap \overline{ \bold
P_{S^2}(S, T)}$. Both of them are codimension-1 hypersurfaces in
$\bold R^E$.

\medskip
\noindent {\bf Theorem 9.8.} \it Suppose $(S, T)$ is a closed
triangulated surface so that $E$ is the set of all edges in the
triangulation.

(a) Let $\lambda \leq -1$. The space $\Phi_{\lambda}(P_{S^2}(S,
T))$ of all $\phi_{\lambda}$ edge invariants of spherical
polyhedral metrics on $(S, T)$ is the open set in $\bold R^E$
whose boundary consists of images under $\Phi_{\lambda}$ of the
bubble degenerated spherical surfaces.

(b) The space $\Phi_0(P_{S^2}(S, T))$ of all $\phi_0$ edge
invariants of spherical polyhedral metrics on $(S, T)$ is the open
set in $\bold R^E$ bounded by the hypersurfaces $\Cal W(i=j+k)$,
$\Cal W(i+j+k)$ and the following set of linear inequalities: for
any disjoint sets $I, J \subset E$ so that no triangle $\sigma \in
T^{(2)}$ has exactly three edges in $J$, or exactly two edges in
$I \cup J$,
$$ \sum_{ e \in I} z(e) -\sum_{e \in J} z(e) > \pi |F(I)| -\pi |G(
I, J)| + \pi(|I| - |J|), \tag 9.9$$ where $F(I)$ consists of all
triangles with all three edges in $I$ and $G(I, J)$ consists of
all triangles with two edges in $J$ and one edge in $I$.

\rm
\medskip

Note that if all triangles are of type (a1), then the spherical
polyhedral metric shrinks to a point. In this case, the image
under $\Phi_{\lambda}$ of these degenerated metrics lies in the
hypersurface $\Phi_{\lambda}(P_{E^2}(S, T))$, i.e.,
$\Phi_{\lambda}(P_{E^2}(S, T))$ is one of the hypersurfaces
appeared in the bubble degenerated spherical metrics.

Theorem 9.8(b) generalizes the main results appeared in [Lu1].

\medskip

\noindent 9.8. {\it A proof of theorem 9.8(b)}

\medskip
By theorem 6.2(a), the map $\Phi_0: P_{S^2}(S, T) \to \bold R^E$
is a smooth embedding. To prove theorem 9.8(b), we need to
analysis the boundary of the open set  $\Omega =\Phi_0(P_{S^2}(S,
T))$ in $\bold R^E$. To this end, take a sequence $\{l^{(m)}\}$ in
$P_{S^2}(S, T)$ converging to $p \in \partial P_{S^2}(S, T)$ so
that the angles of each corner in  metrics $l^{(m)}$ in $T$ are
convergent and $\Phi_0(l^{(m)})$ converge to a point $w \in
\partial \Omega$. If all edge lengths in the degenerated metric $p$
are in the open interval $(0, \pi)$, then all degenerated
triangles in the metric $p$ are of types (a5) or (a6) due to the
classification in subsection 9.7. Thus by definition $b \in \Cal
W(i=j+k)$ or $b \in \Cal W(i+j+k)$ for some $e_i, e_j, e_k$
forming edges of a triangle in $T$. Now if some edge lengths in
the metric $p$ are $0$ or $\pi$, let
$$ I =\{ e \in E | p(e) =0\}$$ and
$$ J =\{ e \in E | p(e) = \pi\}.$$
We have $ I \cap J = \emptyset$, $I \cup J \neq \emptyset$ and
there are no triangle $\sigma \in T^{(2)}$ with all edges in $J$,
or exactly two edges in $I \cup J$. We claim that (9.9) becomes an
equality for this choice of $I, J$. Furthermore, we shall prove
that (9.9) holds for all metrics in $P_{S^2}(S, T)$.

In the discussion below, unless mentioned otherwise, all edge
lengths and angles are measured in the degenerated metric $p$.
Consider a triangle $\sigma$ with an edge in $I \cup J$.  Let
$\theta_i, \theta_j, \theta_k$ be the inner angles and $e_i, e_j,
e_k$ be the edges in $\sigma$ so that $\theta_r$ faces $e_r$.
There are four possibilities: (I) all edges of $\sigma$ are in
$I$; (II) one edge of $\sigma$ is in $I$ and the other two edges
are in $J$;  (III) one edge of $\sigma$ is in $I$ and the other
two are not in $I \cup J$; and (IV) one edge of $\sigma$ is in $J$
and the other two are not in $I \cup J$. We will analysis the
angles $\theta_r$ in each of these four cases.

\noindent Case I, all $e_r$'s are in $I$, thus the triangle
$\sigma$ is of type (a1). We obtain
$$ \theta_i + \theta_j + \theta_k = \pi. \tag 9.10$$
Note that the left-hand-side of (9.10) is strictly greater than
$\pi$ for non-degenerated spherical triangles.

\noindent Case II, $e_i \in I$ and $e_j, e_k \in J$. Thus the
triangle $\sigma$  is of type (a2). Then by the classification,
$$ \theta_i -\theta_j -\theta_k =-\pi. \tag 9.11$$
Note that the left-hand-side of (9.11) is strictly greater than
$-\pi$ for non-degenerated spherical triangles.

\noindent Case III, $e_i \in I$ and $e_j, e_k \notin I \cup J. $
Then $\sigma$ is of type (a3) so that,
$$ \theta_i =0. \tag 9.12$$
Note that the left-hand-side of (9.12) is strictly greater than
$0$ for non-degenerated triangles.

\noindent Case IV, $e_i \in J$ and $e_j, e_k \notin I \cup J$.
Then $\sigma$ is of type (a4) and
$$ -\theta_i =-\pi. \tag 9.13$$
Note that the left-hand-side of (9.13) is strictly greater than
$-\pi$ for non-degenerated triangles.

Now the left-hand-side of (9.9) can be expressed as
$$\sum_{ e \in I} \phi_0(e) -\sum_{e \in J} \phi_0(e)
=\sum_{e \in I}( \alpha + \beta) -\sum_{e \in J} (\alpha + \beta)
+ \pi(|I|-|J|) \tag 9.14$$ where $\alpha, \beta$ are angles facing
the edge $e$.

Break the first two summations in the right-hand-side of  (9.14)
into groups according to the triangles of types I, II, III, IV.
Then,
$$\sum_{e \in I}( \alpha + \beta) -\sum_{e \in J} (\alpha +
\beta)$$
$$=\sum_{\sigma \in \text{ case I} } (\theta_i + \theta_j + \theta_k) +
\sum_{\sigma \in \text{ case II} } (\theta_i - \theta_j -\theta_k)
+ \sum_{\sigma \in \text{ case III} } \theta_i +\sum_{\sigma \in
\text{ case IV} } -\theta_i \tag 9.15
$$
By equalities (9.10)-(9.13), the expression (9.15) is
 $\pi |F(I)|
-\pi |G(I \cup J)|$. This verifies that the condition (9.9)
becomes equality for the degenerated metric $p$. On the other
hand, for a non-degenerated spherical triangle, the
left-hand-sides of (9.10)-(9.13) become strictly greater than the
right-hand-side. Thus the above argument shows (9.15) is strictly
greater than $\pi |F(I)| -\pi |G(I \cup J)|$ for any metric in
$P_{S^2}(S, T)$, i.e., (9.9) holds for non-degenerated metrics.

This ends the proof of theorem 9.8 (b).

\medskip
\noindent 9.9. {\it  A proof of theorem 9.8(a)}
\medskip

By theorem 6.2(a), the map $\Phi_{\lambda}: P_{S^2}(S, T) \to
\bold R^E$ is a smooth embedding. To prove theorem 9.8(a), we need
to show that boundary points of $\Phi_{\lambda}(P_{S^2}(S, T))$ in
$\bold R^E$ come from the images of the bubbled metrics under
$\Phi_{\lambda}$. To this end, take a sequence of points $\{
l^{(m)}\}$ in $P_{S^2}(S, T)$ converging to a point $p \in
\partial P_{S^2}(S, T)$ so that inner angles of each corner in  metrics $l^{(m)}$ in $T$
 converge and $\Phi_{\lambda}(l^{(m)})$
converge to $w \in \bold R^E$. The goal is to show that $p$ must
be of bubbled degeneration, i.e., all triangles in the metric $p$
are of types (a1) or (a2).

The strategy of the proof is as follows. Since $\lambda \leq -1$,
both integrals $\int_{\pi/2}^0 \sin^{\lambda}(t) dt$ and
$\int^{\pi}_{\pi/2} \sin^{\lambda}(t) dt$ diverge. It follows
that, in the metric $p$,  if $\alpha$ and $\beta$ are two angles
facing an edge so that $\alpha \in \{0, \pi\}$, then $\beta =\pi
-\alpha \in \{0, \pi\}$. Indeed, this is due to the assumption
that $\int^{\alpha}_{\pi/2} \sin^{\lambda}(t) dt +
\int^{\beta}_{\pi/2} \sin^{\lambda}(t) dt \in \bold R$. We call
this a $\{0, \pi\}$ opposite angle propagation rule. Now if the
metric $p$ is not of bubble type, there is an inner angle $\alpha
\in \{0, \pi\}$. Using the propagation rule and  the
classification of degenerated triangles, we produce a $\{0, 0,
\pi\}$ angled edge cycle in the metric $p$. By lemma 9.2, this is
impossible.

Here is a detailed proof. First, it suffices to show there is no
degenerated triangle of types (a3)-(a6) in the metric $p$. Indeed,
if all triangles in the degenerated metric $p$ are either
non-degenerated or of types (a1) or (a2), then all triangles in
$p$ are of types (a1) or (a2). This is due to the fact that any
triangle adjacent to a triangle of types (a1) or (a2) along an
edge must be degenerated by the definition of types (a1) or (a2)
(since all edge lengths of types (a1) and (a2) are $0$ or $\pi$).
By assumption, these adjacent triangles must be of types (a1) or
(a2). Since the surface is connected and the metric $p$ is
degenerated, it follows all triangles in $p$ are of types (a1) or
(a2).

Next we prove that there is no degenerated triangles of types
(a3)-(a6) by contradiction. Suppose otherwise that there exists a
triangle $\sigma$ in the metric $p$ which is of type (a3), or
(a4), or (a5), or (a6). By the classification in subsection 9.7,
the triangle must have an inner angle $0$ or $\pi$. For
simplicity, if $\alpha, a \in [0, \pi]$, by an \it $[\alpha, a]$
triangle \rm we mean a degenerated spherical triangle with an
inner angle $\alpha$ so that the length of the opposite edge is
$a$. Thus there is a $[0,a]$ or $[\pi,a]$ triangle in the metric
$p$.
 Let $\theta_1 =0, $ or
$\pi$ be the inner angle of the triangle facing an edge $e$ and
$\theta_1'$ be the other angle facing $e$. Then by the discussion
above, $\theta_1' =\pi-\theta_1 \in \{0,\pi\}$.  Therefore, there
are both $[0, a]$ and $[\pi, a]$ degenerated triangles $\sigma$ in
$p$. Let $\theta_1 \in \{0, \pi\},$ $ \theta_2, \theta_3$ be the
inner angles and $l_1=a, l_2, l_3$ be the opposite edge lengths of
the triangle $\sigma$. We will discuss three cases according to
$a=\pi$, $a=0$, or $a \in (0, \pi)$.

\medskip
\noindent Case 1, $a=\pi$. We may assume that $\sigma$ is a $[0,
\pi]$ triangle. According to the classification in subsection 9.7,
the type of a $[0, \pi]$ triangle is (a2) where the inner angles
are $(\theta_1, \theta_2, \theta_3) =(0, 0, \pi)$ and the opposite
edge lengths $(l_1, l_2, l_3) =(\pi, 0, \pi)$. Let $\tau$ be the
triangle adjacent to the $l_3$-th edge of $\sigma$ and let $\beta$
be the angle in $\tau$ facing the $l_3$-th edge. By the discussion
above, $\beta=\pi -\theta_3=0$. Thus $\tau$ is a $[0, \pi]$
triangle. In summary, we see that $[0, \pi]$ triangles propagate
through one of its edges. In particular, there exists an edge
cycle so that each triangle in the cycle is a $[0,\pi]$ triangle.
 Since the inner angles of a $[0, \pi]$
triangle are $0, 0, \pi$, this edge cycle is also a $\{0, 0,
\pi\}$-angled edge cycle. According to lemma 9.2, this is
impossible. As a consequence, there are no $[\pi,\pi]$ triangles
in the metric $p$.

\medskip
\noindent Case 2, $a=0$. We may assume that $\sigma$ is a $[\pi,
0]$ triangle where $\theta_1 =\pi$ and $l_1=0$. According to the
classification in subsection 9.7, the type of $\sigma$ must be
either (a1) or (a2).
 If $\sigma$ is of type (a2), then its
angles and lengths are:  $(\theta_1, \theta_2, \theta_3)=(\pi,
\pi,\pi)$ and $(l_1, l_2, l_3) =(0,\pi, \pi)$. Thus $\sigma$ is
 a $
[\pi,\pi]$ triangle.  This is impossible by case 1.
 Thus $\sigma$ must be of type
(a1) which has angles  $(\theta_1, \theta_2, \theta_3)=(\pi, 0,
0)$ and lengths $(l_1, l_2, l_3)=(0,0,0)$. Let $\tau$ be the
triangle adjacent to $\sigma$ along the $l_2$-th edge. Then due to
$\theta_2=0$, the angle in $\tau$ facing $l_2$-th edge is $\pi$.
It follows that $\tau$ is again a $[\pi,0]$ triangle. Thus we see
that a $[\pi,0]$ triangle propagates through one of its edges. By
the analysis above, each $[\pi, 0]$ triangle of type (a1) has
inner angles $0, 0, \pi$. By the propagation rule, we obtain a
$\{0, 0, \pi \}$-angled edge cycle. This contradicts lemma 9.2.

\medskip
\noindent Case 3, $a \in (0, \pi)$. We may assume that $\sigma$ is
a $[0, a]$ triangle in the metric $p$. According to the
classification of degenerated triangles, the triangle $\sigma$
must be of types (a3) or (a4) or (a5).

If $\sigma$ is of type (a3), then $(\theta_1, \theta_2,
\theta_3)=(0, \pi, 0)$ and $(l_1, l_2, l_3) =(a,a, 0)$. Thus
$\sigma$ is also a $[0, 0]$ triangle.  According to case 2, this
cannot occur.

If $\sigma$ is of type (a4), then $(\theta_1, \theta_2, \theta_3)
=(0, 0, \pi)$ and $(l_1, l_2, l_3) =(a, \pi-a, \pi)$. This implies
that $\sigma$ is a $[\pi,\pi]$ triangle which is impossible by
case 1.

Thus the type of  $\sigma$ must be (a5) so that $(\theta_1,
\theta_2, \theta_3)=(0, 0, \pi)$ and $(l_1, l_2, l_3) =(a, a',
a+a') \in (0, \pi)^3$. Let $\tau$ be a triangle adjacent to
$\sigma$ along the $l_3$-th edge. Due to $\theta_3=\pi$, the inner
angle of $\tau$ facing the $l_3$-th edge must be 0. Thus $\tau$ is
a $[0, a'']$ triangle where $a'' \in (0, \pi)$. Thus a type (a5)
[0, a] triangle propagates through one of its edge. Since by the
discussion above, a type (a5), [0, a] triangle has inner angles
$0,0, \pi$, it follows that there exists a $\{0, 0, \pi\}$-angled
edge cycle in the metric $p$. This contradicts lemma 9.2.

\medskip
\noindent \S 10. {\bf  Applications Teichm\"uller Spaces and Some
Open Problems }

\medskip
We discussion some applications and open problems in this section.

\medskip
\noindent 10.1. {\it The space of all geometric triangulations
with prescribed curvature}

\medskip
The most interesting problem is probably problem 5.1 in \S5. One
supporting evidence comes from Teichm\"uller spaces on surfaces
with boundary so that the boundary lengths are prescribed. This
was discussed in subsection 5.3. The problem for $S=\bold S^2$ was
first investigated by S. S. Cairns in [Ca]. It was also related to
the work of E. Steinitz [SR] on the moduli space of all convex
polytopes in the 3-space of the same combinatorial type. Let $\Pi:
P_{K^2}(S, T) \to \bold R^V$ be the discrete curvature map sending
a metric to its discrete curvature $k$. Cairns was trying to show
that for spherical polyhedral metrics on $(\bold S^2, T)$,
$\Pi^{-1}(0)$ is either homeomorphic to a Euclidean space or is
the empty set. His first proof in 1941 contained a gap and later
in [Ca] he proved that the set $\Pi^{-1}(0)$ is connected. The
question whether $\Pi^{-1}(0)$ is a cell for spherical polyhedral
metrics on the 2-sphere became Cairns conjecture ([BCH]). In the
work of [BCH], E. Bloch, R. Connelly, and D. Henderson proved that
for Euclidean polyhedral metrics on a simplicially triangulated
disk, the space $\Pi^{-1}(0)$ is homeomorphic to a Euclidean space
for any simplicial triangulation.  Another evidence for the
affirmative solution comes from the work of Rivin [Ri] and Leibon
[Le]. They
 show that the space of all Delaunay
$\bold E^2$ or $\bold H^2$ polyhedral metrics with prescribed
discrete curvature is a cell.

The following result, which is a consequence of the works of Rivin
and Leibon,  implies that the
 spaces $\Pi^{-1}(p)$ are smooth manifolds in the Euclidean and
 hyperbolic cases.
\medskip

\noindent {\bf Proposition 10.1.} \it Suppose $(S, T)$ is a closed
triangulated surface. Then

(a) The curvature map $\Pi: P_{H^2}(S, T) \to \bold R^V$ is a
submersion.

(b) The curvature map  $\Pi$ defined on  $P_{E^2}(S, T) $ is a
submersion to the affine space $$\Cal A=\{ z \in \bold R^V |
\sum_{ v \in v} z(v) = 2 \chi(S)\}$$ of $ \bold R^V$ defined by
the Gauss-Bonnet identity .

\rm

\bigskip

\noindent {\bf Proof. } We begin with a lemma of Rivin and Leibon
relating $\phi_0$, $\psi_0$ with the discrete curvature $k$. We
will use the following notation. If $v$ is a vertex and $e$ is an
edge having $v$ as a vertex, we denote it by $e> v$. Given $v$,
the set of elements in  $\{ e \in E | e >v\}$ will be counted with
multiplicity, i.e., if the two end points of $e$ are $v$, then $e$
will be counted twice.

\medskip
\noindent {\bf Lemma 10.2.} \it Suppose $v$ is a vertex.

(a) (Rivin) for a Euclidean polyhedral metric, $\sum_{ e > v}
\phi_0(e) = k(v)-2\pi$,

(b) (Leibon) $\sum_{ e > v} \psi_0(e) = 2\pi- k(v)$. \rm

\medskip
Let $\{e_1, ..., e_m\}$ be the list of all edges counted with
multiplicity having $v$ as a vertex so that $e_i, e_{i+1}$ are
adjacent to a triangle $\sigma_i$ ( $e_{m+1}=e_1$). Let the inner
angles of $\sigma_i$ be $a_i, b_i, c_i$ with $a_i$ facing $e_i$,
$b_i$ facing $e_{i+1}$ and $c_i$ being the angle at $v$. Then
$k(v) = 2\pi- \sum_{i=1}^m c_i$, $\phi_0(e_i) = a_i+b_{i-1} -\pi$
and $\psi_0(e_i) =\frac{1}{2}( c_i + c_{i-1} + a_{i-1} -a_i + b_i
-b_{i-1})$. By summing up $\phi_0(e_i)$'s and $\psi_0(e_i)$'s, and
use $c_i=\pi-a_i-b_i$ in part (a),  one obtains the results.

\medskip
\noindent {\bf Lemma 10.3.} \it The linear map $L: \bold R^E \to
\bold R^V$ sending a vector $ z \in \bold R^E$ to $L(z) \in \bold
R^V$ defined by $L(z)(v) =\sum_{ e > v} z(e)$ is an epimorphism.
\rm

\medskip
\noindent {\bf Proof. } For a finite set $Z$, we identify the dual
space of $\bold R^Z$ with $\bold R^Z$ using the standard basis.
Then the dual map $L^*: \bold R^V \to \bold R^E$ is defined by $
L^*( f) (e) =\sum_{ v <e} f(v)$ for $e \in E$.  It suffices to
show that $L^*$ is injective. To see this, suppose $ f \in \bold
R^V$ so that $L^*(f) =0$, i.e., $f(v) =-f(v')$ whenever $v, v'$
are end points of an edge. Then $f=0$ follows by considering a
triangle with vertices $v, v', v''$. Indeed, we have $f(v) =
-f(v') =f(v'') =-f(v)$. Thus $f(v) =0$.  qed

\medskip
Now the proof of proposition 10.1 follows from the rigidity
theorems of Rivin and Leibon (Theorem 1.1(c), (d)).

Indeed, to prove part (a) of proposition 10.1, consider the affine
map $A: \bold R^E \to \bold R^V$ so that $A(z)(v) =2\pi- \sum_{ e
>v}  z(e)$. Then lemma 10.2 shows that $\Pi= A \circ
\Psi_0$. Now by lemma 10.3, the derivative of $A$ is $-L$ which is
surjective. It follows the derivative $D(\Pi) = -L D(\Psi_0)$. By
Leibon's theorem that $D(\Psi_0)$ is onto. It follows that
$D(\Pi)$ is onto.

To prove part (b) of proposition 10.1, consider the affine map $B
: \bold R^E \to \bold R^V$ defined by $B(z)(v) =2\pi+\sum_{ e
> v} z(e)$. Then lemma 10.2 shows that $\Pi = B \circ
\Phi_0$. By Rivin's rigidity theorem, the rank of $D(\Phi_0)$ is
$|E|-1$. By lemma 10.3, it follows that the rank of $D(\Pi)$ is at
least $|V|-1$. But on the other hand, by the Gauss-Bonnet formula,
$\Pi(P_{E^2}(S, T))$ lies in the affine space $\{ z \in \bold R^V
| \sum_{ v \in V} z(v) = 2\pi \chi(S)  \}$. Thus, the rank of
$D(\Pi)$ is $|V|-1$ and $\Pi$ is a submersion to the affine space.
qed.

\medskip

The special case of problem 5.1 addresses the space $\Pi^{-1}(0)$,
i.e., the space of all geometric triangulations of constant
curvature metrics on a surface. There exists the obvious map
$\phi: \Pi^{-1}(0) \to Teich(S)$ from $\Pi^{-1}(0)$ to the
Teichm\"uller space by forgetting the triangulation. In view of
the works of Cairns and Bloch-Connelly-Henderson, it seems likely
that $\phi$ is a smooth surjective submersion so that its fibers
$\phi^{-1}(p)$ are diffeomorphic to cells. The fiber
$\phi^{-1}(p)$ can be interpreted as the space of all geodesic
triangulations of combinatorial type $T$ in a fixed constant
curvature metric.

There are two related questions on $\phi$. Namely, when is
$\Pi^{-1}(0)$ non-empty and when is $\phi$ surjective? Both of
these questions have been solved by a combination of the works of
various authors. Call a triangulation  \it geometric \rm if there
exists a constant curvature metric on the surface so that each
cell in the triangulation is geodesic, i.e., the triangulation is
isotopic to a geodesic triangulation in some constant curvature
metric. The question whether $\Pi^{-1}(0)$ is non-empty is the
same as asking if the triangulation is geometric. This was solved
in the work of Thurston [Th], Colin de Verdiere [CV1], [CV3],
Koebe [Koe], Marden-Ridin [MaR] and others. See proposition 10.4
below. The surjectivity of $\phi$ has been investigated by Colin
de Verdiere [CV3]. In [CV3], Colin de Verdiere proved, among other
things, that if $T$ is a simplicial triangulation, then $\phi:
\Pi^{-1}(0) \to Teich(S)$ is onto. Moreover, a careful examination
of the method of the proof in [CV3] shows that if $\Pi^{-1}(0)$ is
non-empty, then $\phi$ is onto.

Recall that a triangulation of a space is called \it simplicial
\rm if the triangulation is isomorphic to a simplicial complex.
 We summarize the above discussion into the following.

\medskip
\noindent {\bf Proposition 10.4.} \it Suppose $T$ is a
triangulation of a closed surface $S$.

(a) ([Th], [CV1], [Ko], [MaR]) The  triangulation $T$ is a
geometric triangulation in some constant curvature metric on $S$
if and only if the lift of the triangulation to the universal
cover is simplicial.

(b) ([CV3]) If $T$ is a geometric triangulation in some constant
curvature metric, then $T$ is isotopic to a geodesic triangulation
in any constant curvature metric.

(c) The lifting of $T$ to the universal cover is simplicial if and
only if there are no null homotopic loop in the surface consisting
of at most two edges. \rm

\medskip
As mentioned above, the proof of this proposition is spread out in
various literatures.  For part (a), if the surface is of
non-positive Euler characteristic,  it is in [Th] and [CV1] where
the metric is given by the circle packing metric.
 For $\chi(S)
>0$, it is Koebe-Andreev-Thurston's theorem ([Koe], [An], [MaR]) on
circle packing. See [MaR] for a nice proof of it. The proof of
part (b) is implicit in [CV3].
  Part (c)
is a simple exercise.

One may add an additional equivalence relation to proposition
10.4(a) and (c) for surfaces of non-positive Euler characteristic.
A triangulation $T$ of a closed surface $S$ is said to support an
\it angle structure \rm if one can assign each corner of the
triangulation a positive number, called angle, so that the sum of
the angles at each vertex is $2\pi$, and each triangle with these
angle assignments becomes a $K^2$ geometric triangle where
$K^2=\bold H^2$ if $\chi(S) <0$, $K^2=\bold E^2$ if $\chi(S)=0$
and $K^2=\bold S^2$ if $\chi(S)
>0$. It can be shown ([CV1], [CL])  that for closed triangulated surfaces
of non-positive Euler characteristic, the existence of an angle
structure is equivalent to that the triangulation is geometric.
However, R. Stong [St] has constructed a non-geometric
triangulation of the 2-sphere which supports an angle structure.
See also the related work of [Gu2].

\medskip
\noindent 10.2. {\it Cellular decompositions of the Teichm\"uller
spaces}

\medskip

One interesting consequence of theorems 7.1 and 7.2 concerns the
cell decompositions of the Teichm\"uller space, first observed in
[Mo] for $\psi_0$-edge invariant.

Recall that the \it arc-complex \rm of a compact surface $S$ with
boundary is the following simplical complex, denoted by $A(S)$.
The vertices of $A(S)$ are isotopy classes $[a]$ of proper arcs
$a$ in $S$ which are homotopically non-trivial relative to the
boundary of $S$. A simplex in $A(S)$ is a collection of distinct
vertices $[a_1],..., [a_k]$ so that $ a_i \cap a_j =\emptyset$ for
all $i \neq j$. For instance the isotopy class of an ideal
triangulation corresponds to a simplex of maximal dimension in
$A(S)$. The non-fillable subcomplex $A_{\infty}(S)$ of $A(S)$
consists of those simplexes $([a_1], ..., [a_k])$ with $a_i \cap
a_j =\emptyset$ so that one component of $S-\cup_{i=1}^k a_i$ is
not simply connected. The simplexes in $A(S) -A_{\infty}(S)$ are
called fillable. Let $(|A(S)|-|A_{\infty}(S)|) \times \bold
R_{>0}$ be the geometric realization space whose points are of the
form $x =\sum_{i=1}^k c_i [a_i]$ where $c_i > 0$  so that $([a_1],
..., [a_k])$ is a fillable simplex.  Now take a point $x
=\sum_{i=1}^k c_i [a_i]$ in $(|A(S)|-|A_{\infty}(S)|) \times \bold
R_{>0}$. Let $([a_1], ..., [a_n])$ be an ideal triangulation
containing the fillable simplex $([a_1], ..., [a_k])$. Assign each
edge $[a_i]$ the positive number $z_i=c_i$ if $i \leq k$ and zero
otherwise. Then this assignment $z$ satisfies the positive edge
cycle condition in theorem 7.2 in the ideal triangulation $([a_1],
..., [a_n])$. By theorem 7.2 for $\lambda \geq 0$, there exists a
hyperbolic metric on $S$ whose $\psi_{\lambda}$-coordinate in the
ideal triangulation $([a_1],...,[a_n])$ is $z$. For
$\psi_0$-coordinate, this fact has also been established by Hazel
[Ha].

On the other hand,  the following results  of Ushijima [Us] and
Kojima [Ko] (see also [BP] and [Ha]) show that,

\medskip
\noindent {\bf Theorem 10.5 (Ushijima [Us], [Ko])}. \it For a
compact hyperbolic surface $S$ with totally geodesic boundary,
there is an ideal triangulation so that the $\psi_0$-coordinate of
the metric in the ideal triangulation is non-negative.
Furthermore, the set of all edges in the ideal triangulation with
positive $\psi_0$-coordinate form a fillable simplex in $A(S)$ and
the fillable simplex is unique. \rm

\medskip
Combining with lemma 7.5 that $\psi_0(e)
>0$ (or $\psi_0(e)=0$) if and only if $\psi_{\lambda}(e) >0$ (or
$\psi_{\lambda}(e)=0$), we can replace positivity of the
$\psi_0$-coordinate in Ushijima's theorem by $\psi_{\lambda}$.
 As a consequence, one can define an injective map
$$\Pi_{\lambda}: Teich(S) \to  |A(S) -A_{\infty}(S)| \times \bold
R_{>0}$$ by sending the equivalence class of a  hyperbolic metric
to the point $\sum_{i=1}^n z_i [a_i]$ where $(a_1, ..., a_n)$ is
the ideal triangulation produced in Ushijima's theorem and $z_i$
is the  $\psi_{\lambda}$ coordinate of the metric at the i-th edge
in the ideal triangulation. The discussion above shows that
$\Pi_{\lambda}$ is onto. Thus we obtain,

\medskip
\noindent {\bf Corollary 10.6.} \it For any compact surface with
boundary and of negative Euler characteristic and $\lambda \geq
0$, the map
$$\Pi_{\lambda}: Teich(S) \to |A(S) -A_{\infty}(S)| \times \bold R_{>0}$$
is a homeomorphism equivariant under the action of the mapping
class group. In particular, for each $\lambda$, the map
$\Pi_{\lambda}$ produces a natural cell-decomposition of the
moduli space of surfaces with boundary. \rm

\medskip
We remark that the underlying cell-structures for various
$\lambda$'s are the same. In particular, if $\lambda \neq
\lambda'$, then $\Pi_{\lambda'}^{-1} \Pi_{\lambda} $ is a
self-homeomorphism of the Teichm\"uller space preserving the
cell-structure derived from $\Pi_0$. These self-homeomorphisms of
$Teich(S)$ deserve a further study. Finally, we remark that those
coordinates $\psi_{\lambda}$ with $\lambda <0$ also produce
cellular structures on the Teichm\"uller space $Teich(S)$ due to
Guo's result [Gu1].

\medskip

\noindent 10.3. {\it Derivative cosine laws for homogeneous
spaces}

\medskip
There are similar  cosine laws for other homogeneous spaces. It
should be interesting to know if the derivatives of these cosine
law have some interesting properties and applications.

\medskip
\noindent 10.4. {\it Global rigidity of Euclidean polyhedral
metrics}
\medskip
We do not know if $\phi_{\lambda}$ determines a Euclidean
polyhedral metric when $\lambda >0$, or $\lambda \in (-1, 0)$.
From the point of view of theorem 6.2, it seems highly likely that
$\phi_{\lambda}$ determines the metric up to scaling when $\lambda
>0$.

\medskip
\noindent 10.5. {\it Global rigidity of hyperbolic or spherical
polyhedral metrics}

\medskip
It will be interesting to know for $\lambda \in (-1, 0)$, if
$\phi_{\lambda}$  determines a spherical polyhedral metric up to
equivalence, or $\psi_{\lambda}$ determines the hyperbolic metric
up to equivalence.

\medskip
\noindent 10.6. {\it Non-convex or concave energy and its use}

\medskip
Those non-convex or concave energy functions in theorems 3.2 and
3.4 have the corresponding variational principles on triangulated
surfaces. We do not know any use of these yet. Probably the most
interesting questions in this area are:

(a) Give a new proof of Andreev-Koebe-Thurston's circle packing
theorem on the 2-sphere using variational principle based on Colin
de Verdiere's energy for spherical triangles or its Legendre
transformation. See [BS] for the Legendre transformed energy
function.

(b) Is a hyperbolic polyhedral surface determined by any of the
$\phi_{\lambda}$ edge invariant?

\medskip

\noindent 10.7. {\it Convexity of the space of edge invariants}

\medskip
It will be interesting to know if the spaces
$\Phi_{\lambda}(P_{E^2}(S, T))$ are convex hypersurfaces, or if
$\Psi_{\lambda}(P_{H^2}(S, T))$ are convex sets. Numerical
calculations shows that both $K_{\lambda}(CP_{E^2}(S, T))$ and
$K_{\lambda}(CP_{H^2}(S, T))$ are not convex in general. However,
as Thurston proved in the case $\lambda =0$, there are still cases
which we do not know. Namely, find all $\lambda$ for which these
spaces are convex sets or convex codimension-1 surfaces.

\medskip
\noindent 10.8. {\it Miscellaneous remarks}

It will be interesting to know if these edge invariants
$\phi_{\lambda}, \psi_{\lambda}$ and $k_{\lambda}$ correspond to
some curvatures in Riemannian geometry as triangulations become
finer and converge to a Riemannian metric.

The relationship between the $\psi_{\lambda}$ coordinate with
quantum Teichm\"uller space deserves a further study. See [CF],
[Ka], [BL], [Te] and others for more details. With the apparent
resolution of the geometrization conjecture for 3-manifolds, it
seems the main focus of 3-manifold study will be shifted to
hyperbolic 3-manifolds. According to a conjecture of Thurston,
each closed hyperbolic 3-manifold has a finite cover which is a
surface bundle over the circle. Thus, put all these together, it
seems that up to finite cover, topology and geometry of
3-manifolds are governed by homeomorphisms of surfaces. From this
point of view, quantum Teichm\"uller theory will likely to play an
important role in 2+1 TQFT. We hope that derivative cosine law
will be a part of the grand picture.

\medskip

\medskip
\medskip
\centerline{\bf Appendix A. A Proof of uniqueness of the 1-forms}

\medskip
The goal of this appendix is to prove the uniqueness part of
theorem 1.5.

\medskip
\noindent {\bf Theorem 1.5.} \it For the cosine law function
$y=y(x)$, all closed 1-forms of the form $ w =\sum_{i=1}^3
f(y_i)dg(x_i)$ where $f,g$ are two non-constant smooth functions,
are up to scaling and complex conjugation,
$$\omega_{\lambda} =\sum_{i=1}^3 \int^{y_i} \sin^{\lambda}(t) dt  d( \int^{x_i} \sin^{-\lambda -1}(t)
dt)=\sum_{i=1}^3 \frac{ \int^{y_i} \sin^{\lambda}(t)
dt}{\sin^{\lambda+1}(x_i)} dx_i$$
 for some $\lambda \in \bold C$.

All closed 1-forms of the form $ \sum_{i=1}^3 f(y_i) dg(r_i)$
where $f,g$ are two non-constant smooth functions, are up to
scaling and complex conjugation,
$$\eta_{\lambda} = \sum_{i=1}^3 \int^{y_i} \tan^{\lambda}(t/2) dt  d(\int^{r_i} \cos^{-\lambda-1}(t) dt)=
\sum_{i=1}^3 \frac{ \int^{y_i} \tan^{\lambda}(t/2)
dt}{\cos^{\lambda+1}(r_i)} dr_i.$$ for some $\lambda \in \bold C$.
In particular, all closed 1-forms are  holomorphic or
anti-holomorphic. \rm

\medskip
\noindent {\bf Proof.} Let $\{i,j,k\}=\{1,2,3\}$. The proof of the
uniqueness depends on the following lemma.

\medskip
\noindent {\bf Lemma 2.3.} \it Suppose $y=y(x)$ is the cosine law
function and $f, g$ are two smooth non-constant functions.

(a) If $f(y_i)/g(x_i)$ is independent of the indices for all $x$,
then there are  constants $\lambda, \mu$, $c_1, c_2$ so that $f(t)
=c_1\sin^{\lambda}(t) \sin^{\mu}(\bar t)$ and $g(t)
=c_2\sin^{\lambda}(t)\sin^{\mu}(\bar t)$.

(b) If $r_i = 1/2(x_j+x_k-x_i)$, and $f(y_i)/g(r_i)$ is
independent of the indices for all $r$, then there are constants
 $\lambda, \mu$, $c_1, c_2$ so that $f(t) = c_1\tan^{\lambda}(t) \tan^{\mu}(\bar t)$ and
$g(t) = c_2\cos^{\lambda}(t)\cos^{\mu}(\bar t)$.

\rm
\medskip
\noindent {\bf Proof.} We use $f_{z}$ and $f_{\bar z}$ to denote
the partial derivatives $\frac{\partial f}{\partial z}$ and
$\frac{\partial f}{\partial \bar z}$ respectively. Note that
$\partial y_i/\partial \bar x_j=0$. Take $\frac{\partial}{\partial
x_k}$ to the identity $\frac{f(y_i)}{g(x_i)} =\frac{
f(y_j)}{g(x_j)}$, we obtain
$$ \frac{ f_{z}(y_i)}{g(x_i)} \frac{\partial y_i}{\partial x_k} = \frac{ f_{ z}(y_j)}{g(x_j)} \frac{\partial y_j}{\partial x_k}$$
By the derivative cosine law that $\frac{ \partial y_i/\partial
x_k}{\partial y_j/\partial x_k}
=\frac{\sin(y_i)\cos(y_j)}{\sin(y_j)\cos(y_i)}$ and
$\frac{f(y_i)}{g(x_i)} =\frac{ f(y_j)}{g(x_j)}$, we obtain,
$$ \frac{ f_z(y_i) \sin(y_i)}{f(y_i) \cos(y_i)} = \frac{ f_z(y_j) \sin(y_j)}{f(y_j) \cos(y_j)}$$
The variables $y_i, y_j$ are independent. This shows that there is
a constant $\lambda \in \bold C$ so that
$$\frac{f_z(t)}{f(t)} = \lambda \cot(t),$$
i.e.,
$$\frac{\partial \ln (f(z))}{\partial z} =\frac{ \partial (\lambda \ln \sin(z))}{\partial z}.$$
If we take $\frac{\partial}{\partial \bar x_k}$ to the equation
$\frac{f(y_i)}{g(x_i)} =\frac{ f(y_j)}{g(x_j)}$ and use $\partial
y_i/\partial \bar x_k =0$, we obtain, by the same argument as
above,
$$\frac{ \partial (\ln f(z))}{\partial \bar z} = \frac{\partial (\mu \ln \sin(\bar z))}{\partial \bar z}$$ for some constant $\mu \in \bold C$.
This implies that $f(z) = c_1 \sin^{\lambda}(z) \sin^{\mu}(\bar
z)$.  Now substitute it back to $f(y_i)/g(x_i)$ and use the sine
law, we obtain that $\frac{ g(x_i)}{\sin^{ \lambda}(x_i)
\sin^{\mu}(\bar x_i)}$ is independent of the indices $i$. Thus it
must be a constant. This shows that $g(z) = c_2 \sin^{\lambda}(z)
\sin^{\mu}(\bar z)$ for some constant $c_2$.

The proof of the second part (b) is exactly the same as part (a)
where we use the tangent law that $\tan(y_i/2)/\cos(r_i)$ is
independent of $i$ instead of the sine law. QED.

\medskip
To prove the uniqueness part of theorem 1.5, we write the closed
1-form $w =\sum_{i=1}^3 f(y_i) dg(x_i)$ as
$$w =\sum_{i=1}^3 f(y_i) g_z(x_i) dx_i + f(y_i) g_{\bar z}(x_i) d\bar x_i.$$ The 1-form $w$ is closed if and only if
for $i \neq j$, the expressions  $ \frac{\partial (f(y_i)
g_z(x_i))}{\partial x_j}$ and$ \frac{\partial (f(y_i) g_{\bar
z}(x_i))}{\partial \bar x_j}$ are symmetric in $i,j$ and
$$ \frac{ \partial (f(y_i) g_{z}(x_i)}{\partial \bar x_j} =
\frac{\partial (f(y_j) g_{\bar z}(x_j))}{\partial x_i}.  \tag 1$$
The symmetry of $i,j$ in $ \frac{\partial (f(y_i)
g_z(x_i))}{\partial x_j}$ and theorem 2.1 show that $$f_z(y_i)
g_z(x_i) \sin(x_i)=f_z(y_j)g_z(x_j) \sin(x_j).$$ By lemma 2.3,
there are constants $c_1, c_2, \alpha, \beta$ so that
$$ f_z(t) = c_1 \sin^{\alpha}(t) \sin^{\beta}(\bar t)  \tag 2$$
and
$$ g_z(t) = c_2  \sin^{-\alpha-1}(t) \sin^{-\beta}(\bar t).  \tag 3$$
By the same argument using the symmetry of $i,j$ in
$\frac{\partial (f(y_i) g_{\bar z}(x_i))}{\partial \bar x_j}$, we
obtain
$$ f_{\bar z}(t) = c_3 \sin^{\lambda}(t) \sin^{\mu}(\bar t)  \tag 4$$
and
$$ g_{\bar z} (t) = c_4 \sin^{-\lambda}(t) \sin^{-\mu-1}(\bar t) \tag 5$$
for some constants $c_3, c_4, \lambda, \mu$. Substitute (2)-(5)
into (1), we obtain

$$ c_2c_3 \sin^{\lambda}(y_i)\sin^{\mu}(\bar y_i) \sin^{-\alpha-1}(x_i) \sin^{-\beta+1}(\bar x_i) \cos(\bar y_k)\bar B$$
$$=c_1c_4 \sin^{\alpha}(y_j) \sin^{\beta}(\bar y_j) \sin^{-\lambda+1}(x_j) \sin^{-\mu-1}(\bar x_j)\cos(y_k)B. \tag 6$$
where $\sin(x_i)=B \sin(y_i)$ and $B$ is a function symmetric in
$i,j,k$. We claim (6) implies that $c_1c_2c_3c_4=0$. Indeed,
suppose otherwise that $c_1c_2c_3c_4 \neq 0$. We will derive a
contradiction as follows. Identity (6) can be written as,
$$c_2 c_3 \sin^{\lambda -\alpha-1}(y_i) \sin^{\mu +1-\beta}(\bar y_i) \cos(\bar y_k) B^{-\alpha-1}(\bar B)^{-\beta+2}$$
$$=c_1 c_4 \sin^{-\lambda +\alpha+1}(y_j) \sin^{-\mu-1+\beta}(\bar y_j) \cos(y_k) B^{-\lambda+2} (\bar B)^{-\mu -1} $$
As a consequence, we conclude that
$$ (\sin^{\lambda -\alpha -1}(y_i) \sin^{\mu+1-\beta}(\bar y_i))(\sin^{\lambda -\alpha-1}(y_j)
\sin^{\mu+1-\beta}(\bar y_j)) \frac{\cos(\bar y_k)}{\cos(y_k)}
\tag 7$$ is independent of the indices $i,j,k$. In particular,
identity (7) is equal to
$$ (\sin^{\lambda -\alpha -1}(y_i) \sin^{\mu+1-\beta}(\bar y_i))(\sin^{\lambda -\alpha-1}(y_k) \sin^{\mu+1-\beta}(\bar y_k)) \frac{\cos(\bar y_j)}{\cos(y_j)}.$$
This shows that
 $$ (\sin^{\lambda -\alpha -1}(y_k) \sin^{\mu+1-\beta}(\bar y_k)) \frac{\cos( y_k)}{\cos(\bar y_k)}
= (\sin^{\lambda -\alpha -1}(y_j) \sin^{\mu+1-\beta}(\bar
y_j))\frac{\cos( y_j)}{\cos(\bar y_j)}$$ Since $y_j, y_k$ are
independent variables, both sides must be constant. But that is
impossible.

As a consequence, we see that $c_1c_2c_3c_4 =0$. Since we assume
that $f$ and $g$ are non-constant functions, we have $|c_1|+|c_3|
\neq 0$ and $|c_2|+|c_4| \neq 0$.   Now if $c_3=0$, then $c_1 \neq
0$ due to $|c_1|+|c_3|>0$. But  (6) shows that $c_1c_4=0$. Since
$c_1 \neq 0$, we must have $c_4=0$. This shows, by (4) and (5)
that $f_{\bar z}=g_{\bar z}=0$, i.e., $f$ and $g$ are holomorphic.
By (2) and (3), due to the holomorphic property of $f,g$, it
follows that $\beta =-0$, i.e., $f(z) = c_1 \int^z
\sin^{\alpha}(t) dt$ and $g(z)= c_2 \int^z \sin^{-\alpha -1}(t)
dt$. The same argument shows that if $c_1=0$, then $f,g$ are
anti-holomorphic given by (4) and (5) with $\lambda=0$. This
establishes theorem 1.5 for $w_{\lambda}$ family.

 The proof for the forms $w =\sum_{i=1}^3 f(y_i) d
g(r_i)$ is exactly the same by using the tangent law that
$\tan(y_i/2)/ \cos(r_i)$ is independent of the indices  and lemma
2.3(b).

\medskip

\medskip
\centerline{\bf Appendix B. Derivative Cosine Law of Second Kind}

%\medskip

%\centerline{\bf Feng Luo}

%\medskip

%\centerline{\bf July 29, 2006, Da Lain, China}

\medskip

Suppose that $y=y(x)$ is the cosine law function so that

$$ \cos(y_i) =\frac{ \cos(x_i) + \cos(x_j) \cos(x_k)}{\sin(x_j)
\sin(x_k)} \tag 1$$ where $\{i,j,k\}=\{1,2,3\}$. This convention
of $\{i,j,k\}=\{1,2,3\}$ is assumed in this appendix.

 Then we know
that

$$ \cos(x_i) =\frac{ \cos(y_i) - \cos(y_j) \cos(y_k)}{\sin(y_j)
\sin(y_k)}  \tag 2$$

\medskip
Identity (2) shows that
$$\cos(y_i) = \cos(y_j)\cos(y_k) + \sin(y_j) \sin(y_k) \cos(x_i)
\tag 3$$

We consider $y_i=y_i(y_j, y_k, x_i)$ and $x_j=x_j(y_j, y_k, x_i)$
as functions of $y_j, y_k$ and $x_i$.
  Let $A^*_{ijk}=\sin(y_i)
\sin(y_j)$$ \sin(x_k)$ and $A_{ijk} = \sin(x_i) \sin(x_j)
\sin(y_k)$. Both $A^*_{ijk}$ and $A_{ijk}$ are independent of the
indices due to the sine law.

\medskip
\noindent {\bf Derivative cosine law II.} \it The derivatives of
functions $y_i=y_i(y_j, y_k, x_i)$ and $x_j=x_j(y_j, y_k, x_i)$
satisfy,

$$\frac{\partial y_i}{\partial y_j} =\cos(x_k)  \tag 4$$

$$\frac{\partial y_i}{\partial x_i} =\frac{A^*_{ijk}}{\sin(y_i)} =\frac{A_{ijk}}{\sin(x_i)} \tag 5$$

$$\frac{\partial x_j}{\partial y_k}= -\sin(x_j) \cot(y_i)  \tag
6$$

$$\frac{\partial x_j}{\partial y_j} =\frac{ \sin(x_k)}{\sin(y_i)}
\tag 7$$

$$ \frac{\partial x_j}{ \partial x_i} = -\frac{\sin(x_j)
\cos(x_k)}{\sin(x_i)} \tag 8$$

 \rm

\medskip
\noindent {\bf Proof.} Take derivative $\partial /\partial x_i$ to
(3), we have
$$-\sin(y_i)\frac{\partial y_i}{\partial x_i}
=-\sin(y_j)\sin(y_k)\sin(x_i)$$
 Divide it by
$-\sin(y_i)$ we obtain (5).

\medskip
To see (4), take $\partial/\partial y_j$ to (3). We obtain
$$-\sin(y_i) \frac{\partial y_i}{\partial y_j} = -\sin(y_j)
\cos(y_k) +\cos(y_j) \sin(y_k) \cos(x_i).  \tag 9$$  Let $c_i
=\cos(x_i)$ and $s_i =\sin(x_i)$. By the sine law, then (9) can be
written as,

$$\frac{ \partial y_i}{\partial y_j} =\frac{\sin(y_j)}{\sin(y_i)}
\cos(y_k)  -\cos(y_j) \cos(x_i) \frac{ \sin(y_k)}{\sin(y_i)}$$
$$=\frac{s_j \cos(y_k)}{s_i} - \frac{ \cos(y_j) \cos(x_i)
s_k}{s_i}$$
$$=\frac{1}{s_i} ( s_j \frac{c_k+ c_i c_j}{s_i s_j} - s_kc_i \frac{ c_j+c_i
c_k}{s_i s_k})$$
$$=\frac{1}{s_i^2}(c_k+c_ic_j - c_i c_j - c_i^2 c_k)$$
$$=\frac{1}{s_i^2}(c_k s_i^2)$$
$$=c_k$$
This verifies (4).

To see the partial derivatives of $x_j= x_j(y_j, y_k, x_i)$, we
use the sine law $$\sin(x_j) =\sin(x_i) \sin(y_j)/\sin(y_i) \tag
10$$

Take partial derivative of (10) with respect to $y_k$. We obtain,

$$ \cos(x_j) \frac{\partial x_j}{\partial y_k}
=-\sin(x_i)\sin(y_j)\frac{ \cos(y_i)}{\sin^2(y_i)}\frac{\partial
y_i}{\partial y_k}$$ By (4), we obtain that
$$\frac{\partial x_j}{\partial y_k}
= -\frac{ \sin(x_i) \sin(y_j)\cos(y_i)}{\sin^2(y_i)}$$
$$=-\sin(x_j) \cot(y_i)$$ where the last equation is due to the sine law.
This establishes (6)

To see (7), we take the partial derivative with respect to $x_i$
of (10). It becomes,
$$\cos(x_j)\frac{\partial x_j}{\partial x_i}
=\frac{\sin(y_j)}{\sin^2(y_i)} (\cos(x_i) \sin(y_i) -\sin(x_i)
\cos(y_i) \partial y_i/\partial x_i)$$ Using identity (5) and the
sine law, the above is
$$ \frac{\cos(x_i) \sin(y_i) -\cos(y_i) A_{ijk}}{\sin^2(y_i)}
\sin(y_j)$$
$$=\frac{ \cos(x_i) \sin(y_i) -\cos(y_i) \sin(y_i) \sin(x_j)
\sin(x_k)}{\sin^2(y_i)} \sin(y_j)$$
$$=\frac{ \cos(x_i) -(\cos(x_i) + \cos(x_j)\cos(x_k))}{\sin(y_i)}
\sin(y_j)$$
$$=-\frac{ \cos(x_j)\cos(x_k)}{\sin(y_i)} \sin(y_j)$$
$$=-\frac{ \cos(x_j) \cos(x_k)}{\sin(x_i)}{\sin(x_j)}$$
Now divide both side by $\cos(x_j)$, we obtain identity (8).

Finally, take partial derivative with respect to $y_j$ to (10).
Use (4) and the sine law, we obtain,

$$\cos(x_j) \frac{\partial x_j}{\partial y_j} = \sin(x_i) \frac{
\cos(y_j) \sin(y_i) -\sin(y_j) \cos(y_i) \partial y_i/\partial
y_j}{\sin^2(y_i)}$$
$$=\frac{\sin(x_i)}{\sin^2(y_i)}( \cos(y_j)\sin(y_i) -\sin(y_j)
\cos(y_i)\cos(x_k)) \tag 11$$ Let $C_r =\cos(y_r)$ and $S_r
=\sin(y_r)$. Then by (2), equation (11) becomes
$$=\frac{\sin(x_i)}{S_i^2}(C_jS_i-S_jC_i
\frac{C_k-C_iC_j}{S_iS_j})$$
$$=\frac{\sin(x_i)}{S_i^3}(C_jS_i^2-C_iC_k+C_i^2C_j)$$
$$=\frac{\sin(x_i)S_k}{S_i^2}(\frac{C_j-C_iC_k}{S_i S_k})$$
$$=\frac{S_k \sin(x_i)\cos(x_j)}{S^2_i}$$
$$=\frac{\sin(x_k) \cos(x_j)}{\sin(y_i)}.$$
Divide both sides by $\cos(x_j)$, we obtain (7). qed

\medskip
We remark that identity (6) for Euclidean triangles was in [CL,
lemma A1(d)].

\medskip

 \noindent {\bf Corollary B2.} \it Let $\lambda
\in \bold C$.

(a)  Consider $y_j, y_k, x_i$ as variables and $x_i$ fixed. Then
the  differential 1-form
$$\frac{\int^{x_j} \sin^{\lambda}( t) dt}{\sin^{\lambda+1}(y_j)}
dy_j + \frac{\int^{x_k} \sin^{\lambda}( t)
dt}{\sin^{\lambda+1}(y_k)} dy_k$$ is closed.

(b)  ([GL]) Consider $y_j, y_k, x_i$ as variables and $x_i$ fixed.
Then the differential 1-form
$$(\int^{x_k} \sin^{\lambda}( t) dt) \sin^{\lambda-1}(y_j)
dy_j + (\int^{x_j} \sin^{\lambda}( t) dt) \sin^{\lambda-1}(y_k)
dy_k$$ is closed.

(c) ([GL]) Consider $x_i,x_j, x_k$ as variables and $x_i$ fixed.
Then the differential 1-form
$$(\int^{y_k} \sin^{\lambda}( t) dt) \sin^{\lambda+1}(x_j)
dx_j + (\int^{y_j} \sin^{\lambda}( t) dt) \sin^{\lambda+1}(x_k)
dx_k$$ is closed.

\rm

\medskip
The proof is a simple application of identities (6) and (7) in
above theorem. We omit the detail. The integrations of the 1-forms
for $\lambda=0,-1$ in part (b) for geometric triangles were first
discovered by Bobenko-Springborn [BS].  Bobenko-Springborn  showed
the integration of the 1-form for $\lambda=0$ can be identified
with the dilogarithmic function.  In the work of [GL], a further
study of the applications of the derivative cosine law of second
kind are carried out.

\medskip
\medskip
\centerline{\bf Appendix C. Relationship to the Lobachevsky
Function}

\medskip

In the special cases of $\lambda =\pm 1$ or $0$, some of
integrations $\int^u w_{\lambda}$ and $\int^u \eta_{\lambda}$ in
theorem 1.5 and corollary B2 or their Legendre transformations
have been found explicitly by various authors. We give a brief
summary in this appendix.

Following Milnor [Mi], let $\Lambda(z) = \int^z_0 -\ln(2 \sin(t))
dt$ be the (complex valued) Lobachevsky function defined as a
multi-valued complex analytic function (depending on the choice of
the branch of $\ln(t)$ and the path). This function is related to
the dilogarithm function (see [Mi]).

\medskip
Let $y=y(x)$ be the cosine function defined by (1.8), $x_i = r_j +
r_k$ for $\{i,j,k\}=\{1,2,3\}$ and $r=(r_1, r_2, r_3)$. Then we
have,

\medskip
\noindent {\bf Proposition C1.} \it The following identities hold
up to addition of a constant.

(a) ([Lu1]) $$\int^{x} \sum_{i=1}^3 \ln \tan(y_i/2) dx_i =-
\sum_{i=1}^3 \Lambda(\pi/2 - r_i) + \Lambda(\pi/2 -r_1 - r_2 -
r_3) + \sqrt{-1} \pi (\sum_{i=1}^3 r_i ) \tag 1$$

(b) (Leibon [Le]) $$\int^r 2\ln \sin(y_i/2) dr_i = \sum_{i=1}^3[
\Lambda(\pi/2 -r_i) +\Lambda(r_i + r_{i+1}) +\sqrt{-1} \pi r_i] +
\Lambda(\pi/2 - r_1 -r_2-r_3)  \tag 2$$ where $r_4=r_1$.

(c) (Bobenko-Springborn [BS]) Consider $x_1, x_2, y_3$ as
variables and fixing $y_3$   The integration
$$ \int^{(x_1, x_2)} \ln \tan(y_1/2) dx_2+ \ln \tan(y_2/2) dx_1 \tag 3$$
$$= \Lambda(\pi/2 -r_1) + \Lambda(\pi/2-r_2) -\Lambda(\pi/2-r_3)
+\Lambda(\pi/2 -r_1 -r_2 -r_3)+ \sqrt{-1}\pi/2(x_1+x_2) + c,$$
where the constant $c$ depends only on $y_3$.

 \rm
\medskip
\noindent {\bf Proof.} The proof is straight forward by checking
the derivatives of the both sides. In part (a), the partial
derivative with respect to $x_i$ of the left-hand-side is $\ln(
\tan(y_i/2))$ by definition. By the tangent law (2.10), we have
$$ 2 \ln(\tan(y_i/2)) = \ln \cos(r_i) + \ln \cos(r_1+r_2+r_3) -\ln
\cos(r_j) -\ln \cos(r_k) +\sqrt{-1} \pi.$$ The right-hand-side of
the above equation is the $x_i$-th partial derivative of the
right-hand-side of (1) by the definition of the Lobachevsky
function.

In part (b), we use the following identity that
$$ \sin^2(y_i/2) = \frac{1-\cos(y_i)}{2} = \frac{ \sin(x_j)
\sin(x_k) -\cos(x_i) -\cos(x_j) \cos(x_k)}{2\sin(x_j) \sin(x_k)}$$
$$=-\frac{\cos(r_i)\cos(r_1+r_2+r_3)}{ \sin(r_i+r_k)
\sin(r_i+r_j)}$$ Now take the logarithm of this function and
compare with the partial derivatives of the right-hand-side of
(2).

The proof of (3) is the same as above. We omit the details. qed
\medskip

These integrations in the cases of spherical or hyperbolic
triangles have geometric interpretations. To be more precise, for
$x$, $y$ to be the inner angles and edge lengths of a spherical
triangle, the integration in proposition C1(a) is the volume of
the ideal hyperbolic octahedron which is the convex hull of the
six intersection points of the three circles at the sphere at
infinity forming a  triangle of inner angles $x_1, x_2, x_3$ (see
[Lu1]). If $x, y$ are the inner angles and edge lengths of a
hyperbolic triangle, Leibon [Le] showed that integration in
proposition C1(b) is the volume of the ideal prism which is the
convex hull of the six intersection points at the sphere at
infinity of the three circles forming a triangle of inner angles
$x_1, x_2, x_3$. For spherical triangle of inner angles $x_1, x_2,
x_3$, the integral in proposition C1(b) was shown by P. Doyle [Le]
to be the volume of the hyperbolic tetrahedron with exactly three
vertices at infinite and a finite vertex $v$ so that the dihedral
angles at the edges from $v$ are $x_1, x_2, x_3$.

\medskip

\centerline{\bf References}

\medskip

[An] Andreev, E. M., Convex polyhedra in Lobachevsky spaces.
(Russian) Mat. Sb. (N.S.) 81 (123) 1970 445--478.

[BCH] Bloch, Ethan D.; Connelly, Robert; Henderson, David W., The
space of simplexwise linear homeomorphisms of a convex $2$-disk.
Topology 23 (1984), no. 2, 161--175.

[BL] Bonahon, F. \& Liu, X., Representations of the quantum
Teichm\"uller space, and invariants of surface diffeomorphisms,
arXiv:math.GT/0407086.

[BE] Bowditch, B. H.; Epstein, D. B. A., Natural triangulations
associated to a surface. Topology 27 (1988), no. 1, 91--117.

[Br] Br\"agger, W.,  Kreispackungen und Triangulierungen. Enseign.
Math., 38:200-217,1992.

[Bu]  Buser, Peter, Geometry and spectra of compact Riemann
surfaces. Progress in Mathematics, 106. Birkhauser Boston, Inc.,
Boston, MA, 1992.

[BS] Bobenko, Alexander I.; Springborn, Boris A., Variational
principles for circle patterns and Koebe's theorem. Trans. Amer.
Math. Soc. 356 (2004), no. 2, 659--689.

[Ca] Cairns, Stewart S., Isotopic deformations of geodesic
complexes on the 2-sphere and on the plane. Ann. of Math. (2) 45,
(1944). 207--217.

[CKP]  Cohn, Henry; Kenyon, Richard; Propp, James, A variational
principle for domino tiling. J. Amer. Math. Soc. 14 (2001), no. 2,
297--346.

[CL] Chow, Bennett; Luo, Feng, Combinatorial Ricci flows on
surfaces. J. Differential Geom. 63 (2003), no. 1, 97--129.

[CV1] Colin de Verdiere, Yves, Un principe variationnel pour les
empilements de cercles. Invent. Math. 104 (1991), no. 3, 655--669.

[CV2] Colin de Verdiere, Yves, private communication.

[CV3] Colin de Verdiere, Yves, Comment rendre g\'eod\'esique une
triangulation d'une surface?  Enseign. Math. (2) 37 (1991), no.
3-4, 201--212.

[FC] Vladimir V. Fok, Leonid O. Chekhov, Quantum Teichm¨uller
spaces (Russian) Teoret. Mat. Fiz. 120 (1999), 511–528;
translation in Theoret. and Math. Phys. 120 (1999), 1245–1259

[Ga] Garrett, Brett, Circle packings and polyhedral surfaces.
Discrete Comput. Geom. 8 (1992), no. 4, 429--440.

[Gu1] Guo, Ren, On parameterizations of Teichm\"uller spaces of
surfaces with boundary, arXiv: math.GT/0612221.

[Gu2] Guo, Ren,  Geometric angle structures on triangulated
surfaces, \newline arXiv:math.GT/0601486.

[GL] Guo, Ren; Luo, Feng, Applications of the derivative of the
cosine law to polyhedral surfaces, in preparation.

[Ha] Hazel, Graham, Triangulating TeichmÄuller space using the
Ricci flow, Ph. D. thesis, UC San Diego, 2004,
http://math.ucsd.edu/~thesis/thesis/ghazel/ghazel.pdf

[IT]  Imayoshi, Y.; Taniguchi, M., An introduction to
Teichm\"uller spaces. Translated and revised from the Japanese by
the authors. Springer-Verlag, Tokyo, 1992.

[Ka] Rinat Kashaev, Quantization of Teichm\"uller spaces and the
quantum dilogarithm, Lett. Math. Phys. 43 (1998), 105–115.

[Ko] Kojima, Sadayoshi, Polyhedral decomposition of hyperbolic
$3$-manifolds with totally geodesic boundary. Aspects of
low-dimensional manifolds, 93--112, Adv. Stud. Pure Math., 20,
Kinokuniya, Tokyo, 1992.

[Le] Leibon, Gregory, Characterizing the Delaunay decompositions
of compact hyperbolic surfaces. Geom. Topol. 6 (2002), 361--391.

[Lu1] Luo, Feng,  A characterization of spherical polyhedron
surfaces,  J. Differential Geom. 74, (2006), no. 3,
407-424.

[Lu2] Luo, Feng, On Teichm\"uller spaces of Surfaces with
boundary, to appear in Duke Math. Jour., arXiv:/math.GT/0601364

[Lu3] Luo, Feng,  Volume and angle structures on 3-manifolds, to
appear in Asia Jour. of Math., arXiv:math.GT/0504049.

[Le] Leibon, Gregory, Characterizing the Delaunay decompositions
of compact hyperbolic surfaces. Geom. Topol. 6 (2002), 361--391

[MaR] Marden, Al; Rodin, Burt, On Thurston's formulation and proof
of Andreev's theorem. Computational methods and function theory
(Valparalo, 1989), 103--115, Lecture Notes in Math., 1435,
Springer, Berlin, 1990.

[Mi] Milnor, John, computation of volume, chapter 7 of Thurston's
note on geometry and topology of 3-manifolds, at
www.msri.org/publications/books/gt3m/.

[MiR] Minda, David; Rodin, Burt, Circle packing and Riemann
surfaces. J. Anal. Math. 57 (1991), 221--249.

[Pe] Penner, R. C., Decorated Teichm\"uller theory of bordered
surfaces.  Comm. Anal. Geom. 12 (2004), no. 4, 793--820.

[Ri] Rivin, Igor, Euclidean structures on simplicial surfaces and
hyperbolic volume. Ann. of Math. (2) 139 (1994), no. 3, 553--580.

[Sc] Schlenker, Jean-Marc, Circle patterns on singular surfaces,
arXiv:math.DG/0601531

[Sp] Springborn, B., A variational principle for weighted Delaunay
triangulations and hyperideal polyhedra, arXiv: math.GT/0603097

[SR] Steinitz, Ernst; Rademacher, Hans Vorlesungen \"uber die
Theorie der Polyeder unter Einschluss der Elemente der Topologie.
Reprint der 1934 Auflage. Grundlehren der Mathematischen
Wissenschaften, No. 41. Springer-Verlag, Berlin-New York, 1976.

[St] Stong, Richard, private communication.

[Te] Teschner, J., An analog of a modular functor from quantized
Teichm"uller theory, \newline arXiv:math.QA/0510174

[Th] Thurston, William, Geometry and topology of 3-manifolds,
lecture notes, Math Dept., Princeton University,   1978,
 at www.msri.org/publications/books/gt3m/

\medskip
\noindent

Department of Mathematics

Rutgers University

Piscataway, NJ 08854, USA

email: fluo\@math.rutgers.edu

\end
\newpage

\newpage

%\midspace{0.1cm}
\centerline{\epsfbox{1.1.eps}} \centerline{Figure 1.1}

\end

\end

\end

\end

\end